# SMARANDACHE NEUTROSOPHIC ALGEBRAIC STRUCTURES

**W. B. Vasantha Kandasamy**
e-mail: vasantha@iitm.ac.in
web: http://mat.iitm.ac.in/~wbv

**2006**

# SMARANDACHE NEUTROSOPHIC ALGEBRAIC STRUCTURES

**W. B. Vasantha Kandasamy**

# CONTENTS

# PREFACE

In this book for the first time we introduce the notion of Smarandache neutrosophic algebraic structures. Smarandache algebraic structures had been introduced in a series of 10 books. The study of Smarandache algebraic structures has caused a shift of paradigm in the study of algebraic structures.

Recently, neutrosophic algebraic structures have been introduced by the author with Florentin Smarandache [2006]. Smarandache algebraic structures simultaneously analyses two distinct algebraic structures. For instance, study of S-semigroups analyses semigroups and groups. Thus a need for the introduction of classical theorems for groups to S-neutrosophic loops, groupoids, groups, and semigroups has become an essentiality. S-neutrosophic groups, S-neutrosophic bigroups, S-neutrosophic N-groups, S-neutrosophic semigroups, S-neutrosophic loops and S-neutrosophic mixed N-structures are introduced in this book.

They are illustrated with examples. This book gives a lot of scope for the reader to develop the subject.

This book has seven chapters. In Chapter one, an elaborate recollection of Smarandache structures like S-semigroups, S-loops and S-groupoids is given. It also gives notions about N-ary algebraic structures and their Smarandache analogue, Neutrosophic structures viz. groups, semigroups, groupoids and loops are given in Chapter one to make the book a self-contained one. For the first time S-neutrosophic groups and S-neutrosophic N-groups are introduced in Chapter two and their properties are given. S-neutrosophic semigroups and S-neutrosophic N-semigroups are defined and discussed in Chapter three.

Chapter four defines S-neutrosophic S-loops and S-N-neutrosophic loops. A brief notion about S-neutrosophic groupoids and their generalizations are given Chapter five. Chapter six gives S-neutrosophic mixed N-structures and their duals. Chapter seven gives 68 problems for any interested reader.

I deeply acknowledge the help rendered by Kama and Meena who have formatted this book. I also thank Dr. K.Kandasamy for proof reading.

W.B.VASANTHA KANDASAMY
2006

**Chapter One**

# INTRODUCTION

In this chapter we introduce certain basic concepts to make this book a self contained one. This chapter has 5 sections. In section one the notion of groups and N-groups are introduced. Section two just mentions about semigroups and N-semigroups. In section 3 loops and N-loops are recalled. Section 4 gives a brief description about groupoids and their properties. Section 5 recalls the mixed N algebraic structure.

In this chapter we give most of the important definitions about neutrosophic structures of groups, semigroups, loops and groupoids. We are forced to do it for we have only reference book [51]. So a reader should know at least some thing of it to understand the theory. Secondly atleast basic definitions about Smarandache semigroups and Smarandache groupoids are given.

This chapter has 5 sections. Section 1 gives definitions of N-groups and their related neutrosophic notions. Section two gives briefly about S-semigroups and neutrosophic semigroups and their generalization. Section 3 is devoted to loops i.e., some notions of S-loops and neutrosophic loops are given. The notions needed about groupoids viz. S-groupoids are given in section 4. The final section gives the mixed structures and their neutrosophic analogue.

## 1.1 Groups, N-groups and neutrosophic groups and neutrosophic N-groups

It is a well-known fact that groups are the only algebraic structures with a single binary operation that is mathematically so perfect that an introduction of a richer structure within it is impossible. Now we proceed on to define a group. In this section we just recall the most important definitions about N-groups neutrosophic groups, neutrosophic bigroups, and neutrosophic N-groups, (which we choose to call as their neutrosophic analogue). The classical theorem on group theory is also just given to make the book self-contained one.

**DEFINITION 1.1.1:** *A non empty set of elements G is said to form a group if in G there is defined a binary operation, called the product and denoted by '•' such that*

i. *$a, b \in G$ implies that $a \bullet b \in G$ (closed).*
ii. *$a, b, c \in G$ implies $a \bullet (b \bullet c) = (a \bullet b) \bullet c$ (associative law).*
iii. *There exists an element $e \in G$ such that $a \bullet e = e \bullet a = a$ for all $a \in G$ (the existence of identity element in G).*
iv. *For every $a \in G$ there exists an element $a^{-1} \in G$ such that $a \bullet a^{-1} = a^{-1} \bullet a = e$ (the existence of inverse in G).*

**DEFINITION 1.1.2:** *A subgroup N of a group G is said to be a normal subgroup of G if for every $g \in G$ and $n \in N$, $g\, n\, g^{-1} \in N$.*

Equivalently by $gNg^{-1}$ we mean the set of all $gng^{-1}$, $n \in N$ then N is a normal subgroup of G if and only if $gNg^{-1} \subset N$ for every $g \in G$.

**THEOREM 1.1.1:** N is a normal subgroup of G if and only if $gNg^{-1} = N$ for every $g \in G$.

**DEFINITION 1.1.3:** *Let G be a group. $Z(G) = \{x \in G \mid gx = xg$ for all $g \in G\}$. Then Z(G) is called the center of the group G.*

**DEFINITION 1.1.4:** *Let G be a group, A, B be subgroups of G. If x, y $\in$ G define x $\sim$ y if y = axb for some a $\in$ A and b $\in$ B. We call the set AxB = {axb / a $\in$ A, b $\in$ B} a double coset of A, B in G.*

**DEFINITION 1.1.5:** *Let G be a group. A and B subgroups of G, we say A and B are conjugate with each other if for some g $\in$ G, A = $gBg^{-1}$.*

*Clearly if A and B are conjugate subgroups of G then o(A) = o(B).*

**THEOREM: (LAGRANGE).** If G is a finite group and H is a subgroup of G then o(H) is a divisor of o(G).

**COROLLARY 1.1.1:** *If G is a finite group and a $\in$ G, then o(a) / o(G).*

**COROLLARY 1.1.2:** *If G is a finite group and a $\in$ G, then $a^{o(G)}$ = e.*

In this section we give the two Cauchy theorems one for abelian groups and the other for non-abelian groups. The main result on finite groups is that if the order of the group is n (n < $\propto$) if p is a prime dividing n by Cauchy's theorem we will always be able to pick up an element a $\in$ G such that $a^p$ = e. In fact we can say Sylow's theorem is a partial extension of Cauchy's theorem for he says this finite group G has a subgroup of order $p^{\alpha}$($\alpha \geq 1$, p, a prime).

**THEOREM: (CAUCHY'S THEOREM FOR ABELIAN GROUPS).** *Suppose G is a finite abelian group and p / o(G), where p is a prime number. Then there is an element a $\neq$ e $\in$ G such that $a^p$ = e.*

**THEOREM: (CAUCHY):** *If p is a prime number and p | o(G), then G has an element of order p.*

Though one may marvel at the number of groups of varying types carrying many different properties, except for Cayley's we would not have seen them to be imbedded in the class of groups this was done by Cayley's in his famous theorem. Smarandache semigroups also has a beautiful analog for Cayley's theorem.

**THEOREM: (CAYLEY)** *Every group is isomorphic to a subgroup of A(S) for some appropriate S.*

The Norwegian mathematician Peter Ludvig Mejdell Sylow was the contributor of Sylow's theorems. Sylow's theorems serve double purpose. One hand they form partial answers to the converse of Lagrange's theorem and on the other hand they are the complete extension of Cauchy's Theorem. Thus Sylow's work interlinks the works of two great mathematicians Lagrange and Cauchy. The following theorem is one, which makes use of Cauchy's theorem. It gives a nice partial converse to Lagrange's theorem and is easily understood.

**THEOREM: (SYLOW'S THEOREM FOR ABELIAN GROUPS)** *If G is an abelian group of order o(G), and if p is a prime number, such that* $p^{\alpha} \mid o(G)$, $p^{\alpha+1} \nmid o(G)$, *then G has a subgroup of order* $p^{\alpha}$.

**COROLLARY 1.1.3:** *If G is an abelian group of finite order and* $p^{\alpha} \mid o(G)$, $p^{\alpha+1} \nmid o(G)$, *then there is a unique subgroup of G of order* $p^{\alpha}$.

**DEFINITION 1.1.6:** *Let G be a finite group. A subgroup G of order* $p^{\alpha}$, *where* $p^{\alpha} / o(G)$ *but* $p^{\alpha} \nmid o(G)$, *is called a p-Sylow subgroup of G. Thus we see that for any finite group G if p is any prime which divides o(G); then G has a p-Sylow subgroup.*

**THEOREM (FIRST PART OF SYLOW'S THEOREM):** *If p is a prime number and* $p^{\alpha}/ o(G)$ *and* $p^{\alpha+1} \nmid o(G)$, *G is a finite group, then G has a subgroup of order* $p^{\alpha}$.

**THEOREM: (SECOND PART OF SYLOW'S THEOREM):** *If G is a finite group, p a prime and $p^n \mid o(G)$ but $p^{n+1} \nmid o(G)$, then any two subgroup of G of order $p^n$ are conjugate.*

**THEOREM: (THIRD PART OF SYLOW'S THEOREM):** *The number of p-Sylow subgroups in G, for a given prime, is of the form 1 + kp.*

**DEFINITION 1.1.7:** *Let $\{G, *_1, \ldots, *_N\}$ be a non empty set with N binary operations. $\{G, *_1, \ldots, *_N\}$ is called a N-group if there exists N proper subsets $G_1, \ldots, G_N$ of G such that*

i. $G = G_1 \cup G_2 \ldots \cup G_N$.
ii. $(G_i, *_i)$ *is a group for* $i = 1, 2, \ldots, N$.

*We say proper subset of G if $G_i \not\subseteq G_j$ and $G_j \not\subseteq G_i$ if $i \neq j$ for $1 \leq i, j \leq N$. When N = 2 this definition reduces to the definition of bigroup.*

**DEFINITION 1.1.8:** *Let $\{G, *_1, \ldots, *_N\}$ be a N-group. A subset H $(\neq \phi)$ of a N-group $(G, *_1, \ldots, *_N)$ is called a sub N-group if H itself is a N-group under $*_1, *_2, \ldots, *_N$, binary operations defined on G.*

**THEOREM 1.1.2:** *Let $(G, *_1, \ldots, *_N)$ be a N-group. The subset $H \neq \phi$ of a N-group G is a sub N-group then $(H, *_i)$ in general are not groups for $i = 1, 2, \ldots, N$.*

**DEFINITION 1.1.9:** *Let $(G, *_1, \ldots, *_N)$ be a N-group where $G = G_1 \cup G_2 \cup \ldots \cup G_N$. Let $(H, *_1, \ldots, *_N)$ be a sub N-group of $(G, *_1, \ldots, *_N)$ where $H = H_1 \cup H_2 \cup \ldots \cup H_N$ we say $(H, *_1, \ldots, *_N)$ is a normal sub N-group of $(G, *_1, \ldots, *_N)$ if each $H_i$ is a normal subgroup of $G_i$ for $i = 1, 2, \ldots, N$.*

*Even if one of the subgroups $H_i$ happens to be non normal subgroup of $G_i$ still we do not call H a normal sub-N-group of the N-group G.*

**DEFINITION 1.1.10:** *Let $(G = G_1 \cup G_2 \cup ... \cup G_N, *_1, *_2, ..., *_N)$ and $(K = K_1 \cup K_2 \cup ... \cup K_N, *_1, ..., *_N)$ be any two N- groups. We say a map $\phi : G \to K$ to be a N-group homomorphism if $\phi$ | $G_i$ is a group homomorphism from $G_i$ to $K_i$ for i = 1, 2,..., N. i.e. $\phi\big|_{G_i} : G_i \to K_i$ is a group homomorphism of the group $G_i$ to the group $K_i$; for i = 1, 2, ..., N.*

**DEFINITION 1.1.11:** *Let $G = (G_1 \cup G_2 \cup ... \cup G_N, *_1, *_2, ..., *_N)$ be a non empty set such that $G_i$ are proper subsets of G. $G = (G_1 \cup G_2 \cup ... \cup G_N, *_1, *_2, ..., *_N)$ is called a Smarandache N-group (S-N-group) if the following conditions are satisfied.*

i. $(G_i, *_i)$ *is a group* $(1 \le i \le N)$.
ii. $(G_j, *_j)$ *is a S-semigroup for some* $j \ne i$, $1 \le j \le N$.

**DEFINITION 1.1.12:** *Let $G = (G_1 \cup G_2 \cup ... \cup G_N, *_1, ..., *_N)$ be a S-N-group; a proper subset P of G is said to be a Smarandache sub N-group (S-sub N-group) of G if $P = P_1 \cup P_2 \cup ... \cup P_N$ where $P_i \subseteq G_i$, i =1, 2, ..., N, and at least one of the $P_i$ is a S-semigroup under the operations of $G_i$. In short a proper subset P of G is a S-sub N-group if P itself is a Smarandache-N-group under the operations of G.*

**DEFINITION 1.1.13:** *Let $(G = G_1 \cup G_2 \cup ... \cup G_N, *_1, ..., *_N)$ be a S-N-group. We say G is a Smarandache commutative N-group (S-commutative N-group) if all $G_i$ which are not semigroups are commutative and when $G_j$ are S-semigroups every proper subset which is a group in $G_j$ is commutative.*

**DEFINITION 1.1.14:** *Let $\{G = G_1 \cup G_2 \cup ... \cup G_N, *_1, ..., *_N\}$be a -N-group. We say G is Smarandache weakly commutative N-group (S-weakly commutative N-group) if in all the S-semigroup $G_j$, they have a proper subset which is a commutative group and the rest of $G_i$'s which are groups are commutative.*

**THEOREM 1.1.3:** *Every S-commutative N-group is a S-weakly commutative N-group. But a S-weakly commutative N-group in general need not be a S-commutative N-group.*

**DEFINITION 1.1.15:** *Let* $G = (G_1 \cup G_2 \cup \ldots \cup G_N, *_1, *_2, \ldots, *_N)$ *be a S- N-group. We say G is a Smarandache cyclic N-group (S-cyclic N-group) if* $G_i$ *is a group it must be cyclic and if* $G_j$ *is a S-semigroup it must be a S-cyclic semigroup for* $1 \le i, j \le N$.

**DEFINITION 1.1.16:** *Let* $G = (G_1 \cup G_2 \cup \ldots \cup G_N, *_1, *_2, \ldots, *_N)$ *be a S-N-group. We say G is a Smarandache weakly cyclic N-group (S-weakly cyclic N-group) if every group in the collection* $\{G_i\}$ *is a cyclic group and every S-semigroup in the collection* $\{G_j\}$ *is a S-weakly cyclic semigroup.*

**DEFINITION 1.1.17:** *Let* $G = (G_1 \cup G_2 \cup \ldots \cup G_N, *_1, *_2, \ldots, *_N)$ *be a S-N-group; let* $H = (H_1 \cup H_2 \cup \ldots \cup H_N, *_1, *_2, \ldots, *_N)$ *be a S-sub N group of G.*
*For any* $g \in G$ *we define the Smarandache right coset (S-right coset) of the N-group G as* $Ha = H_1 \cup H_2 \cup \ldots \cup H_i a \cup \ldots \cup H_N$ *if* $a \in H_i$ *alone; if* $a \in \bigcap_{i=1}^{K} H_i$, *then* $Ha = H_1'a \cup H_2' a \cup \ldots \cup H_N' a$; *here* $H_i = H'_i$ *if* $H_i$ *is a group;* $H_i' \subset H_i$ *if* $H_i$ *is a S-semigroup,* $H'_i$ *is a group. Similarly we can define Smarandache left coset (S-left coset) a H. We say H is a Smarandache coset if a H = Ha.*

**DEFINITION 1.1.18:** $G = (G_1 \cup G_2 \cup \ldots \cup G_N, *_1, *_2, \ldots, *_N)$ *where each* $G_i$ *is either* $S(n_i)$ *or* $S_{n_i}$ *i.e. each* $G_i$ *is either a symmetric semigroup or a symmetric group, for* $i = 1, 2, \ldots, N$. *Then G is defined as the Smarandache symmetric N-group (S-symmetric N-group).*

**THEOREM 1.1.4: (SMARANDACHE CAYLEY'S THEOREM FOR S-N-GROUPS):** *Let* $G = G_1 \cup \ldots \cup G_N, *_1, *_2, \ldots, *_N)$ *be a S- N-group. Every S-N-group is embeddable in a S-symmetric N-group for a suitable* $n_i$, $i = 1, 2, \ldots, N$.

**DEFINITION 1.1.19:** *Let $G = (G_1 \cup G_2 \cup \ldots \cup G_N, *_1, *_2, \ldots, *_N)$ be a S-N group. We say an element $x \in G$ has a Smarandache conjugate (S-conjugate) y in G if.*

i. *$xa = ay$ for some $a \in G$.*
ii. *$ab = bx$ and $ac = cy$ for some b, c in G.*

*It is easy to verify if G is a S-symmetric N-group then G has S-elements which are S-conjugate.*

**DEFINITION 1.1.20:** *Let (G, *) be any group, the neutrosophic group is generated by I and G under * denoted by $N(G) = \{\langle G \cup I\rangle, *\}$.*

**THEOREM 1.1.5:** *Let (G, *) be a group, $N(G) = \{\langle G \cup I\rangle, *\}$ be the neutrosophic group.*

i. *N(G) in general is not a group.*
ii. *N(G) always contains a group.*

**DEFINITION 1.1.21:** *Let $N(G) = \langle G \cup I\rangle$ be a neutrosophic group generated by G and I. A proper subset P(G) is said to be a neutrosophic subgroup if P(G) is a neutrosophic group i.e. P(G) must contain a (sub) group.*

**DEFINITION 1.1.22:** *Let N(G) be a finite neutrosophic group. Suppose L is a pseudo neutrosophic subgroup of N(G) and if o(L) / o(N(G)) then we call L to be a pseudo Lagrange neutrosophic subgroup. If all pseudo neutrosophic subgroups of N(G) are pseudo Lagrange neutrosophic subgroups then we call N(G) to be a pseudo Lagrange neutrosophic group.*

*If N(G) has atleast one pseudo Lagrange neutrosophic subgroup then we call N(G) to be a weakly pseudo Lagrange neutrosophic group. If N(G) has no pseudo Lagrange neutrosophic subgroup then we call N(G) to be pseudo Lagrange free neutrosophic group.*

**DEFINITION 1.1.23:** *Let N(G) be a neutrosophic group. We say a neutrosophic subgroup H of N(G) is normal if we can find x and y in N(G) such that H =xHy for all x, y $\in$ N (G) (Note x = y or y = $x^{-1}$ can also occur).*

**DEFINITION 1.1.24:** *A neutrosophic group N(G) which has no nontrivial neutrosophic normal subgroup is called a simple neutrosophic group.*

Now we define pseudo simple neutrosophic groups.

**DEFINITION 1.1.25:** *Let N(G) be a neutrosophic group. A proper pseudo neutrosophic subgroup P of N(G) is said to be normal if we have P = xPy for all x, y $\in$ N(G). A neutrosophic group is said to be pseudo simple neutrosophic group if N(G) has no nontrivial pseudo normal subgroups.*

We do not know whether there exists any relation between pseudo simple neutrosophic groups and simple neutrosophic groups.

Now we proceed on to define the notion of right (left) coset for both the types of subgroups.

**DEFINITION 1.1.26:** *Let L (G) be a neutrosophic group. H be a neutrosophic subgroup of N(G) for n $\in$ N(G), then H n = {hn / h $\in$ H} is called a right coset of H in G.*

**DEFINITION 1.1.27:** *Let N(G) be a neutrosophic group. K be a pseudo neutrosophic subgroup of N(G). Then for a $\in$ N(G), Ka = {ka | k $\in$ K} is called the pseudo right coset of K in N(G).*

**DEFINITION 1.1.28:** *Let N(G) be a finite neutrosophic group. If for a prime $p^{\alpha}$ / o(N(G)) and $p^{\alpha+1} \nmid$ o(N(G)), N(G) has a neutrosophic subgroup P of order $p^{\alpha}$ then we call P a p-Sylow neutrosophic subgroup of N(G).*

*Now if for every prime p such that $p^{\alpha}$ / o(N(G)) and $p^{\alpha+1} \nmid$ o(N(G)) we have an associated p-Sylow neutrosophic subgroup then we call N(G) a Sylow neutrosophic group.*

*If N(G) has atleast one p-Sylow neutrosophic subgroup then we call N(G) a weakly Sylow neutrosophic group. If N(G) has no p-Sylow neutrosophic subgroup then we call N(G) a Sylow free neutrosophic group.*

**DEFINITION 1.1.29:** *Let $B_N$ (G) = {$B(G_1) \cup B(G_2)$, $*_1$, $*_2$} be a non empty subset with two binary operation on $B_N$ (G) satisfying the following conditions:*

i. $B_N$ *(G) = {$B(G_1) \cup B(G_2)$} where $B(G_1)$ and $B(G_2)$ are proper subsets of $B_N$ (G).*
ii. *($B(G_1)$, $*_1$) is a neutrosophic group.*
iii. *($B(G_2)$, $*_2$) is a group.*

*Then we define ($B_N(G)$, $*_1$, $*_2$) to be a neutrosophic bigroup. If both $B(G_1)$ and $B(G_2)$ are neutrosophic groups we say $B_N(G)$ is a strong neutrosophic bigroup. If both the groups are not neutrosophic group we see $B_N(G)$ is just a bigroup.*

**DEFINITION 1.1.30:** *Let $B_N(G)$ = {$B(G_1) \cup B(G_2)$, $*_1$, $*_2$} be a neutrosophic bigroup. A proper subset P = {$P_1 \cup P_2$, $*_1$, $*_2$} is a neutrosophic subbigroup of $B_N(G)$ if the following conditions are satisfied P = {$P_1 \cup P_2$, $*_1$, $*_2$} is a neutrosophic bigroup under the operations $*_1$, $*_2$ i.e. ($P_1$, $*_1$) is a neutrosophic subgroup of ($B_1$, $*_1$) and ($P_2$, $*_2$) is a subgroup of ($B_2$, $*_2$). $P_1 = P \cap B_1$ and $P_2 = P \cap B_2$ are subgroups of $B_1$ and $B_2$ respectively. If both of $P_1$ and $P_2$ are not neutrosophic then we call $P = P_1 \cup P_2$ to be just a bigroup.*

**DEFINITION 1.1.31:** *Let $B_N$ (G) = {$B(G_1) \cup B(G_2)$, $*_1$, $*_2$} be a neutrosophic bigroup. P(G) = {$P(G_1) \cup P(G_2)$, $*_1$, $*_2$} be a neutrosophic bigroup. P(G) = {$P(G_1) \cup P(G_2)$, $*_1$, $*_2$} is said to be a neutrosophic normal subbigroup of $B_N(G)$ if P(G) is a neutrosophic subbigroup and both $P(G_1)$ and $P(G_2)$ are normal subgroups of $B(G_1)$ and $B(G_2)$ respectively.*

**DEFINITION 1.1.32:** *Let $B_N(G)$ = {$B(G_1) \cup B(G_2)$, $*_1$, $*_2$} be a neutrosophic bigroup. Suppose P = {$P(G_1) \cup P(G_2)$, $*_1$, $*_2$} and*

*$K = \{K(G_1) \cup K(G_2), *_1, *_2\}$ be any two neutrosophic subbigroups we say P and K are conjugate if each $P(G_i)$ is conjugate with $K(G_i)$, i = 1, 2, then we say P and K are neutrosophic conjugate subbigroups of $B_N(G)$.*

**DEFINITION 1.1.33:** *Let $B_N(G) = \{B(G_1) \cup B(G_2), *_1, *_2\}$ be any neutrosophic bigroup. The neutrosophic bicentre of the bigroup $B_N(G)$ denoted $C_N(G) = C(G_1) \cup C(G_2)$ where $C(G_1)$ is the centre of $B(G_1)$ and $C(G_2)$ is the centre of $B(G_2)$. If the neutrosophic bigroup is commutative then $C_N(G) = B_N(G)$.*

**DEFINITION 1.1.34:** *A subset $H \neq \phi$ of a strong neutrosophic bigroup $(\langle G \cup I\rangle, *, o)$ is called a strong neutrosophic subbigroup if H itself is a strong neutrosophic bigroup under '*' and 'o' operations defined on $\langle G \cup I\rangle$.*

**THEOREM 1.1.6:** *Let $(\langle G \cup I\rangle, +, o)$ be a strong neutrosophic bigroup. A subset $H \neq \phi$ of a strong neutrosophic bigroup $\langle G \cup I\rangle$ is a neutrosophic subbigroup then (H, +) and (H, o) in general are not neutrosophic groups.*

**THEOREM 1.1.7:** *Let $\{\langle G \cup I\rangle, +, o\}$ be a strong neutrosophic bigroup. Then the subset H ($\neq \phi$) is a strong neutrosophic subbigroup of $\langle G \cup I\rangle$ if and only if there exists two proper subsets $\langle G_1 \cup I\rangle$, $\langle G_2 \cup I\rangle$ of $\langle G \cup I\rangle$ such that*

*i.* *$\langle G \cup I\rangle = \langle G_1 \cup I\rangle \cup \langle G_2 \cup I\rangle$ with $(\langle G_1 \cup I\rangle, +)$ is a neutrosophic group and $(\langle G_2 \cup I\rangle, o)$ a neutrosophic group.*
*ii.* *$(H \cap \langle G_i \cup I\rangle, +)$ is a neutrosophic subbigroup of $\langle G_i \cup I\rangle$ for i = 1, 2.*

**DEFINITION 1.1.35:** *Let $(\langle G \cup I\rangle, +, o)$ and $(\langle K \cup I\rangle, o', \oplus)$ be any two strong neutrosophic bigroups where $\langle G \cup I\rangle = \langle G_1 \cup I\rangle \cup \langle G_2 \cup I\rangle$ and $\langle K \cup I\rangle = \langle K_1 \cup I\rangle \cup \langle K_2 \cup I\rangle$. We say a bimap $\phi = \phi_1 \cup \phi_2$: $\langle G \cup I\rangle \to \langle K \cup I\rangle$ (Here $\phi_1(I) = I$ and $\phi_2(I) = I$) is said to be a strong neutrosophic bigroup bihomomorphism if $\phi_1 = \phi / \langle G_1 \cup I\rangle$ and $\phi_2 = \phi / \langle G_2 \cup I\rangle$ where $\phi_1$ and $\phi_2$ are* *neutrosophic group homomorphism from $\langle G_1 \cup I\rangle$ to $\langle K_1 \cup I\rangle$ and $\langle G_2 \cup I\rangle$ to $\langle K_2 \cup I\rangle$ respectively.*

**DEFINITION 1.1.36:** *Let $G = \langle G_1 \cup I, *, \oplus\rangle$, be a neutrosophic bigroup. We say two neutrosophic strong subbigroups $H = H_1 \cup H_2$ and $K = K_1 \cup K_2$ are conjugate neutrosophic subbigroups of $\langle G \cup I\rangle = \langle G_1 \cup I\rangle \cup (\langle G_2 \cup I\rangle)$ if $H_1$ is conjugate to $K_1$ and $H_2$ is conjugate to $K_2$ as neutrosophic subgroups of $\langle G_1 \cup I\rangle$ and $\langle G_2 \cup I\rangle$ respectively.*

**DEFINITION 1.1.37:** *Let $(\langle G \cup I\rangle, *, o) = (\langle G \cup I\rangle, *) \cup (\langle G \cup I\rangle, o)$ be a neutrosophic bigroup. The normalizer of a in $\langle G \cup I\rangle$ is the set $N(a) = \{x \in \langle G \cup I\rangle \mid x a = ax\} = N_1(a) \cup N_2(a) = \{x_1 \in \langle G_1 \cup I\rangle \mid x_1 a = ax_1\} \cup \{x_2 \in \langle G_2 \cup I\rangle \mid x_2 a = ax_2\}$ if $a \in (G_1 \cup I) \cap (G_2 \cup I)$ if $a \in (\langle G_1 \cup I\rangle)$ and $a \notin (\langle G_2 \cup I\rangle)$, $N_2(a) = \phi$ like wise if $a \notin \langle G_1 \cup I\rangle$ and $a \in (\langle G_2 \cup I\rangle)$ then $N_1(a) = \phi$ and $N_2(a) = N(a)$ is a neutrosophic subbigroup of $\langle G \cup I\rangle$, clearly $N(a) \neq \phi$ for $I \in N(a)$.*

**DEFINITION 1.1.38:** *Let $(\langle G \cup I\rangle, o, *) = (\langle G_1 \cup I\rangle, o) \cup (\langle G_2 \cup I\rangle, *)$ be a neutrosophic bigroup. Let $H = H_1 \cup H_2$ be a strong neutrosophic subbigroup of $\langle G \cup I\rangle$. The right bicoset of H in $\langle G \cup I\rangle$ for some a in $\langle G \cup I\rangle$ is defined to be $Ha = \{h_1 a \mid h_1 \in H_1$ and $a \in G_1 \cap G_2\} \cup \{h_2 a \mid h_2 \in H_2$ and $a \in G_1 \cap G_2\}$ if $a \in G_1$ and $a \notin G_2$ then $Ha = \{h_1 a \mid h_1 \in H_1\} \cup H_2$ and $a \notin G_1$ then $Ha = \{h_2 a \mid h_2 \in H_2\} \cup H_1$ and $a \in G_2$.*

**DEFINITION 1.1.39:** *Let $\langle G \cup I\rangle = \langle G_1 \cup I\rangle \cup \langle G_2 \cup I\rangle$ be a neutrosophic bigroup we say $N = N_1 \cup N_2$ is a neutrosophic normal subbigroup if and only if $N_1$ is a neutrosophic normal subgroup of $\langle G_1 \cup I\rangle$ and $N_2$ is a neutrosophic normal subgroup of $\langle G_2 \cup I\rangle$.*

*We define the strong neutrosophic quotient bigroup $\frac{\langle G \cup I\rangle}{N}$ as $\left[\frac{\langle G_1 \cup I\rangle}{N_1} \cup \frac{\langle G_2 \cup I\rangle}{N_2}\right]$ which is also a neutrosophic bigroup.*

**DEFINITION 1.1.40:** *Let* $(\langle G \cup I\rangle, *_1,\ldots, *_N)$ *be a nonempty set with N-binary operations defined on it. We say* $\langle G \cup I\rangle$ *is a strong neutrosophic N-group if the following conditions are true.*

i. $\langle G \cup I\rangle = \langle G_1 \cup I\rangle \cup \langle G_2 \cup I\rangle \cup \ldots \cup \langle G_N \cup I\rangle$ *where* $\langle G_i \cup I\rangle$ *are proper subsets of* $\langle G \cup I\rangle$.
ii. $(\langle G_i \cup I\rangle, *_i)$ *is a neutrosophic group,* $i = 1, 2,\ldots, N$. *If in the above definition we have*
    a. $\langle G \cup I\rangle = G_1 \cup \langle G_2 \cup I\rangle \cup \langle G_3 \cup I\rangle \ldots \cup G_K \cup G_{K+1} \cup \ldots \cup G_N$.
    b. $(G_i, *_i)$ *is a group for some i or*
iii. $(\langle G_j \cup I\rangle, *_j)$ *is a neutrosophic group for some j.*

*then we call* $\langle G \cup I\rangle$ *to be a neutrosophic N-group.*

**DEFINITION 1.1.41:** *Let* $(\langle G \cup I\rangle = \langle G_1 \cup I\rangle \cup \langle G_2 \cup I\rangle \cup \ldots \cup \langle G_N \cup I\rangle, *_1, \ldots, *_N)$ *be a neutrosophic N-group. A proper subset* $(P, *_1,\ldots, *_N)$ *is said to be a neutrosophic sub N-group of* $\langle G \cup I\rangle$ *if* $P = (P_1 \cup \ldots \cup P_N)$ *and each* $(P_i, *_i)$ *is a neutrosophic subgroup (subgroup) of* $(G_i, *_i)$, $1 \leq i \leq N$.

**DEFINITION 1.1.42:** *Let* $(\langle G \cup I\rangle = \langle G_1 \cup I\rangle \cup \langle G_2 \cup I\rangle \cup \ldots \cup \langle G_N \cup I\rangle, *_1, \ldots, *_N)$ *be a strong neutrosophic N-group. If* $\langle G \cup I\rangle$ *is a Sylow strong neutrosophic N-group and if for every prime p such that* $p^{\alpha} / o(\langle G \cup I\rangle)$ *and* $p^{\alpha+1} \nmid o(\langle G \cup I\rangle)$ *we have a strong neutrosophic sub N-group of order* $p^{\alpha+t}$ $(t \geq 1)$ *then we call* $\langle G \cup I\rangle$ *a super Sylow strong neutrosophic N-group.*

**THEOREM 1.1.8:** *Suppose* $(\langle G \cup I\rangle = \langle G_1 \cup I\rangle \cup \ldots \cup \langle G_N \cup I\rangle, *_1, \ldots, *_N)$ *be a neutrosophic N-group of order p, p a prime. Then the following are true.*

i. $\langle G \cup I\rangle$ *is not a Cauchy neutrosophic N-group.*
ii. $\langle G \cup I\rangle$ *is not a semi Cauchy neutrosophic N group.*
iii. $\langle G \cup I\rangle$ *is not a weakly Cauchy neutrosophic N-group.*

*All elements of finite order are anti Cauchy elements and anti Cauchy neutrosophic elements.*

**DEFINITION 1.1.43:** *Let* $\{\langle G \cup I\rangle = \{\langle G_1 \cup I\rangle \cup \langle G_2 \cup I\rangle \cup ... \cup \langle G_N \cup I\rangle, *_1, ..., *_N\}$ *and* $\{\langle H \cup I\rangle = \langle H_1 \cup I\rangle \cup \langle H_2 \cup I\rangle \cup ... \cup \langle H_N \cup I\rangle, *_1, *_2, *_3, ..., *_N\}$ *be any two neutrosophic N-group such that if* $(\langle G_i \cup I\rangle, *_i)$ *is a neutrosophic group then* $(\langle H_i \cup I\rangle, *_i)$ *is also a neutrosophic group. If* $(G_t, *_t)$ *is a group then* $(H_t, *_t)$ *is a group.*

*A map* $\phi : \langle G \cup I\rangle$ *to* $\langle H \cup I\rangle$ *satisfying* $\phi(I) = I$ *is defined to be a N homomorphism if* $\phi_i = \phi \mid \langle G_i \cup I\rangle$ *(or* $\phi \mid G_i$*) then each* $\phi_i$ *is either a group homomorphism or a neutrosophic group homomorphism, we denote the N-homomorphism by* $\phi = \phi_1 \cup \phi_2 \cup ... \cup \phi_N : \langle G \cup I\rangle \to \langle H \cup I\rangle$*.*

**DEFINITION 1.1.44:** *Let* $\langle G \cup I\rangle = \{\langle G_1 \cup I\rangle \cup \langle G_2 \cup I\rangle \cup ... \cup \langle G_N \cup I\rangle, *_1, ..., *_N\}$ *be a neutrosophic N-group. Let* $H = \{H_1 \cup H_2 \cup ... \cup H_N\}$ *be a neutrosophic sub N-group of* $\langle G \cup I\rangle$*. We say H is a* $(p_1, p_2, ..., p_N)$ *Sylow neutrosophic sub N-group of* $\langle G \cup I\rangle$ *if* $H_i$ *is a* $p_i$ *- Sylow neutrosophic subgroup of* $G_i$*. If none of the* $H_i$*'s are neutrosophic subgroups of* $G_i$ *we call H a* $(p_1, p_2, ..., p_N)$ *Sylow free sub N-group.*

**DEFINITION 1.1.45:** *Let* $(\langle G \cup I\rangle = \langle G_1 \cup I\rangle \cup \langle G_2 \cup I\rangle \cup ... \cup \langle G_N \cup I\rangle, *_1, ..., *_N)$ *be a strong neutrosophic N-group. Suppose* $H = \{H_1 \cup H_2 \cup ... \cup H_N, *_1, ..., *_N\}$ *and* $K = \{K_1 \cup K_2 \cup ... \cup K_N, *_1, ..., *_N\}$ *are two neutrosophic sub N-groups of* $\langle G \cup I\rangle$*, we say K is a strong conjugate to H or H is conjugate to K if each* $H_i$ *is conjugate to* $K_i$ $(i = 1, 2, ..., N)$ *as subgroups of* $G_i$*.*

## 1.2 Semigroups, N-semigroups, S-semigroups, S-N-semigroups and their neutrosophic analogues

We in this section just briefly introduce the notion of S-semigroups and S-N-semigroups. Also give notion of neutrosophic semigroup and neutrosophic N-semigroups. For more literature please refer [50-1]. In this section we just recall the notion of semigroup, bisemigroup and N-semigroups. Also the notion of symmetric semigroups. For more refer [49-50].

**DEFINITION 1.2.1:** *Let (S, o) be a non empty set S with a closed, associative binary operation 'o' on S. (S, o) is called the semigroup i.e., for a, b $\in$ S, a o b $\in$ S.*

**DEFINITION 1.2.2:** *Let S(n) denote the set of all mappings of (1, 2, …, n) to itself S(n) under the composition of mappings is a semigroup. We call S(n) the symmetric semigroup of order $n^n$.*

**DEFINITION 1.2.3:** *Let ($S = S_1 \cup S_2$, *, o) be a non empty set with two binary operations * and o S is a bisemigroup if*

i. $S = S_1 \cup S_2$, *$S_1$ and $S_2$ are proper subsets of S.*
ii. *($S_1$, *) is a semigroup.*
iii. *($S_2$, o) is a semigroup.*

More about bisemigroups can be had from [48-50]. Now we proceed onto define N-semigroups.

**DEFINITION 1.2.4:** *Let $S = (S_1 \cup S_2 \cup \ldots \cup S_N, *_1, *_2, \ldots, *_N)$ be a non empty set with N binary operations. S is a N-semigroup if the following conditions are true.*

i. $S = S_1 \cup S_2 \cup \ldots \cup S_N$ *is such that each $S_i$ is a proper subset of S.*
ii. *($S_i$, $*_i$) is a semigroup for 1, 2, …, N.*

We just give an example.

***Example 1.2.1:*** Let $S = \{S_1 \cup S_2 \cup S_3 \cup S_4, *_1, *_2, *_3, *_4\}$ where

$S_1$ = $Z_{12}$, semigroup under multiplication modulo 12,
$S_2$ = S(4), symmetric semigroup,
$S_3$ = Z semigroup under multiplication and
$S_4$ = $\left\{\begin{pmatrix} a & b \\ c & d \end{pmatrix} \middle| a, b, c, d \in Z_{10}\right\}$ under matrix multiplication.

S is a 4-semigroup.

**DEFINITION 1.2.5:** *The Smarandache semigroup (S-semigroup) is defined to be a semigroup A such that a proper subset of A is a group (with respect to the same induced operation).*

**DEFINITION 1.2.6:** *Let S be a S-semigroup. If every proper subset of A in S, which is a group is commutative then we say the S-semigroup S to be a Smarandache commutative semigroup.*

**DEFINITION 1.2.7:** *Let S be S-semigroup, if S contains at least a proper subset A that is a commutative subgroup under the operations of S then we say S is a Smarandache weakly commutative semigroup.*

**DEFINITION 1.2.8:** *Let S be S-semigroup if every proper subset A of S which is a subgroup is cyclic then we say S is a Smarandache cyclic semigroup.*

**DEFINITION 1.2.9:** *Let S be a S-semigroup if there exists at least a proper subset A of S, which is a cyclic subgroup under the operations of S then we say S is a Smarandache weakly cyclic semigroup.*

**DEFINITION 1.2.10:** *Let S be a S-semigroup. If the number of distinct elements in S is finite, we say S is a finite S-semigroup otherwise we say S is a infinite S-semigroup.*

**THEOREM 1.2.1:** *Let G be a Smarandache commutative semigroup. G in general need not be a Smarandache cyclic semigroup.*

**THEOREM 1.2.2:** *S(n) is the S-semigroup.*

*Proof:* Clearly S(n) is the semigroup of order $n^n$. $S_n$ is a S-semigroup for it contains the symmetric group of degree n, i.e.

$S_n$ is a proper subset which is the group of permutations on (1, 2, 3, …, n). Hence S(n) is a S-semigroup.

Clearly S(n) is a not Smarandache commutative semigroup. Further we have the following engrossing results about S(n).

**DEFINITION 1.2.11:** *Let S be a S-semigroup. A proper subset A of S is said to be a Smarandache subsemigroup of S if A itself is a S-semigroup, that is A is a semigroup of S containing a proper subset B such that B is the group under the operations of S. Note we do not accept A to be a group. A must only be a semigroup.*

**DEFINITION 1.2.12:** *Let S be a S-semigroup. If A be a proper subset of S which is subsemigroup of S and A contains the largest group of S then we say A to be the Smarandache hyper subsemigroup of S.*

**THEOREM 1.2.3:** *Let S be a S-semigroup. Every Smarandache hyper subsemigroup is a Smarandache subsemigroup but every Smarandache subsemigroup is not a Smarandache hyper subsemigroup.*

**DEFINITION 1.2.13:** *Let S be a S-semigroup. We say S is a Smarandache simple semigroup if S has no proper subsemigroup A, which contains the largest subgroup of S or equivalently S has no Smarandache hyper subsemigroup.*

**DEFINITION 1.2.14:** *Let S be a finite S-semigroup. If the order of every subgroup of S divides the order of the S-semigroup S then we say S is a Smarandache Lagrange semigroup.*

**THEOREM 1.2.4:** *$Z_p = \{0, 1, 2, \ldots , p-1\}$ where p is a prime is a S-semigroup under multiplication modulo p. But $Z_p$ is a Smarandache simple semigroup.*

**DEFINITION 1.2.15:** *Let S be a finite S-semigroup. If there exists at least one subgroup A that is a proper subset ($A \subset S$) having the same operations of S whose order divides the order of S then we say that S is a Smarandache weakly Lagrange semigroup.*

**THEOREM 1.2.5: (CAYLEY'S THEOREM FOR S-SEMIGROUP)** *Every S-semigroup is isomorphic to a S-semigroup S(N); of mappings of a set N to itself, for some appropriate set N.*

**DEFINITION 1.2.16:** *Let $S = (S_1 \cup \ldots \cup S_N, *_1, \ldots, *_N)$ where S is a non empty set on which is defined, N-binary operations $*_1, *_2, \ldots, *_N$. S is called a N-semigroup if the following condition are true*

i. $S = S_1 \cup S_2 \cup \ldots \cup S_N$; *$S_i$'s are proper subsets of S.* $1 \leq i \leq N$.
ii. $(S_i, *i)$ *are semigroups,* $i = 1, 2, \ldots, N$.

**DEFINITION 1.2.17:** *Let $S = \{S_1 \cup S_2 \cup \ldots \cup S_N, *_1, *_2, \ldots, *_N\}$ where S is a non empty set and $*_1, \ldots, *_N$ are N-binary operations defined on S. S is said to be Smarandache N semigroup (S-N semigroup) if the following conditions are satisfied*

i. $S = S_1 \cup S_2 \cup \ldots \cup S_N$ *is such that $S_i$'s are proper subsets of S.*
ii. *Some of $(S_i, *i)$ are groups and some of $(S_j, *j)$ are S-semigroups,* $1 \leq i, j \leq N$. $(i \neq j)$

**DEFINITION 1.2.18:** *Let S be a semigroup, the semigroup generated by S and I i.e. $S \cup I$ denoted by $\langle S \cup I \rangle$ is defined to be a neutrosophic semigroup.*

**DEFINITION 1.2.19:** *Let N(S) be a neutrosophic semigroup. A proper subset P of N(S) is said to be a neutrosophic subsemigroup, if P is a neutrosophic semigroup under the operations of N (S). A neutrosophic semigroup N(S) is said to have a subsemigroup if N(S) has a proper subset which is a semigroup under the operations of N(S).*

**DEFINITION 1.2.20:** *A neutrosophic semigroup N(S) which has an element e in N(S) such that $e * s = s * e = s$ for all $s \in N(S)$, is called as a neutrosophic monoid.*

**DEFINITION 1.2.21:** *Let N(S) be a neutrosophic semigroup under a binary operation *. P be a proper subset of N(S). P is said to be a neutrosophic ideal of N(S) if the following conditions are satisfied.*

i. *P is a neutrosophic semigroup.*
ii. *for all $p \in P$ and for all $s \in N(S)$ we have $p * s$ and $s * p$ are in P.*

**DEFINITION 1.2.22:** *Let N(S) be a neutrosophic semigroup. A neutrosophic ideal P of N(S) is said to be maximal if $P \subset J \subset N(S)$, J a neutrosophic ideal then either J = P or J = N(S). A neutrosophic ideal M of N(S) is said to be minimal if $\phi \neq T \subseteq M \subseteq N(S)$ then T = M or T = $\phi$.*

We cannot always define the notion of neutrosophic cyclic semigroup but we can always define the notion of neutrosophic cyclic ideal of a neutrosophic semigroup N(S).

**DEFINITION 1.2.23:** *Let N(S) be a neutrosophic semigroup. P be a neutrosophic ideal of N (S), P is said to be a neutrosophic cyclic ideal or neutrosophic principal ideal if P can be generated by a single element.*

We proceed on to define the notion of neutrosophic symmetric semigroup.

**DEFINITION 1.2.24:** *Let S(N) be the neutrosophic semigroup. If S(N) contains a subsemigroup isomorphic to S(n) i.e. the semigroup of all mappings of the set (1, 2, 3, …, n) to itself under the composition of mappings, for a suitable n then we call S (N) the neutrosophic symmetric semigroup.*

**DEFINITION 1.2.25:** *Let (BN(S), *, o) be a nonempty set with two binary operations * and o. (BN(S), *, o) is said to be a neutrosophic bisemigroup if BN(S) = $P_1 \cup P_2$ where atleast one of ($P_1$, *) or ($P_2$, o) is a neutrosophic semigroup and other is* *just a semigroup.* $P_1$ *and* $P_2$ *are proper subsets of BN(S), i.e.* $P_1 \not\subseteq P_2$.

**DEFINITION 1.2.26:** *Let* $(BN(S) = P_1 \cup P_2; o, *)$ *be a neutrosophic bisemigroup. A proper subset (T, o, *) is said to be a neutrosophic subbisemigroup of BN (S) if*

i. $T = T_1 \cup T_2$ *where* $T_1 = P_1 \cap T$ *and* $T_2 = P_2 \cap T$ *and*
ii. *At least one of* $(T_1, o)$ *or* $(T_2, *)$ *is a neutrosophic semigroup.*

**DEFINITION 1.2.27:** *Let* $(BN(S) = P_1 \cup P_2, o, *)$ *be a neutrosophic strong bisemigroup. A proper subset T of BN (S) is called the strong neutrosophic subbisemigroup if* $T = T_1 \cup T_2$ *with* $T_1 = P_1 \cap T$ *and* $T_2 = P_2 \cap T_2$ *and if both* $(T_1, *)$ *and* $(T_2, o)$ *are neutrosophic subsemigroups of* $(P_1, *)$ *and* $(P_2, o)$ *respectively. We call* $T = T_1 \cup T_2$ *to be a neutrosophic strong subbisemigroup, if atleast one of* $(T_1, *)$ *or* $(T_2, o)$ *is a semigroup then* $T = T_1 \cup T_2$ *is only a neutrosophic subsemigroup.*

**DEFINITION 1.2.28:** *Let (BN (S), *, o) be a strong neutrosophic bisemigroup where* $BN(S) = P_1 \cup P_2$ *with* $(P_1, *)$ *and* $(P_2, o)$ *be any two neutrosophic semigroups.*

*Let J be a proper subset of BN(S) where* $I = I_1 \cup I_2$ *with* $I_1 = J \cap P_1$ *and* $I_2 = J \cap P_2$ *are neutrosophic ideals of the neutrosophic semigroups* $P_1$ *and* $P_2$ *respectively. Then I is called or defined as the strong neutrosophic biideal of B(N(S)).*

**DEFINITION 1.2.29:** *Let* $(BN(S) = P_1 \cup P_2\ *, o)$ *be any neutrosophic bisemigroup. Let J be a proper subset of B(NS) such that* $J_1 = J \cap P_1$ *and* $J_2 = J \cap P_2$ *are ideals of* $P_1$ *and* $P_2$ *respectively. Then J is called the neutrosophic biideal of BN(S).*

**DEFINITION 1.2.30:** *Let* $(BN(S)\ P_1 \cup P_2\ *, o)$ *be a neutrosophic strong bisemigroup.*

*Suppose I is a neutrosophic strong biideal of BN(S) i.e. $I_1 = P_1 \cap I$ and $I_2 = P_2 \cap I$ we say I is a neutrosophic strong maximal biideal of B (N(S)) if $I_1$ is the maximal ideal of $P_1$ and $I_2$ is the maximal ideal of $P_2$.*

*If only one of $I_1$ or $I_2$ alone is maximal then we call I to be a neutrosophic strong quasi maximal biideal of BN (S).*

**DEFINITION 1.2.31:** *Let {S(N), $*_1$, …, $*_N$} be a non empty set with N-binary operations defined on it. We call S(N) a neutrosophic N-semigroup (N a positive integer) if the following conditions are satisfied.*

i. *$S(N) = S_1 \cup \ldots \cup S_N$ where each $S_i$ is a proper subset of S(N) i.e. $S_i \not\subseteq S_j$ or $S_j \not\subseteq S_i$ if $i \neq j$.*
ii. *($S_i$, $*_i$) is either a neutrosophic semigroup or a semigroup for i = 1, 2, …, N.*

**DEFINITION 1.2.32:** *Let $S(N) = \{S_1 \cup S_2 \cup \ldots \cup S_N, *_1, \ldots, *_N\}$ be a neutrosophic N-semigroup. A proper subset $P = \{P_1 \cup P_2 \cup \ldots \cup P_N, *_1, *_2, \ldots, *_N\}$ of S(N) is said to be a neutrosophic N-subsemigroup if (1) $P_i = P \cap S$, i = 1, 2,…, N are subsemigroups of $S_i$ in which atleast some of the subsemigroups are neutrosophic subsemigroups.*

**DEFINITION 1.2.33:** *Let $S(N) = \{S_1 \cup S_2 \cup \ldots \cup S_N, *_1, \ldots, *_N\}$ be a neutrosophic strong N-semigroup. A proper subset $T = \{T_1 \cup T_2 \cup \ldots \cup T_N, *_1, \ldots, *_N\}$ of S(N) is said to be a neutrosophic strong sub N-semigroup if each ($T_i$, $*_i$) is a neutrosophic subsemigroup of ($S_i$, $*_i$) for i = 1, 2,…, N where $T_i = T \cap S_i$.*

**DEFINITION 1.2.34:** *Let $S(N) = \{S_1 \cup S_2 \cup \ldots \cup S_N, *_1, \ldots, *_N\}$ be a neutrosophic strong N-semigroup. A proper subset $J = \{I_1 \cup I_2 \cup \ldots \cup I_N\}$ where $I_t = J \cap S_t$ for t = 1, 2, …, N is said to be a neutrosophic strong N-ideal of S(N) if the following conditions are satisfied.*

i. *Each $I_t$ is a neutrosophic subsemigroup of $S_t$, t = 1, 2, ..., N i.e. $I_t$ is a neutrosophic strong N-subsemigroup of S(N).*
ii. *Each $I_t$ is a two sided ideal of $S_t$ for t = 1, 2, ..., N.*

*Similarly one can define neutrosophic strong N-left ideal or neutrosophic strong right ideal of S(N).*

*A neutrosophic strong N-ideal is one which is both a neutrosophic strong N-left ideal and N-right ideal of S(N).*

**DEFINITION 1.2.35:** *Let $S(N) = \{S_1 \cup S_2 \cup ... \cup S_N, *_1, ..., *_N\}$ be a neutrosophic N-semigroup. A proper subset $P = \{P_1 \cup P_2 \cup ... \cup P_N, *_1, ..., *_N\}$ of S (N) is said to be a neutrosophic N-subsemigroup, if the following conditions are true*

i. *P is a neutrosophic sub N-semigroup of S(N).*
ii. *Each $P_i = P \cap S_i$, i =1, 2, ..., N is an ideal of $S_i$.*

*Then P is called or defined as the neutrosophic N ideal of the neutrosophic N-semigroup S(N).*

**DEFINITION 1.2.36:** *Let $S(N) = \{S_1 \cup S_2 \cup ... \cup S_N, *_1, ..., *_N\}$ be a neutrosophic strong N-semigroup. Let $J = \{I_1 \cup I_2 \cup ... \cup I_N, *_1, ..., *_N\}$ be a proper subset of S(N) which is a neutrosophic strong N-ideal of S(N). J is said to be a neutrosophic strong maximal N-ideal of S(N) if each $I_t \subset S_t$, (t = 1, 2, ..., N) is a maximal ideal of $S_t$.*

*It may so happen that at times only some of the ideals $I_t$ in $S_t$ may be maximal and some may not be in that case we call the ideal J to be a neutrosophic quasi maximal N-ideal of S(N). Suppose $S(N) = \{S_1 \cup S_2 \cup ... \cup S_N, *_1, ..., *_N\}$ is a neutrosophic strong N-semigroup, $J' = \{J_1 \cup J_2 \cup ... \cup J_N, *_1, ..., *_N\}$ be a neutrosophic strong N-ideal of S(N).*

*J' is said to be a neutrosophic strong minimal N-ideal of S (N) if each $J_i \subset S_i$ is a minimal ideal of $S_i$ for i = 1, 2, ..., N. It may so happen that some of the ideals $J_i \subset S_i$ be minimal and some may not be minimal in this case we call J' the neutrosophic strong quasi minimal N-ideal of S(N).*

## 1.3 S-Loops and S-N-loops

We at this juncture like to express that books solely on loops are meager or absent as, R.H.Bruck deals with loops on his book "*A Survey of Binary Systems*", that too published as early as 1958, [3]. Other two books are on "*Quasigroups and Loops*" one by H.O. Pflugfelder, 1990 which is introductory and the other book co-edited by Orin Chein, H.O. Pflugfelder and J.D. Smith in 1990 [25].

So we felt it important to recall almost all the properties and definitions related with loops [3, 47]. We just recall a few of the properties about loops which will make this book a self contained one. In this section we introduce briefly the notion of S-loop, S-N-loops and neutrosophic loops and their generalization.

**DEFINITION 1.3.1:** *A non-empty set L is said to form a loop, if on L is defined a binary operation called the product denoted by '$\bullet$' such that*

i. *For all a, b $\in$ L we have $a \bullet b \in L$ (closure property).*
ii. *There exists an element e $\in$ L such that $a \bullet e = e \bullet a = a$ for all a $\in$ L (e is called the identity element of L).*
iii. *For every ordered pair (a, b) $\in L \times L$ there exists a unique pair (x, y) in L such that ax = b and ya = b.*

**DEFINITION 1.3.2:** *Let L be a loop. A non-empty subset H of L is called a subloop of L if H itself is a loop under the operation of L.*

**DEFINITION 1.3.3:** *Let L be a loop. A subloop H of L is said to be a normal subloop of L, if*

i. *xH = Hx.*
ii. *(Hx)y = H(xy).*
iii. *y(xH) = (yx)H*

*for all x, y $\in$ L.*

**DEFINITION 1.3.4:** *A loop L is said to be a simple loop if it does not contain any non-trivial normal subloop.*

**DEFINITION 1.3.5:** *The commutator subloop of a loop L denoted by L' is the subloop generated by all of its commutators, that is,* $\langle\{x \in L \,/\, x = (y, z)$ *for some* $y, z \in L\}\rangle$ *where for* $A \subseteq L$, $\langle A\rangle$ *denotes the subloop generated by A.*

**DEFINITION 1.3.6:** *If x, y and z are elements of a loop L an associator (x, y, z) is defined by, (xy)z = (x(yz)) (x, y, z).*

**DEFINITION 1.3.7:** *The associator subloop of a loop L (denoted by A(L)) is the subloop generated by all of its associators, that is* $\langle\{x \in L \,/\, x = (a, b, c)$ *for some* $a, b, c \in L\}\rangle$.

**DEFINITION 1.3.8:** *The centre Z(L) of a loop L is the intersection of the nucleus and the Moufang centre, that is Z(L) =* $C(L) \cap N(L)$.

**DEFINITION [35]**: *A normal subloop of a loop L is any subloop of L which is the kernel of some homomorphism from L to a loop.*

Further Pflugfelder [25] has proved the central subgroup Z(L) of a loop L is normal in L.

**DEFINITION [35]**: *Let L be a loop. The centrally derived subloop (or normal commutator- associator subloop) of L is defined to be the smallest normal subloop* $L' \subset L$ *such that* $L / L'$ *is an abelian group.*

*Similarly nuclearly derived subloop (or normal associator subloop) of L is defined to be the smallest normal subloop* $L_1 \subset L$ *such that* $L / L_1$ *is a group. Bruck proves L' and* $L_1$ *are well defined.*

**DEFINITION [35]**: *The Frattini subloop* $\phi(L)$ *of a loop L is defined to be the set of all non-generators of L, that is the set of all* $x \in L$ *such that for any subset S of L,* $L = \langle x, S\rangle$ *implies L =*

*⟨S⟩. Bruck has proved as stated by Tim Hsu $\phi(L) \subset L$ and $L / \phi(L)$ is isomorphic to a subgroup of the direct product of groups of prime order.*

**DEFINITION [22]**: *Let L be a loop. The commutant of L is the set (L) = $\{a \in L / ax = xa \ \forall x \in L\}$. The centre of L is the set of all $a \in C(L)$ such that $a \bullet xy = ax \bullet y = x \bullet ay = xa \bullet y$ and $xy \bullet a = x \bullet ya$ for all $x, y \in L$. The centre is a normal subloop. The commutant is also known as Moufang Centre in literature.*

**DEFINITION [23]:** *A left loop (B, $\bullet$) is a set B together with a binary operation '$\bullet$' such that (i) for each $a \in B$, the mapping $x \to a \bullet x$ is a bijection and (ii) there exists a two sided identity $1 \in B$ satisfying $1 \bullet x = x \bullet 1 = x$ for every $x \in B$. A right loop is defined similarly. A loop is both a right loop and a left loop.*

**DEFINITION [11]** : *A loop L is said to have the weak Lagrange property if, for each subloop K of L, $|K|$ divides $|L|$. It has the strong Lagrange property if every subloop K of L has the weak Lagrange property.*

**DEFINITION 1.3.9:** *A loop L is said to be power-associative in the sense that every element of L generates an abelian group.*

**DEFINITION 1.3.10:** *A loop L is diassociative loop if every pair of elements of L generates a subgroup.*

**DEFINITION 1.3.11:** *A loop L is said to be a Moufang loop if it satisfies any one of the following identities:*

- i. $(xy)(zx) = (x(yz))x$
- ii. $((xy)z)y = x(y(zy))$
- iii. $x(y(xz) = ((xy)x)z$

*for all $x, y, z \in L$.*

**DEFINITION 1.3.12:** *Let L be a loop, L is called a Bruck loop if $x(yx)z = x(y(xz))$ and $(xy)^{-1} = x^{-1}y^{-1}$ for all $x, y, z \in L$.*

**DEFINITION 1.3.13:** *A loop (L, •) is called a Bol loop if ((xy)z)y = x((yz)y) for all x, y, z ∈ L.*

**DEFINITION 1.3.14:** *A loop L is said to be right alternative if (xy)y = x(yy) for all x, y ∈ L and L is left alternative if (xx)y = x(xy) for all x, y ∈ L. L is said to be an alternative loop if it is both a right and left alternative loop.*

**DEFINITION 1.3.15:** *A loop (L, •) is called a weak inverse property loop (WIP-loop) if (xy)z = e imply x(yz) = e for all x, y, z ∈ L.*

**DEFINITION 1.3.16:** *A loop L is said to be semi alternative if (x, y, z) = (y, z, x) for all x, y, z ∈ L, where (x, y, z) denotes the associator of elements x, y, z ∈ L.*

**THEOREM (MOUFANG'S THEOREM):** *Every Moufang loop G is diassociative more generally, if a, b, c are elements in a Moufang loop G such that (ab)c = a(bc) then a, b, c generate an associative loop.*

The proof is left for the reader; for assistance refer Bruck R.H. [3].

**DEFINITION 1.3.17:** *Let L be a loop, L is said to be a two unique product loop (t.u.p) if given any two non-empty finite subsets A and B of L with |A| + |B| > 2 there exist at least two distinct elements x and y of L that have unique representation in the from x = ab and y = cd with a, c ∈ A and b, d ∈ B.*

*A loop L is called a unique product (u.p) loop if, given A and B two non-empty finite subsets of L, then there always exists at least one x ∈ L which has a unique representation in the from x = ab, with a ∈ A and b ∈ B.*

**DEFINITION 1.3.18:** *Let (L, •) be a loop. The principal isotope (L, ∗) of (L, •) with respect to any predetermined a, b ∈ L is defined by x ∗ y = XY, for all x, y ∈ L, where Xa = x and bY = y for some X, Y ∈ L.*

**DEFINITION 1.3.19:** *Let L be a loop, L is said to be a G-loop if it is isomorphic to all of its principal isotopes.*

The main objective of this section is the introduction of a new class of loops with a natural simple operation. As to introduce loops several functions or maps are defined satisfying some desired conditions we felt that it would be nice if we can get a natural class of loops built using integers.

Here we define the new class of loops of any even order, they are distinctly different from the loops constructed by other researchers. Here we enumerate several of the properties enjoyed by these loops.

**DEFINITION [41]:** *Let $L_n(m) = \{e, 1, 2, \ldots, n\}$ be a set where $n > 3$, n is odd and m is a positive integer such that $(m, n) = 1$ and $(m-1, n) = 1$ with $m < n$.*
*Define on $L_n(m)$ a binary operation '$\bullet$' as follows:*

i. $e \bullet i = i \bullet e = i$ *for all* $i \in L_n(m)$
ii. $i^2 = i \bullet i = e$ *for all* $i \in L_n(m)$
iii. $i \bullet j = t$ *where* $t = (mj - (m - 1)i) \pmod n$

*for all $i, j \in L_n(m)$; $i \neq j$, $i \neq e$ and $j \neq e$, then $L_n(m)$ is a loop under the binary operation '$\bullet$'.*

***Example 1.3.1:*** Consider the loop $L_5(2) = \{e, 1, 2, 3, 4, 5\}$. The composition table for $L_5(2)$ is given below:

| • | e | 1 | 2 | 3 | 4 | 5 |
|---|---|---|---|---|---|---|
| e | e | 1 | 2 | 3 | 4 | 5 |
| 1 | 1 | e | 3 | 5 | 2 | 4 |
| 2 | 2 | 5 | e | 4 | 1 | 3 |
| 3 | 3 | 4 | 1 | e | 5 | 2 |
| 4 | 4 | 3 | 5 | 2 | e | 1 |
| 5 | 5 | 2 | 4 | 1 | 3 | e |

This loop is of order 6 which is both non-associative and non-commutative.

*Physical interpretation of the operation in the loop $L_n(m)$*:

We give a physical interpretation of this class of loops as follows: Let $L_n(m)= \{e, 1, 2, \ldots , n\}$ be a loop in this identity element of the loop are equidistantly placed on a circle with e as its centre.

We assume the elements to move always in the clockwise direction.

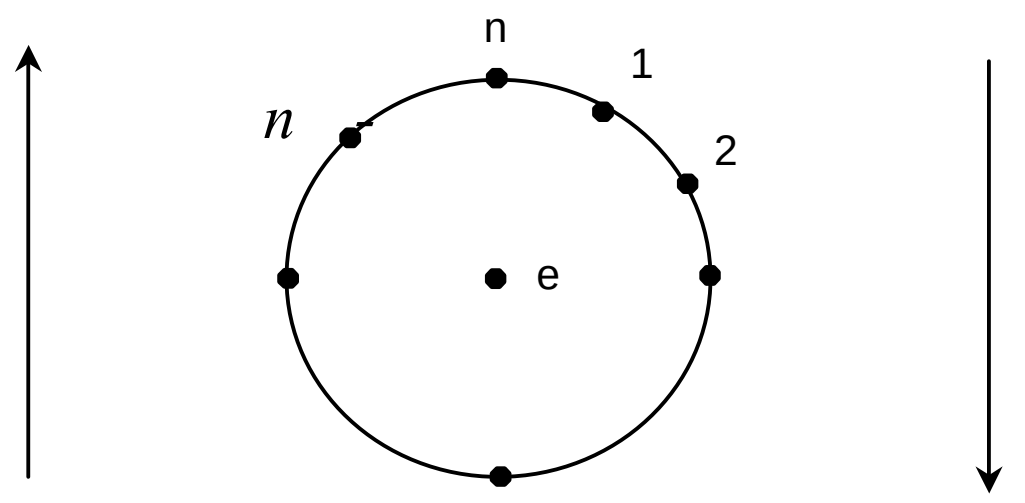

Let $i, j \in L_n(m)$ $(i \neq j, i \neq e, j \neq e)$. If j is the $r^{th}$ element from i counting in the clockwise direction the $i \bullet j$ will be the $t^{th}$ element from j in the clockwise direction where $t = (m-1)r$. We see that in general $i \bullet j$ need not be equal to $j \bullet i$. When $i = j$ we define $i^2 = e$ and $i \bullet e = e \bullet i = i$ for all $i \in L_n(m)$ and e acts as the identity in $L_n(m)$.

***Example 1.3.2***: Now the loop $L_7(4)$ is given by the following table:

| • | e | 1 | 2 | 3 | 4 | 5 | 6 | 7 |
|---|---|---|---|---|---|---|---|---|
| e | e | 1 | 2 | 3 | 4 | 5 | 6 | 7 |
| 1 | 1 | e | 5 | 2 | 6 | 3 | 7 | 4 |
| 2 | 2 | 5 | e | 6 | 3 | 7 | 4 | 1 |
| 3 | 3 | 2 | 6 | e | 7 | 4 | 1 | 5 |
| 4 | 4 | 6 | 3 | 7 | e | 1 | 5 | 2 |
| 5 | 5 | 3 | 7 | 4 | 1 | e | 2 | 6 |
| 6 | 6 | 7 | 4 | 1 | 5 | 2 | e | 3 |
| 7 | 7 | 4 | 1 | 5 | 2 | 6 | 3 | e |

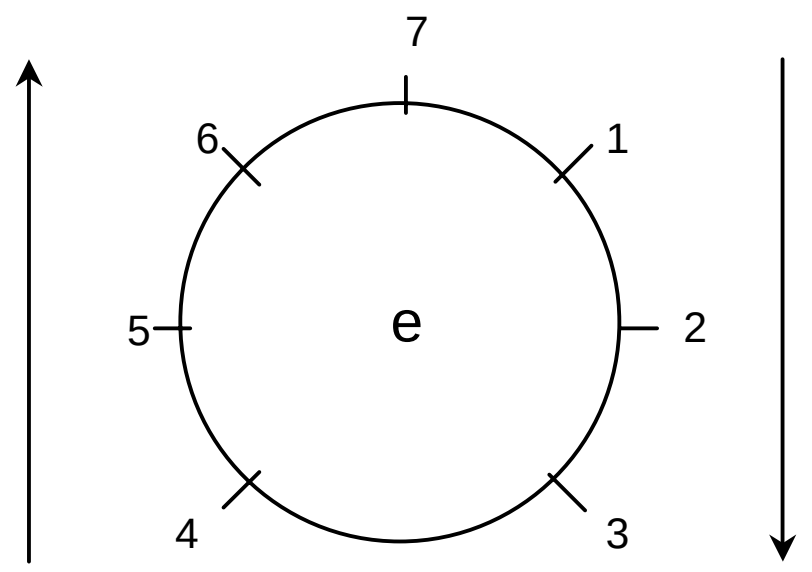

Let 2, 4 ∈ $L_7(4)$. 4 is the 2nd element from 2 in the clockwise direction. So 2.4 will be (4 – 1)2 that is the 6th element from 4 in the clockwise direction which is 3. Hence 2.4 = 3.

**Notation**: Let $L_n$ denote the class of loops. $L_n(m)$ for fixed n and various m's satisfying the conditions m < n, (m, n) = 1 and (m – 1, n) = 1, that is $L_n = \{L_n(m) \mid n > 3, n \text{ odd}, m < n, (m, n) = 1 \text{ and } (m-1, n) = 1\}$.

***Example 1.3.3:*** Let n = 5. The class $L_5$ contains three loops; viz. $L_5(2)$, $L_5(3)$ and $L_5(4)$ given by the following tables:

$L_5(2)$

| • | e | 1 | 2 | 3 | 4 | 5 |
|---|---|---|---|---|---|---|
| e | e | 1 | 2 | 3 | 4 | 5 |
| 1 | 1 | e | 3 | 5 | 2 | 4 |
| 2 | 2 | 5 | e | 4 | 1 | 3 |
| 3 | 3 | 4 | 1 | e | 5 | 2 |
| 4 | 4 | 3 | 5 | 2 | e | 1 |
| 5 | 5 | 2 | 4 | 1 | 3 | e |

$L_5(3)$

| • | e | 1 | 2 | 3 | 4 | 5 |
|---|---|---|---|---|---|---|
| e | e | 1 | 2 | 3 | 4 | 5 |
| 1 | 1 | e | 4 | 2 | 5 | 3 |
| 2 | 2 | 4 | e | 5 | 3 | 1 |
| 3 | 3 | 2 | 5 | e | 1 | 4 |
| 4 | 4 | 5 | 3 | 1 | e | 2 |
| 5 | 5 | 3 | 1 | 4 | 2 | e |

$L_5(4)$

| • | e | 1 | 2 | 3 | 4 | 5 |
|---|---|---|---|---|---|---|
| e | e | 1 | 2 | 3 | 4 | 5 |
| 1 | 1 | e | 5 | 4 | 3 | 2 |
| 2 | 2 | 3 | e | 1 | 5 | 4 |
| 3 | 3 | 5 | 4 | e | 2 | 1 |
| 4 | 4 | 2 | 1 | 5 | e | 3 |
| 5 | 5 | 4 | 3 | 2 | 1 | e |

**THEOREM [27]:** *Let $L_n$ be the class of loops for any $n > 3$, if $n = p_1^{\alpha_1} p_2^{\alpha_2} \ldots p_k^{\alpha_k}$ ($\alpha_i > 1$, for $i = 1, 2, \ldots, k$), then $|L_n| = \prod_{i=1}^{k} (p_i - 2)\, p_i^{\alpha_i - 1}$ where $|L_n|$ denotes the number of loops in $L_n$.*

The proof is left for the reader as an exercise.

**THEOREM [27]**: *$L_n$ contains one and only one commutative loop. This happens when $m = (n + 1) / 2$. Clearly for this m, we have $(m, n) = 1$ and $(m – 1, n) = 1$.*

It can be easily verified by using simple number theoretic techniques.

**THEOREM [27]**: *Let $L_n$ be the class of loops. If $n = p_1^{\alpha_1} p_2^{\alpha_2} \ldots p_k^{\alpha_k}$, then $L_n$ contains exactly $F_n$ loops which are strictly non-commutative where $F_n = \prod_{i=1}^{k} (p_i - 3)\, p_i^{\alpha_i - 1}$.*

The proof is left for the reader as an exercise.

Note: If n = p where p is a prime greater than or equal to 5 then in $L_n$ a loop is either commutative or strictly non-commutative. Further it is interesting to note if n = 3t then the class $L_n$ does not contain any strictly non-commutative loop.

**THEOREM [32]**: *The class of loops $L_n$ contains exactly one left alternative loop and one right alternative loop but does not contain any alternative loop.*

*Proof*: We see $L_n(2)$ is the only right alternative loop that is when m = 2 (Left for the reader to prove using simple number theoretic techniques). When m = n –1 that is $L_n(n-1)$ is the only left alternative loop in the class of loops $L_n$.

From this it is impossible to find a loop in $L_n$, which is simultaneously right alternative and left alternative. Further it is clear from earlier result both the right alternative loop and the left alternative loop is not commutative.

**THEOREM [27]:** *Let $L_n$ be the class of loops:*

i. *$L_n$ does not contain any Moufang loop.*
ii. *$L_n$ does not contain any Bol loop.*
iii. *$L_n$ does not contain any Bruck loop.*

The reader is requested to prove these results using number theoretic techniques.

**THEOREM [41]**: *Let $L_n(m) \in L_n$. Then $L_n(m)$ is a weak inverse property (WIP) loop if and only if $(m^2 - m + 1) \equiv 0 (mod\ n)$.*

*Proof*: It is easily checked that for a loop to be a WIP loop we have "if (xy)z = e then x(yz) = e where x, y, z ∈ L." Both way conditions can be derived using the defining operation on the loop $L_n(m)$.

***Example 1.3.4:*** L be the loop $L_7(3)$ = {e, 1, 2, 3, 4, 5, 6, 7} be in $L_7$ given by the following table:

| • | e | 1 | 2 | 3 | 4 | 5 | 6 | 7 |
|---|---|---|---|---|---|---|---|---|
| e | e | 1 | 2 | 3 | 4 | 5 | 6 | 7 |
| 1 | 1 | e | 4 | 7 | 3 | 6 | 2 | 5 |
| 2 | 2 | 6 | e | 5 | 1 | 4 | 7 | 3 |
| 3 | 3 | 4 | 7 | e | 6 | 2 | 5 | 1 |
| 4 | 4 | 2 | 5 | 1 | e | 7 | 3 | 6 |
| 5 | 5 | 7 | 3 | 6 | 2 | e | 1 | 4 |
| 6 | 6 | 5 | 1 | 4 | 7 | 3 | e | 2 |
| 7 | 7 | 3 | 6 | 2 | 5 | 1 | 4 | e |

It is easily verified $L_7(3)$ is a WIP loop. One way is easy for $(m^2 - m + 1) \equiv 0 \pmod 7$ that is $9 - 3 + 1 = 9 + 4 + 1 \equiv 0 (\mathrm{mod}\ 7)$. It is interesting to note that no loop in the class $L_n$ contain any associative loop.

**THEOREM [27]**: *Let $L_n$ be the class of loops. The number of strictly non-right (left) alternative loops is $P_n$ where $P_n = \prod_{i=1}^{k} (p_i - 3) p_i^{\alpha_i - 1}$ and $n = \prod_{i=1}^{k} p_i^{\alpha_i}$.*

The proof is left for the reader to verify.

Now we proceed on to study the associator and the commutator of the loops in $L_n$.

**THEOREM [27]**: *Let $L_n(m) \in L_n$. The associator $A(L_n(m)) = L_n(m)$.*

For more literature about the new class of loops refer [41, 47].

**DEFINITION 1.3.20:** *Let $(L, *_1, \ldots, *_N)$ be a non empty set with N binary operations $*_i$. L is said to be a N loop if L satisfies the following conditions:*

- *i.* *$L = L_1 \cup L_2 \cup \ldots \cup L_N$ where each $L_i$ is a proper subset of L; i.e., $L_i \not\subseteq L_j$ or $L_j \not\subseteq L_i$ if $i \neq j$ for $1 \leq i, j \leq N$.*
- *ii.* *$(L_i, *_i)$ is a loop for some i, $1 \leq i \leq N$.*
- *iii.* *$(L_j, *_j)$ is a loop or a group for some j, $1 \leq j \leq N$.*

*For a N-loop we demand atleast one $(L_j, *_j)$ to be a loop.*

**DEFINITION 1.3.21:** *Let $(L = L_1 \cup L_2 \cup \ldots \cup L_N, *_1, \ldots, *_N)$ be a N-loop. L is said to be a commutative N-loop if each $(L_i, *_i)$ is commutative, i = 1, 2, …, N. We say L is inner commutative if each of its proper subset which is N-loop under the binary operations of L are commutative.*

**DEFINITION 1.3.22:** *Let* $L = \{L_1 \cup L_2 \cup \ldots \cup L_N, *_1, *_2, \ldots, *_N\}$ *be a N-loop. We say L is a Moufang N-loop if all the loops* $(L_i, *_i)$ *satisfy the following identities.*

i. $(xy)(zx) = (x(yz))x$
ii. $((xy)z)y = x(y(zy))$
iii. $x(y(xz)) = ((xy)x)z$

*for all* $x, y, z \in L_i$, $i = 1, 2, \ldots, N$.

Now we proceed on to define a Bruck N-loop.

**DEFINITION 1.3.23:** *Let* $L = (L_1 \cup L_2 \cup \ldots \cup L_N, *_1, \ldots, *_N)$ *be a N-loop. We call L a Bruck N-loop if all the loops* $(L_i, *_i)$ *satisfy the identities*

i. $(x(yz))z = x(y(xz))$
ii. $(xy)^{-1} = x^{-1}y^{-1}$

*for all* $x, y \in L_i$, $i = 1, 2, \ldots, N$.

**DEFINITION 1.3.24:** *Let* $L = (L_1 \cup L_2 \cup \ldots \cup L_N, *_1, \ldots, *_N)$ *be a N-loop. A non empty subset P of L is said to be a sub N-loop, if P is a N-loop under the operations of L i.e.,* $P = \{P_1 \cup P_2 \cup P_3 \cup \ldots \cup P_N, *_1, \ldots, *_N\}$ *with each* $\{P_i, *_i\}$ *is a loop or a group.*

**DEFINITION 1.3.25:** *Let* $L = \{L_1 \cup L_2 \cup \ldots \cup L_N, *_1, \ldots, *_N\}$ *be a N-loop. A proper subset P of L is said to be a normal sub N-loop of L if*

i. *P is a sub N-loop of L.*
ii. $x_i P_i = P_i x_i$ *for all* $x_i \in L_i$.
iii. $y_i (x_i P_i) = (y_i x_i) P_i$ *for all* $x_i\, y_i \in L_i$.

*A N-loop is said to be a simple N-loop if L has no proper normal sub N-loop.*

Now we proceed on to define the notion of Moufang center.

**DEFINITION 1.3.26:** *Let $L = \{L_1 \cup L_2 \cup \ldots \cup L_N, *_1, \ldots, *_N\}$ be a N-loop. We say $C_N(L)$ is the Moufang N-centre of this N-loop if $C_N(L) = C_1(L_1) \cup C_2(L_2) \cup \ldots \cup C_N(L_N)$ where $C_i(L_i) = \{x_i \in L_i / x_i y_i = y_i x_i$ for all $y_i \in L_i\}$, $i = 1, 2, \ldots, N$.*

**DEFINITION 1.3.27:** *Let L and P to two N-loops i.e. $L = \{L_1 \cup L_2 \cup \ldots \cup L_N, *_1, \ldots, *_N\}$ and $P = \{P_1 \cup P_2 \cup \ldots \cup P_N, o_1, \ldots, o_N\}$. We say a map $\theta : L \to P$ is a N-loop homomorphism if $\theta = \theta_1 \cup \theta_2 \cup \ldots \cup \theta_N$ '$\cup$' is just a symbol and $\theta_i$ is a loop homomorphism from $L_i$ to $P_i$ for each $i = 1, 2, \ldots, N$.*

**DEFINITION 1.3.28:** *Let $L = \{L_1 \cup L_2 \cup \ldots \cup L_N, *_1, \ldots, *_N\}$ be a N-loop. We say L is weak Moufang N-loop if there exists atleast a loop $(L_i, *_i)$ such that $L_i$ is a Moufang loop.*

*Note:* $L_i$ should not be a group it should be only a loop.

**DEFINITION 1.3.29:** *Let $L = \{L_1 \cup L_2 \cup \ldots \cup L_N, *_1, \ldots, *_N\}$ be a N-loop. If x and $y \in L$ are elements of $L_i$ the N-commutator (x, y) is defined as $xy = (yx)(x, y)$, $1 \leq i \leq N$.*

**DEFINITION 1.3.30:** *Let $L = \{L_1 \cup L_2 \cup \ldots \cup L_N, *_1, \ldots, *_N\}$ be a N-loop. If x, y, z are elements of the N-loop L, an associator (x, y, z) is defined only if $x, y, z \in L_i$ for some i $(1 \leq i \leq N)$ and is defined to be $(xy)z = (x(yz))(x, y, z)$.*

**DEFINITION 1.3.31:** *Let $L = \{L_1 \cup L_2 \cup \ldots \cup L_N, *_1, *_2, \ldots, *_N\}$ be a N-loop of finite order. For $\alpha_i \in L_i$ define $R_{\alpha_i}$ as a permutation of the loop $L_i$, $R_{\alpha_i} : x_i \to x_i\alpha_i$. This is true for $i = 1, 2, \ldots, N$ we define the set*

$$\left\{R_{\alpha_1} \cup R_{\alpha_2} \cup \ldots \cup R_{\alpha_N} \,\middle|\, \alpha_i \in L_i; i = 1, 2, \ldots, N\right\}$$

*as the right regular N-representation of the N loop L.*

**DEFINITION 1.3.32:** *Let $L = \{L_1 \cup L_2 \cup \ldots \cup L_N, *_1, \ldots, *_N\}$ be a N-loop. For any pre determined pair $a_i, b_i \in L_i$, $i \in \{1, 2, \ldots, N\}$ a principal isotope $(L, o_1, \ldots, o_N)$, of the N loop L is defined by $x_i \, o_i \, y_i = X_i *_i Y_i$ where $X_i + a_i = x_i$ and $b_i + Y_i = y_i$, $i = 1, 2, \ldots,$*

*N. L is called G-N-loop if it is isomorphic to all of its principal isotopes.*

**DEFINITION 1.3.33:** *The Smarandache loop (S-loop) is defined to be a loop L such that a proper subset A of L is a subgroup (with respect to the same induced operation) that is $\phi \neq A \subset L$.*

**DEFINITION 1.3.34:** *Let L be a loop. A proper subset A of L is said to be a Smarandache subloop (S-subloop) of L if A is a subloop of L and A is itself a S-loop; that is A contains a proper subset B contained in A such that B is a group under the operations of L. We demand A to be a S-subloop which is not a subgroup.*

**DEFINITION 1.3.35:** *Let L be a loop. We say a non-empty subset A of L is a Smarandache normal subloop (S-normal subloop) of L if*

i. *A is itself a normal subloop of L.*
ii. *A contains a proper subset B where B is a subgroup under the operations of L. If L has no S-normal subloop we say the loop L is Smarandache simple (S-simple).*

**DEFINITION 1.3.36:** *Let L and L' be two Smarandache loops with A and A' its subgroups respectively. A map $\phi$ from L to L' is called Smarandache loop homomorphism (S-loop homomorphism) if $\phi$ restricted to A is mapped onto a subgroup A' of L'; that is $\phi : A \to A'$ is a group homomorphism. The concept of Smarandache loop isomorphism and automorphism are defined in a similar way. It is important to observe the following facts:*

i. *The map $\phi$ from L to L' need not be even be a loop homomorphism.*
ii. *Two loops of different order can be isomorphic.*
iii. *Further two loops which are isomorphic for a map $\phi$ may not be isomorphic for some other map $\eta$.*

iv. *If L and L' have at least subgroups A and A' in L and L' respectively which are isomorphic then certainly L and L' are isomorphic.*

**DEFINITION 1.3.37:** *Let L be a loop. If A (A proper subset of L) is a S-subloop of L is such that A is a pseudo commutative loop then we say L is a Smarandache pseudo commutative loop (S-pseudo commutative loop) i.e. for every a, b $\in$ A we have an x $\in$ B such that a(xb) = b(xa) (or (bx)a), B is a subgroup in A. Thus we see for a loop to be a S-pseudo commutative we need not have the loop L to be a pseudo-commutative loop. If L is itself a pseudo commutative loop then trivially L is a S-pseudo commutative loop.*

**DEFINITION 1.3.38:** *Let L be a loop. We say L is a Smarandache associative loop (S-associative loop) if L has a S-subloop A such that A contains a triple x, y, z (all three elements distinct and different from e, the identity element of A) such that $x \bullet (y \bullet z) = (x \bullet y) \bullet z$. This triple is called the S-associative triple.*

**DEFINITION 1.3.39:** *Let L be a loop the Smarandache left nucleus (S-left nucleus) $S(N_\lambda) = \{a \in A / (a, x, y) = e$ for all x, y $\in A\}$ is a subloop of A, where A is a S-subloop of L. Thus we see in case of S-left nucleus we can have many subloops; it is not unique as in case of loops. If L has no S-subloops but L is a S-loop then $S(N_\lambda) = S N_\lambda = N_\lambda$.*

**DEFINITION 1.3.40:** *Let L be a loop, the Smarandache right nucleus (S-right nucleus) $S(N_\rho) = \{ a \in A / (x, y, a) = e$ for all x, y $\in A\}$ is a subloop of L where A is a S-subloop of L. If L has no S-subloops but L is a S-loop then $S(N_\rho) = SN_\rho = N_\rho$.*

**DEFINITION 1.3.41:** *Let L be a loop, the Smarandache middle nucleus (S-middle nucleus). $S(N_\mu) = \{ a \in A / ( x, a, y) = e$ for all x, y $\in A\}$ is a subloop of L where A is a S-subloop of L. If L has no S-subloop but L is a S-loop then $S(N_\mu) = SN_\mu = N_\mu$.*

**DEFINITION 1.3.42:** *The Smarandache nucleus S(N(L)) of a loop L is a subloop given by* $S(N(L)) = SN_\mu \cap SN_\lambda \cap SN_\rho$*. It is important to note that unlike the nucleus we can have several Smarandache nucleus depending on the S-subloops. If L has no S-subloops but L is S-loop then S(N(L)) = N(L).*

**DEFINITION 1.3.43:** *Let L be a loop. The Smarandache Moufang centre (S-Moufang centre) SC(L) is the set of elements in a S-subloop A of L which commute with every element of A, that is* $SC(L) = \{x \in A / xy = yx \text{ for all } y \in A\}$*. If L has no S-subloops but L is a S-loop then we have SC(L) = C(L).If L has many S-subloops then we have correspondingly many S-Moufang centres SC(L).*

**DEFINITION 1.3.44:** *Let L be a loop, the Smarandache centre (S-centre) (SZ(L)) of a loop. L is the intersection of* $SZ(L) = SC(L) \cap S(N(L))$ *for a S-subloop A of L.*

**DEFINITION 1.3.45:** *Let* $L = L_1 \times S_n$ *be the direct product of a loop and a group. We call this the Smarandache mixed direct product of loops (S-mixed direct product of loops). We insist that L is a loop and one of* $L_1$ *or* $S_n$ *is group. We do not take all groups in the S-mixed direct product or all loops in this product. Further we can extend this to any number of loops and groups that is if* $L = L_1 \times L_2 \times \ldots \times L_n$ *where some of the* $L_i$*'s are loops and some of them are groups. Then L is called the S-mixed direct product of n-terms.*

**DEFINITION 1.3.46:** *Let L be a S-loop. We define Smarandache right cosets (S-right cosets) in L as follows:*

*Let* $A \subset L$ *be the subgroup of L and for* $m \in L$ *we have* $Am = \{am/a \in A\}$*. Am is called a S-right coset of A in L.*

*Similarly we can for any subgroup A of L define Smarandache left coset (S-left coset) for some* $m \in L$ *as* $mA = \{ma / a \in A\}$*. If the loop L is commutative then only we have mA = Am. Even if L is S-commutative or S-strongly commutative still we need not have Am = mA for all* $m \in L$*.*

**DEFINITION 1.3.47:** *Let L = $\{L_1 \cup L_2 \cup \ldots \cup L_N, *_1, \ldots, *_N\}$ be a N-loop. We call L a Smarandache N-loop (S-N-loop) if L has a proper subset P; (P = $P_1 \cup P_2 \cup \ldots \cup P_N, *_1, \ldots, *_N$ ) where P is a N-group.*

**DEFINITION 1.3.48:** *Let L = $\{L_1 \cup L_2 \cup \ldots \cup L_N, *_1, \ldots, *_N\}$ be a N-loop. We call a proper subset P of L where P = $\{P_1 \cup P_2 \cup \ldots \cup P_N, *_1, \ldots, *_N\}$ to be a Smarandache sub N-loop (S-sub N-loop) of L if P itself is a S-N-loop.*

**DEFINITION 1.3.49**: *Let L = $\{L_1 \cup L_2 \cup \ldots \cup L_N, *_1, \ldots, *_N\}$ be a N-loop. L is said to be a Smarandache commutative N-loop (S-commutative N-loop) if every proper subset P of L which are N-groups are commutative N-groups of L.*

*If the N-loop L has atleast one proper subset P which is a commutative N-group then we call L a Smarandache weakly commutative N-loop (S-weakly commutative loop).*

**DEFINITION 1.3.50:** *Let L = $\{L_1 \cup L_2 \cup \ldots \cup L_N, *_1, \ldots, *_N\}$ be a N-loop. We say the N-loop L is a Smarandache cyclic N-loop (S-cyclic N-loop) if every proper subset which is a N-group is a cyclic N-group.*

**DEFINITION 1.3.51:** *Let L = $\{L_1 \cup L_2 \cup \ldots \cup L_N, *_1, \ldots, *_N\}$ be a N-loop. A proper S-sub N-loop (P = $P_1 \cup P_2 \cup \ldots \cup P_N, *_1, \ldots, *_N$) of L is said to be a Smarandache normal N-loop (S-normal N-loop) if*

i. $x_i P_i = P_i x_i$
ii. $P_i x_i(y_i) = P_i (x_i y_i)$
iii. $y_i (x_i P_i) = (y_i x_i) P_i$

*for all $x_i, y_i \in P_i$ for i = 1, 2, ..., N.*

*If the N-loop L has no proper S-normal sub N-loop then we call the N-loop to be Smarandache simple N-loop (S-simple N-loop).*

**DEFINITION 1.3.52:** *Let L = $\{L_1 \cup L_2 \cup \ldots \cup L_N, *_1, \ldots, *_N\}$ be a N-loop and P = $(P_1 \cup P_2 \cup \ldots \cup P_N, *_1, \ldots, *_N)$ be a S-sub N-*

*loop of L. ($P \subset L$) for any N-pair of elements $x_i$, $y_i \in P_i$ ($i \in$ {1, 2, …, N}) the N-commutator ($x_i$, $y_i$) is defined by $x_i y_i = (y_i x_i)(x_i, y_i)$.*

*The Smarandache commutator sub N-loop (S-commutator sub N-loop) of L denoted by S ($L^s$) is the S-sub N-loop generated by all its commutators.*

**DEFINITION 1.3.53:** *Let $L = \{L_1 \cup L_2 \cup \ldots \cup L_N, *_1, \ldots, *_N\}$ be a N-loop. Let $P \subset L$ be a S-sub N-loop of L. If x, y, z are elements of the S-sub N-loop P of L, an associator (x, y, z) is defined by $(xy)z = (x (y z)) (x, y, z)$.*

*The associator S-sub N-loop of the S-sub N-loop P of L denoted by $A\left(L_N^S\right)$ is the S-sub N-loop generated by all the associators, that is $\langle \{x \in P \mid x = (a, b, c)$ for some $a, b, c \in P\}\rangle$. If $A\left(L_N^S\right)$ happens to be a S-sub N-loop then only we call $A\left(L_N^S\right)$ the Smarandache associator (S-associator) of L related to the S-sub N-loop P.*

**DEFINITION 1.3.54:** *Let $L = \{L_1 \cup L_2 \cup \ldots \cup L_N, *_1, \ldots, *_N\}$ be a N-loop. If for $a, b \in L$ with $ab = ba$ in L, we have $(ax_i) b = (bx_i) a$ (or $b (x_i a)$) for all $x_i \in L_i$ (if $a, b \in L_i$), then we say the pair is pseudo commutative.*

*If every commutative pair in every $L_i$ is pseudo commutative then we call L a pseudo commutative N-loop. If we have for every pair of commuting elements a, b in a proper subset P of the N-loop L (which is a S-sub N-loop). If a, b is a pseudo commutative pair for every x in P then we call the pair a Smarandache pseudo commutative pair (S-pseudo commutative pair).*

**DEFINITION 1.3.55:** *Let $L = \{L_1 \cup L_2 \cup \ldots \cup L_N, *_1, \ldots, *_N\}$ be a N-loop. The pseudo commutator of L denoted by $P(L) = \langle\{p \in L$ / $a(xb) = p ([bx] a)\}\rangle$; we define the Smarandache pseudo commutator (S-pseudo commutator) of $P \subset L$ (P a S-sub-N-loop of L) to be the S-sub N-loop generated by $\langle\{p \in P \mid a (xb) = p ($*

*[b x] a), a, b $\in$ P}〉; denoted by $P_S^N$ (L). If $P_S^N$ (L) is not a S-sub N loop then we say $P_S^N$ (L) = $\phi$.*

**DEFINITION 1.3.56:** *Let L = {$L_1 \cup L_2 \cup \ldots \cup L_N$, $*_1, \ldots, *_N$} be a N-loop. An associative triple a, b, c $\in P \subset L$ where P is a S-sub N- loop of L is said to be Smarandache pseudo associative (S-pseudo associative) if (a b) (x c) = (a x) (b c) for all x $\in$ P. If (a b) (x c) = (a x) (b c) for some x $\in$ P we say the triple is Smarandache pseudo associative (S-pseudo associative) relative to those x in P.*

*If in particular every associative triple in P is S-pseudo associative then we say the N-loop is a Smarandache pseudo associative N-loop (S- pseudo associative N-loop).*

*Thus for a N-loop to be a S-pseudo associative N-loop it is sufficient that one of its S-sub N-loops are a S-pseudo associative N-loop.*

*The Smarandache pseudo associator of a N-loop (S-pseudo associator of a N-loop) L denoted by $S\left(AL_p^N\right)$ = 〈{t $\in$ P / (a b) (t c) = (at) (bc) where a (bc) = (a b) c for a, b, c $\in$ P}〉 where P is a S-sub N-loop of L.*

**DEFINITION 1.3.57:** *Let L = {$L_1 \cup L_2 \cup \ldots \cup L_N$, $*_1, \ldots, *_N$} be a N-loop. L is said to be a Smarandache Moufang N loop (S-Moufang N-loop) if there exists S-sub N-loop P of L which is a Moufang N-loop.*

*A N-loop L is said to be a Smarandache Bruck N-loop (S-Bruck N-loop) if L has a proper S-sub N-loop P, where P is a Bruck N-loop.*

*Similarly a N-loop L is said to be a Smarandache Bol N-loop (S-Bol N-loop) if L has a proper S-sub N loop) K, where K is a Bol N-loop.*

*We call a N-loop L to be a Smarandache Jordan N-loop (S-Jordan N-loop) if L has a S-sub N-loop T, where T is a Jordan N-loop.*

On similar lines we define Smarandache right (left) alternative N-loop.

Now proceed on to define the Smarandache N-loop homomorphism.

**DEFINITION 1.3.58:** *Let $L = \{L_1 \cup L_2 \cup \ldots \cup L_N, *_1, \ldots, *_N\}$ and $K = \{ K_1 \cup K_2 \cup \ldots \cup K_N, o_1, o_2, \ldots, o_N\}$ be two N-loops. We say a map $\phi = \phi_1 \cup \phi_2 \cup \ldots \cup \phi_N$ from L to K is a Smarandache N-loop homomorphism (S-N loop homomorphism) if $\phi$ is a N-group homomorphism from P to T where P and T are proper subsets of L and K respectively such that P and T are N-groups of L and K.*

Thus for Smarandache $\phi$ - homomorphism of N-loops we need not have $\phi$ to be defined on the whole of L. It is sufficient if $\phi$ is defined on a proper subset P of L which is a N-group.

Now we proceed on to define Smarandache Moufang center of a N loop.

**DEFINITION 1.3.59:** *Let $L = \{L_1 \cup L_2 \cup \ldots \cup L_N, *_1, \ldots, *_N\}$be a N-loop. Let P be a S-sub N-loop of L. The Smarandache Moufang Centre (S-Moufang Centre) of L is defined to be the Moufang centre of the S-sub N-loop; $P = (P_1 \cup P_2 \cup \ldots \cup P_N, *_1, \ldots, *_N)$.*

Thus for a given N-loop L we can have several S-Moufang centres or no Moufang centre if the N-loop L has no S-Sub N-loop.

Now we proceed on to define the notion of S-center.

**DEFINITION 1.3.60:** *Let L be a N-loop. P be a S-sub N-loop of L. The Smarandache N-centre (S N center) $SZ_N$ (L) is the center of the sub-N-loop $P \subset L$.*

Thus even in case of S-centre for a N-loop we may have several S-centers depending on the number of S-sub N-loops.

Now notion of Smarandache middle, left and right nucleus of a N-loop is defined.

**DEFINITION 1.3.61:** *Let L be a N-loop, P be a S-sub N-loop of L. To define Smarandache nucleus (S-nucleus) related to P we*

*consider the Smarandache left nucleus (S-left nucleus) of L to be the left nucleus of the S-sub N-loop P denoted by* $N_\lambda^{P_N}$.

*Similarly the Smarandache right nucleus (S-right nucleus) of the N-loop L to be the right nucleus of the S-sub N-loop P denoted by* $N_P^{P_N}$.

*We can define the Smarandache middle nucleus (S-middle nucleus) of the N-loop L to be the middle nucleus of the S-sub N-loop P, denoted by* $N_m^{P_N}$*. Thus the Smarandache nucleus (S-nucleus) is* $[SN(L)]_N = N_\lambda^{P_N} \cap N_P^{P_N} \cap N_M^{P_N}$.

**DEFINITION 1.3.62:** *A neutrosophic loop is generated by a loop L and I denoted by* $\langle L \cup I\rangle$*. A neutrosophic loop in general need not be a loop for* $I^2 = I$ *and I may not have an inverse but every element in a loop has an inverse.*

**DEFINITION 1.3.63:** *Let* $\langle L \cup I\rangle$ *be a neutrosophic loop. A proper subset* $\langle P \cup I\rangle$ *of* $\langle L \cup I\rangle$ *is called the neutrosophic subloop, if* $\langle P \cup I\rangle$ *is itself a neutrosophic loop under the operations of* $\langle L \cup I\rangle$.

**DEFINITION 1.3.64:** *Let* $\langle L_n(m) \cup I\rangle = \{e, 1, 2, \ldots, n, e.I, 1I, \ldots, nI\}$*, where* $n > 3$*, n is odd and m is a positive integer such that* $(m, n) = 1$ *and* $(m – 1, n) = 1$ *with* $m < n$*. Define on* $\langle L_n(m) \cup I\rangle$ *a binary operation ' . ' as follows.*

i. $e.i. = i.e. = i$ *for all* $i \in L_n(m)$.
ii. $i^2 = e$ *for all* $i \in L_n(m)$.
iii. $iI.\ iI = e\,I$ *for all* $i \in L_n(m)$.
iv. $i.\ j = t$ *where* $t = (mj – (m – 1)i) \pmod n$ *for all* $i, j \in L_n(m)$, $i \neq j$, $i \neq e$ *and* $j \neq e$.
v. $iI.\ jI = tI$ *where* $t = (mj – (m – 1) i) \pmod n$ *for all* $iI, jI \in \langle L_n(m) \cup I\rangle$. $\langle L_n(m) \cup I\rangle$ *is a neutrosophic loop of order* $2(n + 1)$.

**DEFINITION 1.3.65:** *Let* $(\langle B \cup I\rangle, *_1, *_2)$ *be a non empty neutrosophic set with two binary operations* $*_1, *_2$, $\langle B \cup I\rangle$ *is a neutrosophic biloop if the following conditions are satisfied.*

i. *$\langle B \cup I\rangle = P_1 \cup P_2$ where $P_1$ and $P_2$ are proper subsets of $\langle B \cup I\rangle$.*
ii. *$(P_1, *_1)$ is a neutrosophic loop.*
iii. *$(P_2, *_2)$ is a group or a loop.*

**DEFINITION 1.3.66:** *Let $(\langle B \cup I\rangle, *_1, *_2)$ be a neutrosophic biloop. A proper subset P of $\langle B \cup I\rangle$ is said to be a neutrosophic subbiloop of $\langle B \cup I\rangle$ if $(P, *_1, *_2)$ is itself a neutrosophic biloop under the operations of $\langle B \cup I\rangle$.*

**DEFINITION 1.3.67:** *Let $(B = B_1 \cup B_2, *_1, *_2)$ be a neutrosophic biloop. A neutrosophic subbiloop $H = (H_1 \cup H_2, *_1, *_2)$ is said to be a neutrosophic normal subbiloop of B if*

i. *$x H_1 = H_1 x$*
ii. *$(H_1 x) y = H_1 (xy)$*
iii. *$y (x H_1) = (y x) H_1$ for all $x, y \in B_1$*
iv. *$H_2$ is a normal subgroup of $B_2$.*

*We call a neutrosophic biloop to be a simple neutrosophic biloop if it has no nontrivial neutrosophic normal subbiloops.*

**DEFINITION 1.3.68:** *Let $(B = B_1 \cup B_2, *_1, *_2)$ be a neutrosophic biloop if only one of the neutrosophic loop or the group is simple then we call the neutrosophic biloop to be a semi-simple neutrosophic biloop.*

**DEFINITION 1.3.69:** *Let $(B = B_1 \cup B_2, *_1, *_2)$ be non empty set with two binary operations. B is said to be strong neutrosophic biloop.*

i. *$B = B_1 \cup B_2$ is a union of proper subsets of B.*
ii. *$(B_1, *_1)$ is a neutrosophic loop.*
iii. *$(B_2, *_2)$ is a neutrosophic group.*

**DEFINITION 1.3.70:** *Let $B = (B_1 \cup B_2, *_1, *_2)$ be a proper set on which is defined two binary operations $*_1$ and $*_2$. We call B a* *neutrosophic biloop of type II if the following conditions are satisfied*

.

i. $B = B_1 \cup B_2$ *where* $B_1$ *and* $B_2$ *are proper subsets of* $B$.
ii. $(B_1, *_1)$ *is a neutrosophic loop.*
iii. $(B_2, *_2)$ *is a loop.*

*(Clearly we call a neutrosophic biloop of type I is one in which* $B_1$ *is a neutrosophic loop and* $B_2$ *is a group).*

**DEFINITION 1.3.71:** *Let* $S(B) = \{S(B_1) \cup \ldots \cup S(B_N), *_1, \ldots, *_N\}$ *be a non empty neutrosophic set with N binary operations. S(B) is a neutrosophic N-loop if* $S(B) = S(B_1) \cup \ldots \cup S(B_N)$, $S(B_i)$ *are proper subsets of S(B) for* $1 \leq i \leq N$*) and some of* $S(B_i)$ *are neutrosophic loops and some of the* $S(B_j)$ *are groups.*

*The order of the neutrosophic N-loop is the number of elements in S(B). If the number of elements in S(B) is in finite then we say S(B) is an infinite neutrosophic N-loop of infinite order. Thus even if one of the* $S(B_i)$ *is infinite we see the neutrosophic N-loop S(B) is infinite.*

**DEFINITION 1.3.72:** *Let* $S(B) = \{S(B_1) \cup S(B_2) \cup \ldots \cup S(B_N), *_1, \ldots, *_N\}$ *be a neutrosophic N-loop. A proper subset* $(P, *_1, \ldots, *_N)$ *of S(B) is said to be a neutrosophic sub N loop of S (B) if P itself is a neutrosophic N-loop under the operations of S (B).*

**DEFINITION 1.3.73:** *Let* $S(L) = \{L_1 \cup L_2 \cup \ldots \cup L_N, *_1, \ldots, *_N\}$ *where*

i. $L = L_1 \cup \ldots \cup L_N$ *is such that each* $L_i$ *is a proper subset of L,* $1 \leq i \leq N$.
ii. *Some of* $(L_i, *_i)$ *are a Neutrosophic loops.*
iii. *Some of* $(L_j, *_j)$ *are just loops.*
iv. *Some of* $(L_K, *_K)$ *are groups and rest of*
v. $(L_t, *_t)$ *are neutrosophic groups.*

*Then we call* $L = L_1 \cup L_2 \cup \ldots \cup L_N$ , $*_1, \ldots, *_N\}$ *to be a neutrosophic N-loop of level II.*

**DEFINITION 1.3.74:** *Let $L = \{L_1 \cup L_2 \cup \ldots \cup L_N, *_1, *_2, \ldots, *_N\}$ be a neutrosophic N-loop. A proper subset H of L is said to be a neutrosophic normal sub N-loop of L if the following conditions are satisfied.*

i. *H is a neutrosophic sub N-loop of L*
ii. *$xH = H x$*
*$(H x) y = H (xy)$*
*$y (x H) = (y x) H$*

*for all $x, y \in L$.*

*If the neutrosophic N-loop L has no trivial normal sub-N-loop, we call L to be a simple neutrosophic N-loop.*

**DEFINITION 1.3.75:** *Let $\{\langle L \cup I\rangle = L_1 \cup L_2 \cup \ldots \cup L_N, *_1, *_2, \ldots, *_N\}$, be a nonempty set with N-binary operations where*

i. *$\langle L \cup I\rangle = L_1 \cup L_2 \cup \ldots \cup L_N$ where each $L_i$ is a proper subset of $\langle L \cup I\rangle$; $1 \leq i \leq N$*
ii. *$(L_i, *_i)$ is a neutrosophic loop, at least for some i.*
iii. *$(L_j, *_j)$ is a neutrosophic group.*

*Then we call $\{\langle L \cup I\rangle, *_1, \ldots, *_N\}$ to be a strong neutrosophic N-loop.*

**DEFINITION 1.3.76:** *Let $\{\langle L \cup I\rangle = L_1 \cup L_2 \cup \ldots \cup L_N, *_1, \ldots, *_N\}$, be a neutrosophic N-loop. Let P be a proper subset of $\langle L \cup I\rangle$ such that*

$$P = \left\{P_{L_1} \cup \ldots \cup P_{L_t}, *_{i_1}, \ldots, *_{i_t} \mid 1 \leq i_1, \ldots, i_t \leq N \text{ and } 1 < t < N\right\}.$$

*If P is a neutrosophic t-loop ($*_{i_p} = *_j$, $1 \leq j \leq N$; $1 \leq i_p \leq t$ and $P_{i_p} = P \cap L_j$) then we call P a neutrosophic (N – t) deficit N-subloop of $\langle L \cup I\rangle$.*

**DEFINITION 1.3.77:** *Let $\{\langle L \cup I\rangle = L_1 \cup L_2 \cup \ldots \cup L_N, *_1, \ldots, *_N\}$ be a neutrosophic N-loop of finite order. A neutrosophic (N – t) deficit sub N-loop P of $\langle L \cup I\rangle$ is said to be Lagrange neutrosophic (N – t) deficit sub N-loop if $o(P) / o(\langle L \cup I\rangle)$. If* *every neutrosophic (N – t) deficit sub N-loop of ⟨L ∪ I⟩ is Lagrange neutrosophic then we call ⟨L ∪ I⟩ a Lagrange (N – t) deficit neutrosophic N-loop. If ⟨L ∪ I⟩ has atleast one Lagrange (N – t) deficit neutrosophic sub N-loop then we call ⟨L ∪ I⟩ a weak Lagrange (N – t) deficit neutrosophic N-loop.*

## 1.4 S-groupoids, neutrosophic groupoids and their generalizations

In this section we just recall the notion of groupoids. We also give some new classes of groupoids constructed using the set of modulo integers. This book uses in several examples from these new classes of groupoids. For more about groupoids please refer [45]. In this section we just give a brief introduction to groupoids, N-groupoids and their Smarandache analogue, neutrosophic groupoids and their generalizations.

**DEFINITION 1.4.1:** *Given an arbitrary set P a mapping of $P \times P$ into P is called a binary operation on P. Given such a mapping $\sigma: P \times P \to P$ we use it to define a product $*$ in P by declaring $a * b = c$ if $\sigma(a, b) = c$.*

**DEFINITION 1.4.2:** *A non empty set of elements G is said to form a groupoid if in G is defined a binary operation called the product denoted by $*$ such that $a * b \in G$ for all $a, b \in G$.*

**DEFINITION 1.4.3:** *A groupoid G is said to be a commutative groupoid if for every $a, b \in G$ we have $a * b = b * a$.*

**DEFINITION 1.4.4:** *A groupoid G is said to have an identity element e in G if $a * e = e * a = a$ for all $a \in G$.*

**DEFINITION 1.4.5:** *Let $(G, *)$ be a groupoid a proper subset $H \subset G$ is a subgroupoid if $(H, *)$ is itself a groupoid.*

**DEFINITION 1.4.6:** *A groupoid G is said to be a Moufang groupoid if it satisfies the Moufang identity $(xy)(zx) = (x(yz))x$ for all x, y, z in G.*

**DEFINITION 1.4.7:** *A groupoid G is said to be a Bol groupoid if G satisfies the Bol identity ((xy) z) y = x ((yz) y) for all x, y, z in G.*

**DEFINITION 1.4.8:** *A groupoid G is said to be a P-groupoid if (xy) x = x (yx) for all x, y $\in$ G.*

**DEFINITION 1.4.9:** *A groupoid G is said to be right alternative if it satisfies the identity (xy) y = x (yy) for all x, y $\in$ G. Similarly we define G to be left alternative if (xx) y = x (xy) for all x, y $\in$ G.*

**DEFINITION 1.4.10:** *A groupoid G is alternative if it is both right and left alternative, simultaneously.*

**DEFINITION 1.4.11:** *Let (G, $*$) be a groupoid. A proper subset H of G is said to be a subgroupoid of G if (H, $*$) is itself a groupoid.*

**DEFINITION 1.4.12:** *A groupoid G is said to be an idempotent groupoid if $x^2 = x$ for all x $\in$ G.*

**DEFINITION 1.4.13:** *Let G be a groupoid. P a non empty proper subset of G, P is said to be a left ideal of the groupoid G if 1) P is a subgroupoid of G and 2) For all x $\in$ G and a $\in$ P, xa $\in$ P. One can similarly define right ideal of the groupoid G. P is called an ideal if P is simultaneously a left and a right ideal of the groupoid G.*

**DEFINITION 1.4.14:** *Let G be a groupoid A subgroupoid V of G is said to be a normal subgroupoid of G if*

i. *aV = Va*
ii. *(Vx)y = V(xy)*
iii. *y(xV) = (yx)V*

*for all x, y, a $\in$ V.*

**DEFINITION 1.4.15:** *A groupoid G is said to be simple if it has no non trivial normal subgroupoids.*

**DEFINITION 1.4.16:** *A groupoid G is normal if*

i. $xG = Gx$
ii. $G(xy) = (Gx)y$
iii. $y(xG) = (yx)G$ *for all* $x, y \in G$.

**DEFINITION 1.4.17:** *Let G be a groupoid H and K be two proper subgroupoids of G, with* $H \cap K = \phi$. *We say H is conjugate with K if there exists a* $x \in H$ *such that* $H = x K$ *or* $Kx$ *('or' in the mutually exclusive sense).*

**DEFINITION 1.4.18:** *Let* $(G_1, \theta_1), (G_2, \theta_2), \ldots, (G_n, \theta_n)$ *be n groupoids with* $\theta_i$ *binary operations defined on each* $G_i$, $i = 1, 2, 3, \ldots, n$. *The direct product of* $G_1, \ldots, G_n$ *denoted by* $G = G_1 \times \ldots \times G_n = \{(g_1, \ldots, g_n) \mid g_i \in G_i\}$ *by component wise multiplication on G, G becomes a groupoid.*

*For if* $g = (g_1, \ldots, g_n)$ *and* $h = (h_1, \ldots, h_n)$ *then* $g \bullet h = \{(g_1\theta_1h_1, g_2\theta_2h_2, \ldots, g_n\theta_nh_n)\}$. *Clearly,* $gh \in G$. *Hence G is a groupoid.*

**DEFINITION 1.4.19:** *Let G be a groupoid we say an element* $e \in G$ *is a left identity if* $ea = a$ *for all* $a \in G$. *Similarly we can define right identity of the groupoid G, if* $e \in G$ *happens to be simultaneously both right and left identity we say the groupoid G has an identity.*

**DEFINITION 1.4.20:** *Let G be a groupoid. We say a in G has right zero divisor if* $a * b = 0$ *for some* $b \neq 0$ *in G and a in G has left zero divisor if* $b * a = 0$. *We say G has zero divisors if* $a \bullet b = 0$ *and* $b * a = 0$ *for* $a, b \in G \setminus \{0\}$ *A zero divisor in G can be left or right divisor.*

**DEFINITION 1.4.21:** *Let G be a groupoid. The center of the groupoid* $C(G) = \{x \in G \mid ax = xa \text{ for all } a \in G\}$.

**DEFINITION 1.4.22:** *Let G be a groupoid. We say a, b $\in$ G is a conjugate pair if a = bx (or xa for some x $\in$ G) and b = ay (or ya for some y $\in$ G).*

**DEFINITION 1.4.23:** *Let G be a groupoid of order n. An element a in G is said to be right conjugate with b in G if we can find x, y $\in$ G such that a • x = b and b • y = a (x $*$ a = b and y $*$ b = a). Similarly, we define left conjugate.*

**DEFINITION 1.4.24:** *Let $Z^+$ be the set of integers. Define an operation $*$ on $Z^+$ by x $*$ y = mx + ny where m, n $\in Z^+$, m < $\propto$ and n < $\propto$ (m, n) = 1 and m $\neq$ n. Clearly $\{Z^+, *, (m, n)\}$ is a groupoid denoted by $Z^+$ (m, n). We have for varying m and n get infinite number of groupoids of infinite order denoted by $\mathbf{Z}^+$.*

Here we define a new class of groupoids denoted by Z(n) using $Z_n$ and study their various properties.

**DEFINITION 1.4.25:** *Let $Z_n$ = {0, 1, 2, ... , n – 1} n $\geq$ 3. For a, b $\in Z_n \setminus \{0\}$ define a binary operation $*$ on $Z_n$ as follows. a $*$ b = ta + ub (mod n) where t, u are 2 distinct elements in $Z_n \setminus \{0\}$ and (t, u) =1 here ' + ' is the usual addition of 2 integers and ' ta ' means the product of the two integers t and a. We denote this groupoid by $\{Z_n, (t, u), *\}$ or in short by $Z_n$ (t, u).*

It is interesting to note that for varying t, u $\in Z_n \setminus \{0\}$ with (t, u) = 1 we get a collection of groupoids for a fixed integer n. This collection of groupoids is denoted by Z(n) that is Z(n) = $\{Z_n$, (t, u), $*$ | for integers t, u $\in Z_n \setminus \{0\}$ such that (t, u) = 1}. Clearly every groupoid in this class is of order n.

***Example 1.4.1:*** Using $Z_3$ = {0, 1, 2}. The groupoid $\{Z_3$, (1, 2), $*\}$ = ($Z_3$ (1, 2)) is given by the following table:

| $*$ | 0 | 1 | 2 |
|---|---|---|---|
| 0 | 0 | 2 | 1 |
| 1 | 1 | 0 | 2 |
| 2 | 2 | 1 | 0 |

Clearly this groupoid is non associative and non commutative and its order is 3.

**THEOREM 1.4.1:** *Let $Z_n = \{0, 1, 2, \ldots, n\}$. A groupoid in Z (n) is a semigroup if and only if $t^2 \equiv t \pmod n$ and $u^2 \equiv u \pmod n$ for t, $u \in Z_n \setminus \{0\}$ and (t, u) = 1.*

**THEOREM 1.4.2:** *The groupoid $Z_n$ (t, u) is an idempotent groupoid if and only if $t + u \equiv 1 \pmod n$.*

**THEOREM 1.4.3:** *No groupoid in Z (n) has {0} as an ideal.*

**THEOREM 1.4.4:** *P is a left ideal of $Z_n$ (t, u) if and only if P is a right ideal of $Z_n$ (u, t).*

**THEOREM 1.4.5:** *Let $Z_n$ (t, u) be a groupoid. If n = t + u where both t and u are primes then $Z_n$ (t, u) is simple.*

**DEFINITION 1.4.26:** *Let $Z_n = \{0, 1, 2, \ldots, n-1\}$ $n \geq 3$, $n < \infty$. Define $*$ a closed binary operation on $Z_n$ as follows. For any a, $b \in Z_n$ define $a * b = at + bu \pmod n$ where (t, u) need not always be relatively prime but $t \neq u$ and t, $u \in Z_n \setminus \{0\}$.*

**THEOREM 1.4.6:** *The number of groupoids in $Z^*(n)$ is (n – 1) (n – 2).*

**THEOREM 1.4.7:** *The number of groupoids in the class Z (n) is bounded by (n – 1) (n – 2).*

**THEOREM 1.4.8:** *Let $Z_n$ (t, u) be a groupoid in $Z^*(n)$ such that (t, u) = t, n = 2m, t / 2m and t + u = 2m. Then $Z_n$ (t, u) has subgroupoids of order 2m / t or n / t.*

*Proof:* Given n is even and t + u = n so that u = n – t.

Thus $Z_n$ (t, u) = $Z_n$ (t, n – t). Now using the fact $t \cdot Z_n = \left\{0, t, 2t, 3t, \ldots, \left(\frac{n}{t} - 1\right)t\right\}$ that is t . $Z_n$ has only n / t elements and these n / t elements from a subgroupoid. Hence $Z_n$ (t, n – t) where (t, n – t) = t has only subgroupoids of order n / t.

**DEFINITION 1.4.27:** *Let $Z_n = \{0, 1, 2, \ldots, n - 1\}$ $n \geq 3$, $n < \infty$. Define $*$ on $Z_n$ as $a * b = ta + ub \pmod{n}$ where t and $u \in Z_n \setminus \{0\}$ and t can also equal u. For a fixed n and for varying t and u we get a class of groupoids of order n which we denote by $Z^{**}(n)$.*

**DEFINITION 1.4.28:** *Let $Z_n = \{0, 1, 2, \ldots, n - 1\}$ $n \geq 3$, $n < \infty$. Define $*$ on $Z_n$ as follows. $a * b = ta + ub \pmod{n}$ where $t, u \in Z_n$. Here t or u can also be zero.*

**DEFINITION 1.4.29:** *Let $G = (G_1 \cup G_2 \cup \ldots \cup G_N; *_1, \ldots, *_N)$ be a non empty set with N-binary operations. G is called the N-groupoid if some of the $G_i$'s are groupoids and some of the $G_j$'s are semigroups, $i \neq j$ and $G = G_1 \cup G_2 \cup \ldots \cup G_N$ is the union of proper subsets of G.*

**DEFINITION 1.4.30:** *Let $G = (G_1 \cup G_2 \cup \ldots \cup G_N; *_1, *_2, \ldots, *_N)$ be a N-groupoid. The order of the N-groupoid G is the number of distinct elements in G. If the number of elements in G is finite we call G a finite N-groupoid. If the number of distinct elements in G is infinite we call G an infinite N-groupoid.*

**DEFINITION 1.4.31:** *Let $G = \{G_1 \cup G_2 \cup \ldots \cup G_N, *_1, *_2, \ldots, *_N\}$ be a N-groupoid. We say a proper subset $H = \{H_1 \cup H_2 \cup \ldots \cup H_N, *_1, *_2, \ldots, *_N\}$ of G is said to be a sub N-groupoid of G if H itself is a N-groupoid of G.*

**DEFINITION 1.4.32:** *Let $G = (G_1 \cup G_2 \cup \ldots \cup G_N, *_1, *_2, \ldots, *_N)$ be a finite N-groupoid. If the order of every sub N-groupoid H divides the order of the N-groupoid, then we say G is a Lagrange N-groupoid.*

It is very important to mention here that in general every N-groupoid need not be a Lagrange N-groupoid.

Now we define still a weaker notion of a Lagrange N-groupoid.

**DEFINITION 1.4.33:** *Let $G = (G_1 \cup G_2 \cup \ldots \cup G_N, *_1, \ldots, *_N)$ be a finite N-groupoid. If G has atleast one nontrivial sub N-groupoid H such that $o(H) / o(G)$ then we call G a weakly Lagrange N-groupoid.*

**DEFINITION 1.4.34:** *Let $G = (G_1 \cup G_2 \cup \ldots \cup G_N, *_1, *_2, \ldots, *_N)$ be a N groupoid. A sub N groupoid $V = V_1 \cup V_2 \cup \ldots \cup V_N$ of G is said to be a normal sub N-groupoid of G; if*

i. $a V_i = V_i a$, *$i = 1, 2, \ldots, N$ whenever $a \in V_i$*
ii. $V_i (x y) = (V_i x) y$, *$i = 1, 2, \ldots, N$ for $x, y \in V_i$*
iii. $y (x V_i) = (xy) V_i$, *$i = 1, 2, \ldots, N$ for $x, y \in V_i$.*

*Now we say a N-groupoid G is simple if G has no nontrivial normal sub N-groupoids.*

Now we proceed on to define the notion of N-conjugate groupoids.

**DEFINITION 1.4.35:** *Let $G = (G_1 \cup G_2 \cup \ldots \cup G_N, *_1, \ldots, *_N)$ be a N-groupoid. Let $H = \{H_1 \cup \ldots \cup H_N; *_1, \ldots, *_N\}$ and $K = \{K_1 \cup K_2 \cup \ldots \cup K_N, *_1, \ldots, *_N\}$ be sub N-groupoids of $G = G_1 \cup G_2 \cup \ldots \cup G_N$; where $H_i$, $K_i$ are subgroupoids of $G_i$ ($i = 1, 2, \ldots, N$).*

*Let $K \cap H = \phi$. We say H is N-conjugate with K if there exists $x_i \in H_i$ such that $x_i K_i = H_i$ (or $K_i x_i = H_i$) for $i = 1, 2, \ldots, N$ 'or' in the mutually exclusive sense.*

**DEFINITION 1.4.36:** *Let $G = (G_1 \cup G_2 \cup \ldots \cup G_N, *_1, *_2, \ldots, *_N)$ be a N-groupoid. An element x in G is said to be a zero divisor if their exists a y in G such that $x *_i y = y *_i x = 0$ for some i in $\{1, 2, \ldots, N\}$.*

We define N-centre of a N-groupoid G.

**DEFINITION 1.4.37:** *Let $G = (G_1 \cup G_2 \cup \ldots \cup G_N, *_1, *_2, \ldots, *_N)$ be a N-groupoid. The N-centre of $(G_1 \cup G_2 \cup \ldots \cup G_N,$*

*$*_1, \ldots, *_N$) denoted by NC (G) = $\{x_1 \in G_1 \mid x_1 a_1 = a_1 x_1$ for all $a_1 \in G_1\} \cup \{x_2 \in G_2 \mid x_2 a_2 = a_2 x_2$ for all $a_2 \in G_2\} \cup \ldots \cup \{x_N \in G_N \mid x_N a_N = a_N x_N$ for all $a_N \in G_N\}$ = $\bigcup_{i=1}^{N} \{x_i \in G \mid x_i a_i = a_i x_i \text{ for all } a_i \in G_i\}$ = NC (G).*

**DEFINITION 1.4.38:** *A Smarandache groupoid G is a groupoid which has a proper subset S, S $\subset$ G such that S under the operations of G is a semigroup.*

**DEFINITION 1.4.39:** *Let G be a Smarandache groupoid (SG) if the number of elements in G is finite we say G is a finite SG, otherwise G is said to be an infinite SG.*

**DEFINITION 1.4.40:** *Let G be a Smarandache groupoid. G is said to be a Smarandache commutative groupoid if there is a proper subset, which is a semigroup, is a commutative semigroup.*

**DEFINITION 1.4.41:** *Let (G, $*$) be a Smarandache groupoid (SG). A non-empty subset H of G is said to be a Smarandache subgroupoid if H contains a proper subset K $\subset$ H such that K is a semigroup under the operation $*$.*

**DEFINITION 1.4.42:** *A Smarandache left ideal A of the Smarandache groupoid G satisfies the following conditions:*

i. *A is a Smarandache subgroupoid*
ii. *$x \in G$ and $a \in A$ then $xa \in A$.*

*Similarly, we can define Smarandache right ideal. If A is both a Smarandache right ideal and Smarandache left ideal simultaneously then we say A is a Smarandache ideal.*

**DEFINITION 1.4.43:** *Let G be a SG. V be a Smarandache subgroupoid of G. We say V is a Smarandache seminormal groupoid if*

i. *aV = X for all a* $\in$ *G*
ii. *Va = Y for all a* $\in$ *G*

*where either X or Y is a Smarandache subgroupoid of G but X and Y are both subgroupoids.*

**DEFINITION 1.4.44:** *Let G be a SG and V be a Smarandache subgroupoid of G. V is said to be a Smarandache normal groupoid if aV = X and Va = Y for all a* $\in$ *G where both X and Y are Smarandache subgroupoids of G.*

**DEFINITION 1.4.45***:* *Let G be a SG. H and P be any two subgroupoids of G. We say H and P are Smarandache semiconjugate subgroupoids of G if*

i. *H and P are Smarandache subgroupoids of G*
ii. *H = xP or Px or*
iii. *P = xH or Hx for some x* $\in$ *G.*

**DEFINITION 1.4.46:** *Let G be a Smarandache groupoid. H and P be subgroupoids of G. We say H and P are Smarandache conjugate subgroupoids of G if*

i. *H and P are Smarandache subgroupoids of G*
ii. *H = xP or Px and*
iii. *P = xH or Hx.*

**DEFINITION 1.4.47:** *Let* $G = \{G_1 \cup G_2 \cup \ldots \cup G_N, *_1, \ldots, *_N\}$ *be a non empty set with N-binary operations and* $G_i$ *are proper subsets of G. We call G a Smarandache N-groupoid (S-N-groupoid) if some of the* $G_j$ *are S-groupoids and the rest of them are S-semigroups; i.e. each* $G_i$ *is either a S-semigroup or a S-groupoid.*

**DEFINITION 1.4.48:** *Let* $G = \{G_1 \cup G_2 \cup \ldots \cup G_N, *_1, \ldots, *_N\}$ *be a N-groupoid. A non empty proper subset* $K = (K_1 \cup K_2 \cup \ldots \cup K_N, *_1, \ldots, *_N)$ *is said to be a Smarandache sub-N-groupoid (S-sub N-groupoid) if K itself is a S-N-groupoid.*

**DEFINITION 1.4.49:** *Let $G = \{G_1 \cup G_2 \cup \ldots \cup G_N, *_1, \ldots, *_N\}$ be a N-groupoid. We say G is a Smarandache commutative N-groupoid (S-commutative N-groupoid) if every S-sub N-groupoid of G is commutative. If at least one S-sub-N-groupoid which is commutative then we call G a Smarandache weakly commutative N-groupoid (S-weakly commutative N-groupoid).*

**DEFINITION 1.4.50:** *Let $G = \{G_1 \cup G_2 \cup \ldots \cup G_N, *_1, *_2, \ldots, *_N\}$ be a N-groupoid. We call G an idempotent N-groupoid if every element of G is an idempotent. We call the S-N-groupoid to be a Smarandache idempotent N-groupoid (S-idempotent N-groupoid) if every element in every S- sub N-groupoid in G is a S-idempotent.*

**DEFINITION 1.4.51:** *Let $G = \{G_1 \cup G_2 \cup \ldots \cup G_N, *_1, *_2 \ldots, *_N\}$ be a N-groupoid. A non empty proper subset $P = (P_1 \cup P_2 \cup \ldots \cup P_N, *_1, *_2, \ldots, *_N)$ of G is said to be a Smarandache left N-ideal (S-left N ideal) of the N-groupoid G if*

i. *P is a S-sub N-groupoid.*
ii. *Each $(P_i, *_i)$ is a left ideal of $G_i$, $1 \leq i \leq N$.*

*On similar lines we can define S-right N-ideal. If P is both a S-right ideal and S-left ideal then we say P is a Smarandache N-ideal (S-N-ideal) of G.*

**DEFINITION 1.4.52:** *Let $G = \{G_1 \cup G_2 \cup \ldots \cup G_N , *_1, *_2, \ldots, *_N\}$ be a N-groupoid. A S-sub N-groupoid $K = (K_1 \cup K_2 \cup \ldots \cup K_N, *_1, \ldots, *_N)$ is said to be a Smarandache normal sub N-groupoid (S-normal sub N-groupoid) of G if*

i. $a_i P_i = P_i a_i$, $a_i \in P_i$ .
ii. $P_i (x_i y_i) = (P_i x_i) y_i$ . $x_i , y_i \in P_i$ .
iii. $x_i (y_i P_i) = (x_i y_i)P_i$ .

*for all $a_i$, $x_i y_i \in G$.*

*We call the Smarandache simple N-groupoid (S-simple N-groupoid) if it has no non trivial S- normal sub- N-groupoids.*

**DEFINITION 1. 4.53:** *Let* $G = \{G_1 \cup G_2 \cup \ldots \cup G_N, *_1, *_2, \ldots, *_N\}$ *be a finite S- N-groupoid.*

*We say G to be a Smarandache normal N-groupoid (S-normal N-groupoid) if the largest S-sub N-groupoid of G is a normal N-groupoid. We call a S-sub N-groupoid to be the largest if the number of elements in it is the maximum.*

For more please refer [35-40]

**DEFINITION 1.4.54:** *Let* $G = \{G_1 \cup G_2 \cup \ldots \cup G_N, *_1, *_2, \ldots, *_N\}$ *be a N-groupoid. Let* $H = \{H_1 \cup H_2 \cup \ldots \cup H_N, *_1, *_2, \ldots, *_N\}$ *and* $K = \{K_1 \cup K_2 \cup \ldots \cup K_N, *_1, *_2, \ldots, *_N\}$ *be any two S-sub N-groupoids of G. We say H is Smarandache N-conjugate (S-N-conjugate) with K if there exists* $x_i \in H_i$ *and* $x_j \in H_j$ *with* $x_i K_i = H_i$ *(or* $K_i x_i$*) and* $x_j K_j = H_j$ *(or* $K_j x_j$*) 'or' in the mutually exclusive sense.*

**DEFINITION 1.4.55:** *Let (G, *) be a groupoid. A neutrosophic groupoid is defined as a groupoid generated by* $\{\langle G \cup I\rangle\}$ *under the operation *.*

**DEFINITION 1.4.56:** *Let* $\{\langle G \cup I\rangle, *\}$ *be a neutrosophic groupoid. The number of distinct elements in* $\{\langle G \cup I\rangle, *\}$ *is called the order of the neutrosophic groupoid. If* $\{\langle G \cup I\rangle, *\}$ *has infinite number of elements then we say the neutrosophic groupoid is infinite. If* $\{\langle G \cup I\rangle, *\}$ *has only a finite number of elements then we say* $\{\langle G \cup I\rangle, *\}$ *is a finite neutrosophic groupoid.*

**DEFINITION 1.4.57:** *Let* $\{\langle G \cup I\rangle, *\}$ *be a neutrosophic groupoid. A proper subset H of* $\langle G \cup I\rangle$ *is said to be a neutrosophic subgroupoid of* $\langle G \cup I\rangle$ *if (H, *) itself is a neutrosophic groupoid.*

*A non empty subset P of the neutrosophic groupoid* $\langle G \cup I\rangle$ *is said to be a left neutrosophic ideal of the neutrosophic groupoid* $\langle G \cup I\rangle$ *if*

*i. P is a neutrosophic subgroupoid.*

ii. *For all* $x \in \langle G \cup I \rangle$ *and* $a \in P$, $x * a \in P$.

*One can similarly define right neutrosophic ideal of a neutrosophic groupoid* $\langle G \cup I \rangle$. *We say P is a neutrosophic ideal of the neutrosophic groupoid* $\langle G \cup I \rangle$ *if P is simultaneously a left and a right neutrosophic ideal of* $\langle G \cup I \rangle$.

**DEFINITION 1.4.58:** *Let* $\langle G \cup I \rangle$ *be a neutrosophic groupoid. A neutrosophic subgroupoid V of* $\langle G \cup I \rangle$ *is said to be a neutrosophic normal subgroupoid of* $\langle G \cup I \rangle$ *if*

i. $aV = Va$
ii. $(Vx)\, y = V\,(xy)$
iii. $y\,(xV) = (yx)\, V$

*for all* $x, y, a \in \langle G \cup I \rangle$.

*A neutrosophic groupoid is said to be neutrosophic simple if it has no nontrivial neutrosophic normal subgroupoids.*

**DEFINITION 1.4.59:** *Let* $\langle G \cup I \rangle$ *be a neutrosophic groupoid. Let H and K be two proper neutrosophic subgroupoids of G with* $H \cap K = \phi$ *we say H is neutrosophic conjugate with K if there exists a* $x \in H$ *such that* $H = xK$ *(or Kx) (or in the mutually exclusive sense).*

**DEFINITION 1.4.60:** *Let* $\{\langle G \cup I \rangle, *_1\}, \{\langle G \cup I \rangle, *_2\}, \ldots, \{\langle G \cup I \rangle, *_n\}$ *be n neutrosophic groupoids* $*_i$ *binary operations defined on each* $\langle G_i \cup I \rangle$, $i = 1, 2, \ldots, n$. *The direct product of* $\langle G_1 \cup I \rangle, \langle G_2 \cup I \rangle, \ldots, \langle G_n \cup I \rangle$ *denoted by* $\langle G \cup I \rangle = \langle G_1 \cup I \rangle \times \langle G_2 \cup I \rangle \times \ldots \times \langle G_n \cup I \rangle = \{(g_1, g_2, \ldots, g_n) \mid g_i \in \langle G_i \cup I \rangle\ 1 \le i \le n\}$, *component wise multiplication of* $G_i$ *makes* $\langle G \cup I \rangle$ *a neutrosophic groupoid.*

*For if* $g = (g_i, \ldots, g_n)$ *and* $h = (h_1, h_2, \ldots, h_n)$ *in* $\langle G \cup I \rangle$ *then* $g * h = (g_1 *_1 h_1, g_2 *_2 h_2, \ldots, g_n *_n h_n)$. *Clearly* $g * h \in \{\langle G \cup I \rangle, *\}$. *Thus* $\{\langle G \cup I \rangle, *\}$ *is a neutrosophic groupoid.*

**DEFINITION 1.4.61:** *Let* $\{\langle G \cup I \rangle, *\}$ *be a neutrosophic groupoid, we say an element* $e \in \langle G \cup I \rangle$ *is a left identity if* $e * a$

*= a for all a ∈ ⟨G ∪ I⟩, similarly right identity of a neutrosophic groupoid can be defined. If e happens to be simultaneously both right and left identity we say the neutrosophic groupoid has an identity.*

*Similarly we can say an element $0 \neq a \in \langle G \cup I\rangle$ has a right zero divisor if $a * b = 0$ for some $b \neq 0$ in $\langle G \cup I\rangle$ and $a_1$ in $\langle G \cup I\rangle$ has left zero divisor if $b_1 * a_1 = 0$ (both $a_1$ and $b_1$ are different from zero). We say $\langle G \cup I\rangle$ has zero divisor if $a * b = 0$ and $b * a = 0$ for $a, b \in \langle G \cup I\rangle \setminus \{0\}$.*

**DEFINITION 1.4.62:** *Let $\{\langle G \cup I\rangle, *\}$ be a neutrosophic groupoid, the neutrosophic centre of the groupoid $\langle G \cup I\rangle$ is $C(\langle G \cup I\rangle) = \{a \in \langle G \cup I\rangle\} \setminus a * x = x * a$ for all $x \in \langle G \cup I\rangle\}$.*

**DEFINITION 1.4.63:** *Let $\{\langle G \cup I\rangle, *\}$ be a neutrosophic groupoid of order n. ($n < \infty$). We say $a, b \in \langle G \cup I\rangle$ is a conjugate pair if $a = b * x$ (or $x * b$ for some $x \in \langle G \cup I\rangle$) and $b = a * y$ (or $y * a$ for some $y \in \langle G \cup I\rangle$). An element a in $\langle G \cup I\rangle$ is said to be right conjugate with b in $\langle G \cup I\rangle$ if we can find $x, y \in \langle G \cup I\rangle$ such that $a * x = b$ and $b * y = a$ ($x * a = b$ and $y * b = a$).*

**DEFINITION 1.4.64:** *Let (BN(G), *, o) be a non empty set with two binary operations * and o. (BN(G), *, o) is said to be a neutrosophic bigroupoid if*

*$BN(G) = G_1 \cup G_2$ where at least one of $(G_1, *)$ or $(G_2, o)$ is a neutrosophic groupoid and other is just a groupoid. $G_1$ and $G_2$ are proper subsets of BN(G); i.e., $G_1 \not\subseteq G_2$.*

**DEFINITION 1.4.65:** *Let $N(G) = \{G_1 \cup G_2 \cup \ldots \cup G_N, *_1, \ldots, *_N\}$ be a non-empty set with N-binary operations, N(G) is called a neutrosophic N-groupoid if some of the $G_i$'s are neutrosophic groupoids and some of them are neutrosophic semigroups and $N(G) = G_1 \cup G_2 \cup \ldots \cup G_N$ is the union of the proper subsets of N(G).*

*It is important to note that a $G_i$ is either a neutrosophic groupoid or a neutrosophic semigroup. We call a neutrosophic N-groupoid to be a weak neutrosophic N-groupoid if in the* *union $N(G) = G_1 \cup G_2 \cup ... \cup G_N$ some of the $G_i$'s are neutrosophic groupoids, some of the $G_j$'s are neutrosophic semigroups and some of the $G_k$'s are groupoids or semigroups 'or' not used in the mutually exclusive sense.*

*The order of the neutrosophic N-groupoids are defined as that of N-groupoids. Further we call a neutrosophic N-groupoid to be commutative if each $(G_i, *_i)$ is commutative for $i = 1, 2, ..., N$.*

## 1.5 Mixed N-algebraic Structures and neutrosophic mixed N-algebraic structure

In this section we just recall some of the definition of mixed algebraic structure, S- mixed algebraic structure, and neutrosophic S- algebraic structure

In this section we proceed onto define the notion of N-group-semigroup algebraic structure and other mixed substructures and enumerate some of its properties.

**DEFINITION 1.5.1:** *Let $G = \{G_1 \cup G_2 \cup ... \cup G_N, *_1, ..., *_N\}$ where some of the $G_i$'s are groups and the rest of them are semigroups. $G_i$'s are such that $G_i \not\subseteq G_j$ or $G_j \not\subseteq G_I$ if $i \neq j$, i.e. $G_i$'s are proper subsets of G. $*_1, ..., *_N$ are N binary operations which are obviously are associative then we call G a N-group semigroup.*

We can also redefine this concept as follows:

**DEFINITION 1.5.2:** *Let G be a non empty set with N-binary operations $*_1, ..., *_N$. We call G a N-group semigroup if G satisfies the following conditions:*

- *i.* *$G = G_1 \cup G_2 \cup ... \cup G_N$ such that each $G_i$ is a proper subset of G (By proper subset $G_i$ of G we mean $G_i \not\subseteq G_j$ or $G_j \not\subseteq G_i$ if $i \neq j$. This does not necessarily imply $G_i \cap G_j = \phi$).*
- *ii.* *$(G_i, *_i)$ is either a group or a semigroup, $i = 1, 2,..., N$.*

iii. *At least one of the $(G_i, *_i)$ is a group.*
iv. *At least one of the $(G_j, *_j)$ is semigroup $i \neq j$.*

*Then we call $G = \{G_1 \cup G_2 \cup \ldots \cup G_N, *_1, \ldots, *_N\}$ to be a N-group semigroup ( $1 \leq i, j \leq N$).*

**DEFINITION 1.5.3:** *Let $G = \{G_1 \cup G_2 \cup \ldots \cup G_N, *_1, \ldots, *_N\}$ be a N-group-semigroup. We say G is a commutative N-group semigroup if each $(G_i, *_i)$ is a commutative structure, i = 1, 2, …, N.*

**DEFINITION 1.5.4:** *Let $G = \{G_1 \cup G_2 \cup \ldots \cup G_N, *_1, \ldots, *_N\}$ be a N-group. A proper subset P of G where $(P_1 \cup P_2 \cup \ldots \cup P_N, *_1, \ldots, *_N)$ is said to be a N-subgroup of the N-group semigroup G if each $(P_i, *_i)$ is a subgroup of $(G_i, *_i)$; i = 1, 2, …, N.*

**DEFINITION 1.5.5:** *Let $G = \{G_1 \cup G_2 \cup \ldots \cup G_N, *_1, \ldots, *_N\}$ be a N-group semigroup where some of the $(G_i, *_i)$ are groups and rest are $(G_j, *_j)$ are semigroups, $1 \leq i, j \leq N$. A proper subset P of G is said to be a N-subsemigroup if $P = \{P_1 \cup P_2 \cup \ldots \cup P_N, *_1, \ldots, *_N\}$ where each $(P_i, *_i)$ is only a semigroup under the operation $*_i$.*

Now we proceed on to define the notion of N-subgroup semigroup of a N-group semigroup.

**DEFINITION 1.5.6:** *Let $G = \{G_1 \cup G_2 \cup \ldots \cup G_N, *_1, \ldots, *_N\}$ be a N-group semigroup. Let P be a proper subset of G. We say P is a N-subgroup semigroup of G if $P = \{P_1 \cup P_2 \cup \ldots \cup P_N, *_1, \ldots, *_N\}$ and each $(P_i, *_i)$ is a group or a semigroup.*

**DEFINITION 1.5.7**: *Let $G = \{G_1 \cup G_2 \cup \ldots \cup G_N, *_1, \ldots, *_N\}$ be a N-group semigroup. We call a proper subset P of G where $P = \{P_1 \cup P_2 \cup \ldots \cup P_N, *_1, \ldots, *_N\}$ to be a normal N-subgroup semigroup of G if $(G_i, *_i)$ is a group then $(P_i, *_i)$ is a normal subgroup of $G_i$ and if $(G_j, *_j)$ is a semigroup then $(P_j, *_j)$ is an ideal of the semigroup $G_j$. If G has no normal N-subgroup semigroup then we say N-group semigroup is simple.*

**DEFINITION 1.5.8:** *Let $L = \{L_1 \cup L_2 \cup \ldots \cup L_N, *_1, \ldots, *_N\}$ be a non empty set with N-binary operations defined on it. We call L a N-loop groupoid if the following conditions are satisfied:*

i. *$L = L_1 \cup L_2 \cup \ldots \cup L_N$ where each $L_i$ is a proper subset of L i.e. $L_i \not\subseteq L_j$ or $L_j \not\subseteq L_i$ if $i \neq j$, for $1 \leq i, j \leq N$.*
ii. *$(L_i, *_i)$ is either a loop or a groupoid.*
iii. *There are some loops and some groupoids in the collection $\{L_1, \ldots, L_N\}$.*

*Clearly L is a non associative mixed N-structure.*

**DEFINITION 1.5.9:** *Let $L = \{L_1 \cup L_2 \cup \ldots \cup L_N, *_1, \ldots, *_N\}$ be a N-loop groupoid. L is said to be a commutative N-loop groupoid if each of $\{L_i, *_i\}$ is commutative.*

Now we give an example of a commutative N-loop groupoid.

**DEFINITION 1.5.10:** *Let $L = \{L_1 \cup L_2 \cup \ldots \cup L_N, *_1, \ldots, *_N\}$ be a N-loop groupoid. A proper subset P of L is said to be a sub N-loop groupoid if $P = \{P_1 \cup P_2 \cup \ldots \cup P_N, *_1, \ldots, *_N\}$ be a N-loop groupoid. A proper subset $P = \{P_1 \cup P_2 \cup \ldots \cup P_N, *_1, \ldots, *_N\}$ is such that if P itself is a N-loop groupoid then we call P the sub N-loop groupoid of L.*

**DEFINITION 1.5.11:** *Let $L = \{L_1 \cup L_2 \cup \ldots \cup L_N, *_1, \ldots, *_N\}$ be a N-loop groupoid. A proper subset $G = \{G_1 \cup G_2 \cup \ldots \cup G_N, *_1, \ldots, *_N\}$ is called a sub N-group if each $(G_i, *_i)$ is a group.*

**DEFINITION 1.5.12:** *Let $L = \{L_1 \cup L_2 \cup \ldots \cup L_N, *_1, \ldots, *_N\}$ be a N-loop groupoid. A proper subset $T = \{T_1 \cup T_2 \cup \ldots \cup T_N, *_1, \ldots, *_N\}$ is said to be a sub N-groupoid of the N-loop groupoid if each $(T_i, *_i)$ is a groupoid.*

**DEFINITION 1.5.13:** *Let $L = \{L_1 \cup L_2 \cup \ldots \cup L_N, *_1, \ldots, *_N\}$ be a N-loop groupoid. A non empty subset $S = \{S_1 \cup S_2 \cup \ldots \cup S_N, *_1, \ldots, *_N\}$ is said to be a sub N-loop if each $\{S_i, *_i\}$ is a loop.*

**DEFINITION 1.5.14:** *Let* $L = \{L_1 \cup L_2 \cup \ldots \cup L_N, *_1, \ldots, *_N\}$ *be a N-loop groupoid. A non empty subset* $W = \{W_1 \cup W_2 \cup \ldots \cup W_N, *_1, \ldots, *_N\}$ *of L said to be a sub N-semigroup if each* $\{W_i, *_i\}$ *is a semigroup.*

**DEFINITION 1.5.15:** *Let* $L = \{L_1 \cup L_2 \cup \ldots \cup L_N, *_1, \ldots, *_N\}$ *be a N-loop groupoid. Let* $R = \{R_1 \cup R_2 \cup \ldots \cup R_N, *_1, \ldots, *_N\}$*be a proper subset of L. We call R a sub N-group groupoid of the N-loop groupoid L if each* $\{R_i, *_i\}$ *is either a group or a groupoid.*

**DEFINITION 1.5.16:** *Let* $L = \{L_1 \cup L_2 \cup \ldots \cup L_N, *_1, \ldots, *_N\}$ *be a N-loop groupoid of finite order.* $K = \{K_1 \cup K_2 \cup K_3, *_1, \ldots, *_N\}$ *be a sub N-loop groupoid of L. If every sub N-loop groupoid divides the order of the N-loop groupoid L then we call L a Lagrange N-loop groupoid. If no sub N-loop groupoid of L divides the order of L then we say L is a Lagrange free N-loop groupoid.*

**DEFINITION 1.5.17**: *Let* $L = \{L_1 \cup L_2 \cup \ldots \cup L_N, *_1, \ldots, *_N\}$ *be a N-loop groupoid. We call L a Moufang N-loop groupoid if each* $(L_i, *_i)$ *satisfies the following identities:*

i. $(xy)\,(zx) = (x\,(yz))x$.
ii. $((x\,y)z)y = x\,(y\,(2y))$.
iii. $x\,(y\,(xz)) = (xy)x)z$ *for* $x, y, z \in L_i$, $1 \leq i \leq N$.

*Thus for a N-loop groupoid to be Moufang both the loops and the groupoids must satisfy the Moufang identity.*

**DEFINITION 1.5.18:** *Let* $L = \{L_1 \cup L_2 \cup \ldots \cup L_N, *_1, \ldots, *_N\}$ *be a N-loop groupoid. A proper subset P* $(P = P_1 \cup P_2 \cup \ldots \cup P_N, *_1, \ldots, *_N)$ *of L is a normal sub N-loop groupoid of L if*

i. *If P is a sub N-loop groupoid of L.*
ii. $x_i\,P_i = P_i\,x_i$ *(where* $P_i = P \cap L_i$*).*
iii. $y_i\,(x_i\,P_i) = (y_i\,x_i)\,P_i$ *for all* $x_i, y_i \in L_i$.

*This is true for each* $P_i$ *, i.e., for* $i = 1, 2, \ldots, N$.

**DEFINITION 1.5.19:** *Let $L = \{L_1 \cup L_2 \cup \ldots \cup L_N, *_1, \ldots, *_N\}$ and $K = \{K_1 \cup K_2 \cup \ldots \cup K_N, *_1, \ldots, *_N\}$ be two N-loop groupoids such that if $(L_i, *_i)$ is a groupoid then $\{K_i, *_i\}$ is also a groupoid. Likewise if $(L_j, *_j)$ is a loop then $(K_j, *_j)$ is also a loop true for $1 \leq i, j \leq N$. A map $\theta = \theta_1 \cup \theta_2 \cup \ldots \cup \theta_N$ from L to K is a N-loop groupoid homomorphism if each $\theta_i : L_i \to K_i$ is a groupoid homomorphism and $\theta_j : L_j \to K_j$ is a loop homomorphism $1 \leq i, j \leq N$.*

**DEFINITION 1.5.20:** *Let $L = \{L_1 \cup L_2 \cup \ldots \cup L_N, *_1, \ldots, *_N\}$ be a N-loop groupoid and $K = \{K_1 \cup K_2 \cup \ldots \cup K_M, *_1, \ldots, *_M\}$ be a M-loop groupoid.*

*A map $\phi = \phi_1 \cup \phi_2 \cup \ldots \cup \phi_{N'}$ from L to K is called a pseudo N-M-loop groupoid homomorphism if each $\phi_i$: $L_t \to K_s$ is either a loop homomorphism or a groupoid homomorphism, $1 \leq t \leq N$ and $1 \leq s \leq M$, according as $L_t$ and $K_s$ are loops or groupoids respectively (we demand $N \leq M$ for if $M > N$ we have to map two or more $L_i$ onto a single $K_j$ which can not be achieved easily).*

**DEFINITION 1.5.21:** *Let $L = \{L_1 \cup L_2 \cup \ldots \cup L_N, *_1, \ldots, *_N\}$be a N-loop groupoid. We call L a Smarandache loop N-loop groupoid (S-loop N-loop groupoid) if L has a proper subset $P = \{P_1 \cup P_2 \cup \ldots \cup P_N, *_1, \ldots, *_N\}$ such that each $P_i$ is a loop i.e. P is a N-loop.*

Now we proceed on to define the mixed N-algebraic structures, which include both associative and non associative structures. Here we define them and give their substructures and a few of their properties.

**DEFINITION 1.5.22:** *Let A be a non empty set on which is defined N-binary closed operations $*_1, \ldots, *_N$. A is called as the N-group-loop-semigroup-groupoid (N-glsg) if the following conditions, hold good.*

i. *$A = A_1 \cup A_2 \cup ... \cup A_N$ where each $A_i$ is a proper subset of A (i.e. $A_i \not\subseteq A_j \not\subseteq$ or $A_j \not\subseteq A_i$ if $(i \neq j)$).*

ii. *$(A_i, *_i)$ is a group or a loop or a groupoid or a semigroup (or used not in the mutually exclusive sense) $1 \leq i \leq N$. A is a N - glsg only if the collection $\{A_1, ..., A_N\}$ contains groups, loops, semigroups and groupoids.*

**DEFINITION 1.5.23:** *Let $A = \{A_1 \cup ... \cup A_N, *_1, ..., *_N\}$ where $A_i$ are groups, loops, semigroups and groupoids. We call a non empty subset $P = \{P_1 \cup P_2 \cup ... \cup P_N, *_1, ..., *_N\}$ of A, where $P_i = P \cap A_i$ is a group or loop or semigroup or groupoid according as $A_i$ is a group or loop or semigroup or groupoid. Then we call P to be a sub N-glsg.*

**DEFINITION 1.5.24**: *Let $A = \{A_1 \cup A_2 \cup ... \cup A_N; *_1, ..., *_N\}$ be a N-glsg. A proper subset $T = \left\{T_{i_1} \cup ... \cup T_{i_K}, *_{i_1}, ..., *_{i_K}\right\}$ of A is called the sub K-group of N-glsg if each $T_{i_t}$ is a group from $A_r$ where $A_r$ can be a group or a loop or a semigroup of a groupoid but has a proper subset which is a group.*

**DEFINITION 1.5.25:** *Let $A = \{A_1 \cup A_2 \cup ... \cup A_N; *_1, ..., *_N\}$ be a N-glsg.*

*A proper subset $T = \left\{T_{i_1} \cup ... \cup T_{i_r}\right\}$ is said to be sub r-loop of A if each $T_{i_j}$ is a loop and $T_{i_j}$ is a proper subset of some $A_p$. As in case of sub K-group r need not be the maximum number of loops in the collection $A_1, ..., A_N$.*

**DEFINITION 1.5.26:** *Let $A = \{A_1 \cup A_2 \cup ... \cup A_N; *_1, ..., *_N\}$ be a N-glsg. Let $P = \{P_1 \cup P_2 \cup ... \cup P_N\}$ be a proper subset of A where each $P_i$ is a semigroup then we call P the sub u-semigroup of the N-glsg.*

**DEFINITION 1.5.27:** *Let $A = \{A_1 \cup A_2 \cup ... \cup A_N; *_1, ..., *_N\}$ be a N-glsg. A proper subset $C = \{C_1 \cup C_2 \cup ... \cup C_t\}$ of A is said to be a sub-t-groupoid of A if each $C_i$ is a groupoid.*

**DEFINITION 1.5.28:** *Let $A = \{A_1 \cup A_2 \cup \ldots \cup A_N; *_1, \ldots, *_N\}$ be a N-glsg. Suppose A contains a subset $P = P_{L_1} \cup \ldots \cup P_{L_k}$ of A such that P is a sub K-group of A. If every P-sub K-group of A is commutative we call A to be a sub-K-group commutative N-glsg.*

*If atleast one of the sub-K-group P is commutative we call A to be a weakly sub K-group-commutative N-glsg. If no sub K-group of A is commutative we call A to be a non commutative sub-K-group of N-glsg.*

For more about these notions please refer [50].

**DEFINITION 1.5.29:** *Let $(G = G_1 \cup G_2 \cup \ldots \cup G_N, *_1, \ldots, *_N)$ be a non empty set on which is defined N-binary operations. $(G, *_1, \ldots, *_N)$ is defined as a Smarandache N-group semigroup (S-N-group semigroup) if*

- *i.* $G = G_1 \cup G_2 \cup \ldots \cup G_N$ *where each* $G_i$ *is a proper subset of G* ($G_i \not\subseteq G_j$, $G_j \not\subseteq G_j$, $i \neq j$).
- *ii.* *Some of the* $(G_i, *_i)$ *are groups.*
- *iii.* *Some of the* $(G_j, *_j)$ *are S-semigroups.*
- *iv.* *Some of the* $(G_k, *_k)$ *are just semigroups* ($1 \leq i, j, k \leq N$).

**DEFINITION 1.5.30:** *Let $G = \{G_1 \cup G_2 \cup \ldots \cup G_N, *_1, \ldots, *_N\}$ be a N-group semigroup. G be a S-N-group semigroup. We call G a Smarandache commutative N-group semigroup (S-commutative N-group semigroup) if $G_i$ are commutative groups, every proper subgroup of the S-semigroup is commutative and every semigroup which are not S-semigroups in G is also commutative.*

**DEFINITION 1.5.31:** *Let $G = \{G_1 \cup G_2 \cup \ldots \cup G_N, *_1, \ldots, *_N\}$ be a N-group semigroup. We call G a Smarandache subcommutative N-group semigroup (S-subcommutative N-group semigroup) if the following conditions are satisfied.*

- *i.* *G has non trivial sub N-groups.*

ii. *Every sub N-group of G is commutative.*

**DEFINITION 1.5.32:** *Let $G = \{G_1 \cup G_2 \cup \ldots \cup G_N, *_1, \ldots, *_N\}$ be a N-group semigroup. We say G is Smarandache cyclic N-group (S-cyclic N-group) if $\{G_i, *_i\}$ are cyclic groups and $\{G_j, *_j\}$ is a S-cyclic semigroup. We say G to be a Smarandache weakly cyclic N-group (S-weakly cyclic N-group) if every subgroup of $(G_i, *_i)$ are cyclic and $(G_j, *_j)$ is a S-weakly cyclic semigroup.*

**DEFINITION 1.5.33:** *Suppose $G = \{G_1 \cup G_2 \cup \ldots \cup G_N, *_1, \ldots, *_N\}$ be a S-N-group. The Smarandache hyper N-group semigroup (S- hyper N-group semigroup) P is defined as follows.*

i. *If $(G_i, *_i)$ are groups then they have no subgroups.*
ii. *If $(G_j, *_j)$ is a S-semigroup then $(A_j, *_j)$ is the largest subgroup of $(G_j, *_j)$.*
iii. *If $(G_k, *_k)$ are semigroups then $(T_k, *_k)$ are the largest ideals of $(G_k, *_k)$, where $P = \bigcup_{\substack{over\\relevant\\i}} (G_i, *_i) \bigcup_{\substack{over\\relevant\\j}} (A_j, *_j) \bigcup_{\substack{over\\relevant\\K}} (T_k, *_k)$.*

*Thus $P = P_1 \cup P_2 \cup \ldots \cup P_N$ where $P_i = (G_i, *_i)$ or $P_j = (A_j, *_j)$ or $P_k = (T_k, *_k)$. $1 \leq i, j, k \leq N$.*

*We call a S-N group semigroup G to be Smarandache simple N-group semigroup (S-simple N-group semigroup) if G has no S-hyper N-group semigroup.*

**DEFINITION 1.5.34:** *Let $G = \{G_1 \cup G_2 \cup \ldots \cup G_N, *_1, \ldots, *_N\}$ be a S-N-group semigroup. Let $P = \{P_1 \cup P_2 \cup \ldots \cup P_N, *_1, \ldots, *_N\}$ be a subset of G. If P itself is a S-N-group semigroup then we call P a Smarandache sub N-group semigroup (S- sub N-group semigroup) of G. Now it may happen in case of finite N-group semigroup G, $o(P) / o(G)$ or $o(P) \nmid o(G)$.*

**DEFINITION 1.5.35:** *Let $G = \{G_1 \cup G_2 \cup \ldots \cup G_N, *_1, \ldots, *_N\}$ be a S-N-group semigroup. A proper subset $P = \{P_1 \cup P_2 \cup P_3 \cup \ldots \cup P_N, *_1, \ldots, *_N\}$ is called the Smarandache P-hyper sub-N-* *group semigroup (S-P-hyper sub-N-group semigroup) if each $P_i \subset G_i$ is a maximal subgroup of $G_i$ and $P_i$ is a maximal subgroup of $G_i$, if $G_i$ is a group; or if $P_j \subset G_j$ is the largest group if $G_j$ is a S-semigroup.*

**DEFINITION 1.5.36:** *Let $G = \{G_1 \cup G_2 \cup G_3 \cup \ldots \cup G_N, *_1, \ldots, *_N\}$ be a S-N-group semigroup. Let $H = \{H_1 \cup H_2 \cup \ldots \cup H_N, *_1, \ldots, *_N\}$ be a S-sub N group semigroup of G, we define Smarandache right coset (S-right coset) $Ha = \{H_i a \cup H_1 \cup \ldots \cup \hat{H}_i \cup \ldots \cup H_N \mid$ if $a \in G_i$ and $\hat{H}_i$ denotes $H_i$ is not present in the union; this is true for any i, i.e. whenever $a \in G_i\}$.*

*Suppose a is in some r number of $G_i$'s (say) $G_i, G_j, \ldots, G_t, G_K$ then we have $Ha = \{H_1 \cup \ldots \cup a \cup H_j a \cup \ldots \cup H_i a \cup H_K a \cup \ldots \cup H_N\}$ if we assume $H_i$ are subgroups. If $H_i$'s are just S-semigroups then we take $H'_i$ to be the subgroup of the S-semigroup $H_i$.*

**DEFINITIONS 1.5.37:** *Let $G = \{G_1 \cup G_2 \cup \ldots \cup G_N, *_1, \ldots, *_N\}$ be a S-N- group semigroup. Suppose A and B be sub N groups of G. Then we define the double coset of x in G with respect to A, B if $A\,x\,B = \{A_1 x_1 B_1 \cup \ldots \cup A_N x_N B_N \mid x = x_1 \cup x_2 \cup \ldots \cup x_N\}$.*

**DEFINITION 1.5.38:** *Let $G = \{G_1 \cup G_2 \cup \ldots \cup G_N, *_1, \ldots, *_N\}$ be a S-group semigroup. We call a S-sub N-group H of S-N-group semigroup G where $A = \{A_1 \cup A_2 \cup \ldots \cup A_N, *_1, \ldots, *_N\}$ to be a Smarandache normal sub N-group semigroup (S-normal sub N-group semigroup) of G, if $x \in G$ and $x = (x_1 \cup x_2 \cup \ldots \cup x_N) \in G$ then $x\,A = (x_1 A_1 \cup \ldots \cup x_N A_N) \subseteq A$ and $Ax = Ax_1 \cup \ldots \cup A_N x_N \subseteq A$ for all $x \in G$.*

**DEFINITION 1.5.39:** *Let $L = \{L_1 \cup L_2 \cup \ldots \cup L_N, *_1, \ldots, *_N\}$ be a N-loop groupoid. We call L a Smarandache group N-loop groupoid (S-group N-loop groupoid) if L has a proper subset P where $P = \{P_1 \cup P_2 \cup \ldots \cup P_N, *_1, \ldots, *_N\}$ is a N-group i.e., each $P_i$ is a group and $P_i \subset L_i$; $P_i = P \cap L_i$.*

**DEFINITION 1.5.40:** *Let* $L = \{L_1 \cup L_2 \cup \ldots \cup L_N, *_1, \ldots, *_N\}$ *be a N-loop groupoid. We call L a Smarandache N-loop groupoid (S-N-loop groupoid) if each* $L_i$ *is either a S-loop or a S-groupoid. i.e. L has a proper subset set* $P = \{P_1 \cup P_2 \cup \ldots \cup P_N, *_1, \ldots, *_N\}$ *such that each* $P_i = P \cap L_i$ $(i = 1, 2, \ldots, N)$ *is either a group or a semigroup.*

**DEFINITION 1.5.41:** *Let* $L = \{L_1 \cup L_2 \cup \ldots \cup L_N, *_1, \ldots, *_N\}$ *be a N-loop groupoid. Let* $P = \{P_1 \cup P_2 \cup \ldots \cup P_N, *_1, \ldots, *_N\}$ *be a proper subset of L and P be a sub N-loop groupoid, if P has a proper subset* $T = \{T_1 \cup T_2 \cup \ldots \cup T_N, *_1, \ldots, *_N\}$ *such that each* $T_i$ *is a group or a semigroup then P is called a Smarandache sub N-loop groupoid (S-sub N-loop groupoid) i.e. if P is a sub N-loop groupoid then P must be S-N loop groupoid for P to be a S-sub N-loop groupoid.*

**DEFINITION 1.5.42:** *Let* $L = \{L_1 \cup L_2 \cup \ldots \cup L_N, *_1, \ldots, *_N\}$ *be a N-loop groupoid. We call L a Smarandache groupoid N-loop groupoid (S-groupoid N-loop groupoid) if L contains a proper subset* $Y = \{Y_1 \cup Y_2 \cup \ldots \cup Y_N, / *_1, \ldots, *_N\}$ *such that each* $Y_i$ *is a S-N groupoid.*

**DEFINITION 1.5.43:** *Let* $L = \{L_1 \cup L_2 \cup \ldots \cup L_N, *_1, \ldots, *_N\}$ *and* $K = \{K_1 \cup K_2 \cup \ldots \cup K_N, *_1, \ldots, *N\}$ *be S-N-loop groupoids. A map* $\phi = \{\phi_1 \cup \ldots \cup \phi_N\}$ *from L to K is called the Smarandache homomorphism of N-loop groupoids (S- homomorphism of N-loop groupoids) if each* $\phi_i$ *is a group homomorphism of a semigroup homomorphism according as* $L_i$ *is a S-loop or a S-groupoid respectively.*

**DEFINITION 1.5.44:** *Let* $L = \{L_1 \cup L_2 \cup \ldots \cup L_N, *_1, \ldots, *_N\}$ *be a N-loop groupoid suppose* $L_{i_1}, \ldots, L_{i_K}$ *be K groupoid in the set* $\{L_1, L_2, \ldots, L_N\}$*. We call L a Smarandache groupoid N-loop groupoid (S-groupoid N-loop groupoid) only if each groupoid* $\{L_{i_t} \mid 1 \leq t \leq K\}$ *is a S-groupoid.*

**DEFINITION 1.5.45:** *Let* $L = \{L_1 \cup L_2 \cup \ldots \cup L_N, *_1, \ldots, *_N\}$ *be a S-groupoid N-loop groupoid. A proper subset of the collection* *of all groupoids G =* $\left\{L_{i_1} \cup ... \cup L_{i_K}\right\}$ *of* $\{L_1 \cup ... \cup L_N\}$ *is said to be Smarandache subgroupoid N-loop groupoid (S-subgroupoid N-loop groupoid) if G has a proper subset T =* $\left\{T_{i_1} \cup T_{i_2} \cup ... \cup T_{i_K}\right\}$ *such that T is itself a Smarandache sub K-groupoid.*

**DEFINITION 1.5.46:** *Let* $L = \{L_1 \cup L_2 \cup ... \cup L_N, *_1, ..., *_N\}$ *be a S-groupoid N-loop groupoid. Let G =* $\left\{L_{i_1} \cup ... \cup L_{i_K}\right\}$ *be the collection of groupoids in L. A non empty proper subset P =* $\left\{P_{i_1} \cup ... \cup P_{i_K}\right\}$ *of the K-groupoid G is said to be a K-left ideal of the K-groupoid G if*

i. *P is a sub K-groupoid.*
ii. *For all* $x_i \in L_{i_t}$ *and* $a_i \in P_{i_t}$, $x_i a_i \in P_{i_t}$, *t = 1, 2, ..., K.*

*One can similarly define K-right ideal of the K-groupoid. We say P is an K-ideal of the K-groupoid G if P is simultaneously a K-left and a K-right ideal of G.*

**DEFINITION 1.5.47:** *Let* $L = \{L_1 \cup L_2 \cup ... \cup L_N, *_1, ..., *_N\}$ *be a N-loop groupoid. Let G =* $G_{i_1} \cup ... \cup G_{i_K}$ *be the K-groupoid; i.e.* $G_{i_1}, ..., G_{i_K}$ *be the collection of all groupoids in the set* $\{L_1, ..., L_N\}$*. A sub K-groupoid V =* $\left\{V_{i_1} \cup V_{i_2} \cup ... \cup V_{i_K}\right\}$ *of G is said to be normal K subgroupoid of G if*

$$a_i V_{i_t} = V_{i_t} a_i$$
$$(V_{i_t} x_i) y_i = V_{i_t} (x_i y_i)$$
$$y_i (x_i V_{i_t}) = (y_i x_i) V_{i_1}$$

*for all* $a_i, x_i, y_i \in G_{i_t}$. *The S-groupoid N-loop groupoid is K-simple if it has no nontrivial normal K-subgroupoids.*

**DEFINITION 1.5.48:** *Let L = {$L_1 \cup L_2 \cup \ldots \cup L_N$, $*_1$, …, $*_N$} be a S-groupoid N-loop groupoid. Let G = $G_{i_1} \cup \ldots \cup G_{i_K}$ be the K-groupoid (i.e. $G_{i_1}, \ldots, G_{i_K}$ be the collection of all groupoids from the set $L_1$, …, $L_N$). We call L a Smarandache normal groupoid N-loop groupoid (S-normal groupoid N-loop groupoid) if*

*i.* $xG = Gx$ *where*

$X = x_{i_1} \cup \ldots \cup x_{i_K}$ *and*

$G = G_{i_1} \cup G_{i_2} \cup \ldots \cup G_{i_K}$ *i.e.*

$$xG = \left(x_{i_1} G_{i_1} \cup \ldots \cup x_{i_K} G_{i_K}\right) = \left(G_{i_1} x_{i_1} \cup \ldots \cup G_{i_K} x_{i_K}\right) = Gx.$$

*ii.* $G(xy) = (Gx)\, y$; $x, y \in G$ *where*

$x = x_{i_1} \cup \ldots \cup x_{i_K}$ *and*

$y = y_{i_1} \cup \ldots \cup y_{i_K}$.

$$\begin{aligned} G(xy) &= \left(G_{i_1} \cup \ldots \cup G_{i_K}\right)\left(x_{i_1} \cup \ldots \cup x_{i_K}\right)\left(y_{i_1} \cup \ldots \cup y_{i_K}\right) \\ &= \left(G_{i_1} (x_{i_1} y_{i_1}) \cup \ldots \cup G_{i_K} (x_{i_K} y_{i_K})\right) \\ &= \left\{\left(G_{i_1} x_{i_1}\right) y_{i_1} \cup \ldots \cup \left(G_{i_K} x_{i_K}\right) y_{i_K}\right\} \\ &= (Gx)\, y. \end{aligned}$$

*iii.* $y\,(xG) = (yx)\, G$ *where*

$y = y_{i_l} \cup \ldots \cup y_{i_K}$ *and*

$x = x_{i_l} \cup \ldots \cup x_{i_K}$

$$\begin{aligned} y(xG) &= \left(y_{i_1} \cup \ldots \cup y_{i_K}\right)\left[\left(x_{i_l} \cup \ldots \cup x_{i_K}\right)\left(G_{i_1} \cup \ldots \cup G_{i_K}\right)\right] \\ &= y_{i_1}\left(x_{i_1} G_{i_1}\right) \cup y_{i_2}\left(x_{i_2} G_{i_2}\right) \cup \ldots \cup y_{i_K} (x_{i_K} G_{i_K}) \\ &= \left(y_{i_1} x_{i_1}\right) G_{i_1} \cup \left(y_{i_2} x_{i_2}\right) G_{i_2} \cup \ldots \cup (y_{i_K} x_{i_K}) G_{i_K} \end{aligned}$$

$$= (yx)G$$

*for all x, y $\in$ G.*

**DEFINITION 1.5.49:** *Let L = {$L_1 \cup L_2 \cup \ldots \cup L_N$, $*_1, \ldots, *_N$} be a S-groupoid N-loop groupoid. Suppose H and P be two S-K-subgroupoids of*

$$G = \{ L_{i_1} \cup L_{i_2} \cup \ldots \cup L_{i_K} \}$$

*where $\{L_{i_1}, \ldots, L_{i_K}\}$ is the collection of all groupoids from the set {$L_1, L_2, \ldots, L_N$} we say H and P are Smarandache conjugate if there exists*

$$x = \left(x_{i_1} \cup \ldots \cup x_{i_K}\right) \in H$$

*such that H = xP or Px where*

$H = \left\{H_{i_1} \cup \ldots \cup H_{i_K}\right\}$ *and*

$P = \left\{P_{i_1} \cup \ldots \cup P_{i_K}\right\}$ *and*

$H = xP$

$$\text{so } H = \left\{H_{i_1} \cup \ldots \cup H_{i_K}\right\} = \left(x_{i_1} \cup \ldots \cup x_{i_K}\right)\left(P_{i_1} \cup \ldots \cup P_{i_K}\right).$$

$$= \left\{x_{i_1} P_{i_1} \cup \ldots \cup P_{i_K} x_{i_K}\right\} = \left\{P_{i_1} x_{i_1} \cup \ldots \cup P_{i_K} x_{i_K}\right\}$$

*Clearly $H \cap P = \phi$.*

**DEFINITION 1.5.50:** *Let A = {$A_1 \cup A_2 \cup \ldots \cup A_N$, $*_1, \ldots, *_N$} be a N-glsg. A is said to be a Smarandache N-glsg (S-N-glsg) if some $A_i$ are S-loops, some of $A_j$ are S-semigroups and some of $A_K$ are S-groupoids.*

**DEFINITION 1.5.51:** *Let A = {$A_1 \cup A_2 \cup \ldots \cup A_N$, $*_1, *_2, \ldots, *_N$} be a N-glsg. A proper subset P = {$P_1 \cup P_2 \cup \ldots \cup P_N$, $*_1, \ldots, *_N$} of A is said to be a Smarandache sub N-glsg if P itself is a S-N-glsg under the operations of A.*

**DEFINITION 1.5.52:** *Let* $\{\langle M \cup I \rangle = M_1 \cup M_2 \cup \ldots \cup M_N, *_1, \ldots, *_N\}$, $(N \geq 5)$ *we call* $\langle M \cup I \rangle$ *a mixed neutrosophic N-structure if*

i. $\langle M \cup I \rangle = M_1 \cup M_2 \cup \ldots \cup M_N$, *each* $M_i$ *is a proper subset of* $\langle M \cup I \rangle$.
ii. *Some of* $(M_i, *_i)$ *are neutrosophic groups.*
iii. *Some of* $(M_j, *_j)$ *are neutrosophic loops.*
iv. *Some of* $(M_k, *_k)$ *are neutrosophic groupoids.*
v. *Some of* $(M_r, *_r)$ *are neutrosophic semigroups.*
vi. *Rest of* $(M_t, *_t)$ *can be loops or groups or semigroups or groupoids. ('or' not used in the mutually exclusive sense*

*(From this the assumption* $N \geq 5$ *is clear).*

**DEFINITION 1.5.53:** $\{\langle D \cup I \rangle = D_1 \cup D_2 \cup \ldots \cup D_N, *_1, *_2, \ldots, *_N\}$, $N \geq 5$ *be a non empty set on which is defined N-binary operations.*

*We say* $\langle D \cup I \rangle$ *is a mixed dual neutrosophic N-structure if the following conditions are satisfied.*

i. $\langle D \cup I \rangle = D_1 \cup D_2 \cup \ldots \cup D_N$ *where each* $D_i$ *is a proper subset of* $\langle D \cup I \rangle$
ii. *For some i,* $(D_i, *_i)$ *are groups*
iii. *For some j,* $(D_j, *_j)$ *are loops*
iv. *For some k,* $(D_k, *_k)$ *are semigroups*
v. *For some t,* $(D_t, *_t)$ *are groupoids.*
vi. *The rest of* $(D_m, *_m)$ *are neutrosophic groupoids or neutrosophic groups or neutrosophic loops or neutrosophic semigroup 'or' not used in the mutually exclusive sense.*

**DEFINITION 1.5.54:** *Let* $\langle W \cup I \rangle = \{W_1 \cup W_2 \cup \ldots \cup W_N, *_1, *_2, \ldots, *_N\}$ *be a non empty set with N-binary operations* $*_1, \ldots, *_N$. $\langle W \cup I \rangle$ *is said to be a weak mixed neutrosophic structure, if the following conditions are true.*

i. $\langle W \cup I \rangle = W_1 \cup W_2 \cup \ldots \cup W_N$ *is such that each* $W_i$ *is a proper subset of* $\langle W \cup I \rangle$.
ii. *Some of* $(W_i, *_i)$ *are neutrosophic groups or neutrosophic loops.*
iii. *Some of* $(W_j, *_j)$ *are neutrosophic groupoids or neutrosophic semigroups*
iv. *Rest of* $(W_k, *_k)$ *are groups or loops or groupoids or semigroups. i.e. In the collection* $\{W_i, *_i\}$ *all the 4 algebraic neutrosophic structures may not be present.*

*At most 3 algebraic neutrosophic structures are present and atleast 2 algebraic neutrosophic structures are present. Rest being non neutrosophic algebraic structures.*

**DEFINITION 1.5.55:** *Let* $\{\langle V \cup I \rangle = V_1 \cup V_2 \cup \ldots \cup V_N, *_1, \ldots, *_N\}$ *be a non empty set with N-binary operations. We say* $\langle V \cup I \rangle$ *is a weak mixed dual neutrosophic N-structure if the following conditions are true.*

i. $\langle V \cup I \rangle = V_1 \cup V_2 \cup \ldots \cup V_N$ *is such that each* $V_i$ *is a proper subset of* $\langle V \cup I \rangle$
ii. *Some of* $(V_i, *_i)$ *are loops or groups*
iii. *Some of* $(V_j, *_j)$ *are groupoids or semigroups*
iv. *Rest of the* $(V_k, *_k)$ *are neutrosophic loops or neutrosophic groups or neutrosophic groupoids or neutrosophic semigroups.*

**DEFINITION 1.5.56:** *Let* $\{\langle M \cup I \rangle = M_1 \cup M_2 \cup \ldots \cup M_N, *_1, \ldots, *_N\}$ *where* $\langle M \cup I \rangle$ *is a mixed neutrosophic N-algebraic structure. We say a proper subset* $\{\langle P \cup I \rangle = P_1 \cup P_2 \cup \ldots \cup P_N, *_1, \ldots, *_N\}$ *is a mixed neutrosophic sub N-structure if* $\langle P \cup I \rangle$ *under the operations of* $\langle M \cup I \rangle$ *is a mixed neutrosophic N-algebraic structure.*

**DEFINITION 1.5.57:** *Let* $\langle W \cup I \rangle = \{W_1 \cup W_2 \cup \ldots \cup W_N, *_1, *_2, \ldots, *_N\}$ *be a mixed neutrosophic N-structure.*

*We call a finite non empty subset P of* $\langle W \cup I \rangle$, *to be a weak mixed deficit neutrosophic sub N-structure if* $P = \{P_1 \cup P_2 \cup \ldots$ $\cup P_t$, $*_1$, …, $*_t\}$, *$1 < t < N$ with $P_i = P \cap L_k$, $1 \leq i \leq t$ and $1 \leq k \leq N$ and some $P_i$'s are neutrosophic groups or neutrosophic loops some of the $P_j$'s are neutrosophic groupoids or neutrosophic semigroups and rest of the $P_k$'s are groups or loops or groupoids or semigroups.*

For more about the concepts defined please refer [41-51].

**Chapter two**

# SMARANDACHE NEUTROSOPHIC GROUPS AND THEIR PROPERTIES

This chapter has three sections. Study of Smarandache neutrosophic groups is interesting and innovative. So in section one the notion of S-neutrosophic groups are introduced. Also the notion of pseudo S- neutrosophic group is defined. Section two introduces the notion of Smarandache neutrosophic bigroups and the exceptional properties enjoyed by them are brought out.

Examples are given for the better understanding of the problems. The Smarandache neutrosophic N-groups are introduced for the first time and they happen to help in getting a generalized situation, for in these S-neutrosophic N-structure S-semigroups and S-neutrosophic semigroups are also taken. There by make the notion of S-neutrosophic N-groups a most generalized structure [51].

## 2.1 Smarandache neutrosophic groups and their properties

In this section for the first time the new notion of Smarandache neutrosophic group is defined. It is important to mention here that we do not have the concept of Smarandache groups. Substructures like Smarandache neutrosophic subgroup and Smarandache right cosets are defined. However in general the S-right (left) coset do not partition the S-neutrosophic group.

Now we proceed on to define the notion of Smarandache Neutrosophic groups.

**DEFINITION 2.1.1:** *Let N(G) be a neutrosophic group. We say N (G) is a Smarandache neutrosophic group (S-neutrosophic group) if N(G) has a proper subset P where P is a pseudo neutrosophic subgroup of N(G).*

All neutrosophic groups need not be S-neutrosophic groups. For they may not contain a pseudo neutrosophic group in them.

Now a S-neutrosophic group is said to be of finite order if the number of distinct elements in them is finite and it is of infinite order if it has infinite number of elements in them.

**DEFINITION 2.1.2:***Let (N(G), *) be a neutrosophic group a proper subset P of N(G) is said to be a Smarandache neutrosophic subgroup (S-neutrosophic subgroup) if (P, *) is a neutrosophic subgroup and has a proper subset T such that T is a pseudo neutrosophic group under the same binary operation *.*

*Note:* It is a matter of routine to show if a neutrosophic group N(G) has a S-neutrosophic subgroup P then N(G) is itself a S-neutrosophic group.

Now we proceed on to define strong Smarandache neutrosophic group.

**DEFINITION 2.1.3:** *Let N(G) be a neutrosophic group. Suppose N(G) has a proper subset P such that P is a neutrosophic subgroup of N(G) then we call N(G) to be a strong Smarandache neutrosophic group (strong S-neutrosophic group) .*

Also we define the notion of pseudo Smarandache neutrosophic group.

**DEFINITION 2.1.4:** *Let N(G) be a neutrosophic group. Let P be a proper subset of N (G) and if P is a Smarandache neutrosophic semigroup then we call N (G) to be a pseudo Smarandache neutrosophic group (pseudo S-neutrosophic group).*

The interested reader can analyze and find out whether there exists any form of relation among these 3 types of Smarandache neutrosophic groups.

For the notion of S-neutrosophic semigroup refer chapter 3.

Now we proceed on to briefly define some more concepts in finite S-neutrosophic groups before we go for the definition of neutrosophic bigroups.

**DEFINITION 2.1.5:** *Let N(G) be a finite neutrosophic group. A proper subset P of N(G) which is a S-neutrosophic subgroup of N(G) is said to be Lagrange Smarandache neutrosophic subgroup (Lagrange S-neutrosophic subgroup) if, o(P) / o(N(G)).*

*If all S-neutrosophic subgroups are Lagrange S-neutrosophic subgroups then we call N(G) a Lagrange Smarandache neutrosophic group (Lagrange S-neutrosophic group).*

*If N(G) has at least one Lagrange neutrosophic subgroup then we call N(G) a weakly Lagrange Smarandache neutrosophic group (weakly Lagrange S-neutrosophic group). If N(G) has no Lagrange S-neutrosophic subgroup then we call N(G) to be a Lagrange free Smarandache neutrosophic group (Lagrange free S-neutrosophic group).*

We illustrate these by the following example.

***Example 2.1.1:*** Consider N(G) = {0, 1, 2, 3, 4, 5, 6, I, 2I, 3I, 4I 5I, 6I}, a neutrosophic group under multiplication modulo 7.

Clearly N(G) is not a S-neutrosophic group. But N(G) is a strong S-neutrosophic group for P = {1, I, 6, 6I} is a neutrosophic subgroup but N(G) has no proper subset which is a pseudo neutrosophic subgroup. Clearly o(P) $\nmid$ o(N(G)), since the order of N(G) is a prime we can declare N(G) to be Lagrange free strong S-group.

Now we proceed on to give yet another example.

***Example 2.1.2:*** Let N(G) be a finite S-neutrosophic group where N(G) = {1, 2, 3, 4, I, 2I, 3I, 4I, 1 + I, 2 + I, 3 + I, 4 + I, 1 + 2I, 2 +2I, 3+2I, 4 + 2I, 1+3I, 2+3I, 3+3I, 4 + 3I, 1 +4I, 2 + 4I, 3 + 4I, 4 + 4I}. N(G) is a S-neutrosophic group of order 24.

P = {1, I, 4I} is a pseudo neutrosophic subgroup of N(G), o(P) / o(N(G)).

L = {1, I, 4, 4I} is a neutrosophic subgroup of N(G), o(L) / o(N(G)).

Also T = {1, 1 + 3I} is a pseudo neutrosophic subgroup of N(G) and o(T) / o(N(G)). Thus this S-neutrosophic group has both S-neutrosophic subgroups and pseudo S-neutrosophic subgroups.

Now we proceed on to define the notion of Cauchy elements in a Smarandache neutrosophic element.

**DEFINITION 2.1.6:** *Let N(G) be a S-neutrosophic group of finite order if every torsion element x of N(G) and every neutrosophic element n of N(G) are such that if $x^t = 1$ and $n^m = I$ and if t / o(N(G)) and m/o(N(G)) then we call N(G) to be a Cauchy Smarandache neutrosophic group (Cauchy S-neutrosophic group).*

*The element x is called as the Cauchy Smarandache element (Cauchy S-element) of N(G) and the element n is called the Cauchy Smarandache neutrosophic element (Cauchy S-neutrosophic element) of N(G).*

*We can have elements in N(G) which are torsion elements y in N(G) which are such that $y^r = 1$, r $\nmid$ o(N(G)); also neutrosophic elements p in N(G) with $p^s = I$ still, s $\nmid$ o(N(G)).*

*These elements will not be termed as Cauchy elements. If every element of N(G) is either a Cauchy S-element or Cauchy*

*S-neutrosophic element then we call N(G) to be a Cauchy S-neutrosophic group.*

*If N(G) has atleast one Cauchy S-element and Cauchy S-neutrosophic element then we call N(G) to be a weakly Cauchy Smarandache neutrosophic group (weakly Cauchy S-neutrosophic group). If N(G) has no Cauchy elements N(G) is called as Cauchy free Smarandache neutrosophic group (Cauchy free S-neutrosophic group).*

We can have several examples of these, we proceed onto define p-Sylow Smarandache neutrosophic group and strong p-Sylow Smarandache neutrosophic group.

**DEFINITION 2.1.7:** *Let N(G) be a finite S-neutrosophic group. If for a prime p such that $p^{\alpha}$ / o(N(G)) and $p^{\alpha+1} \nmid$ o(N(G)), N(G) has a S-neutrosophic subgroup P of order $p^{\alpha}$ then, P is called the p-Sylow Smarandache neutrosophic subgroup (p-Sylow S-neutrosophic subgroup) of N(G). If for every prime p such that $p^{\alpha}$ / o(N(G)) and $p^{\alpha+1} \nmid$ o(N(G)), we have p-Sylow S-neutrosophic subgroup then we call N(G) a Sylow strong Smarandache neutrosophic group (Sylow strong S-neutrosophic group).*

*If for p a prime $p^{\alpha}$ / o(N(G)) and $p^{\alpha+1} \nmid$ o(N(G)). N(G) has a pseudo neutrosophic subgroup L of order $p^{\alpha}$ then L is called the p-pseudo Sylow Smarandache neutrosophic subgroup (p-pseudo Sylow S-neutrosophic subgroup) of N(G). If for every prime p, such that $p^{\alpha}$ / o(N(G)) and $p^{\alpha+1} \nmid$ o(N(G)) we have p-pseudo Sylow neutrosophic subgroup then we call N(G) a pseudo Sylow Smarandache neutrosophic group (pseudo Sylow S-neutrosophic group).*

*If on the other hand N(G) has at least one p-pseudo Sylow neutrosophic subgroup then we call N(G) to be a weak Sylow Smarandache neutrosophic group (weak Sylow S-neutrosophic group).*

On similar lines we can define weak Sylow strong Smarandache neutrosophic group.

As in case of neutrosophic groups [51] we can also in the case of S-neutrosophic groups define the notion of Smarandache coset and pseudo Smarandache coset.

**DEFINITION 2.1.8:** *Let N(G) be a S-neutrosophic group. Let P be a pseudo neutrosophic subgroup of N(G).*

*The pseudo Smarandache right coset (pseudo S-right coset) of P is defined as Pa = {pa | p ∈ P and a ∈ N(G)}). If L is a S-neutrosophic subgroup of N(G) we for a ∈ N(G) define the Smarandache right coset (S-right coset) of L as La = {la / l ∈ L and a ∈ N(G)}.*

Now we illustrate these situations by the following example.

***Example 2.1.3:*** Let N(G) = {0, 1, 2, 3, 4 I, 2I, 3I, 4I, 1 + I, 2 + I, 3 + I, 4 + I, 1 + 2I, 1 + 3I, 1 + 4I, 2 + 2I, 2 + 3I, 2 + 4I, 3 + 2I, 3 + 3I, 3 + 4I, 4+2I, 4+ 3I, 4+4I} be a S-neutrosophic group and o (N(G)) = 25. Let P = {1, I, 4I} which is a pseudo neutrosophic subgroup of N(G), we will find the pseudo right cosets; P = {I, 1, 4I}, now

| | | |
|---|---|---|
| P. P | = | {1, I, 4I} |
| P. 0 | = | {0} |
| P.1 | = | {1, I, 4I} |
| P.I | = | {I, 4I} |
| P.4I | = | {4I} |
| P.2 | = | {2, 2I, 3I} |
| P.3 | = | {3, 3I, 2I} |
| P.4 | = | {4, 4I} |
| P.2I | = | {2I, 3I} |
| P. 3I | = | {3I, 2I} |
| | = | P.2I |
| P (1 + I) | = | {1 + I, 2I, 3I} |
| P (2 + I) | = | {2 + I, 3I, 2I} |
| P (3 + 1) | = | {3 + I, 4I 3I} |
| P (4 + I) | = | {4 + I, 0} |
| P (2I + 1) | = | {2I + 1, 3I, 2I} |
| P (2I + 2) | = | {2I + 2, 4I, I} |
| P (2I + 3) | = | {2I + 3, 0} |

P (2I + 4) = {2I + 4, 4I}
P (3I + 1) = {3I + 1, 4I, I}
P (3I + 2) = {3I + 2, 0}
P (3I + 3) = {3I+3, I, 4I}
P (3I + 4) = {3I + 4, 2I, 3I}
P (4I + 1) = {4I + 1, 0}
P (4I + 2) = {4I + 2, I, 3I}
P (4I + 3) = {4I + 3, 2I, 3I}
P (4I + 3) = {4I + 3, 2I, 3I}
P (4I + 4) = {4I + 4, 3I, 2I}.

Clearly we see the S-right coset of a pseudo neutrosophic subgroup in general does not partition N(G).

Consider M = {1, I, 4I, 4}, M is a S-neutrosophic subgroup of N (G). Now we find their S-right coset values.

M. 0 = {0}
M. 1 = {1, I, 4, 4I}
M.I = {I, 4I}
M.r = {4, 4I}
M 4 I = {4I, I }
= MI
M (1 + I) = {1+ I, 2I, 4 + 4I, 3I}
M (2 + I) = {2 + I, 3I, 3 + 4I, 2I }
M (3 + I) = {3 + I, 4I, 2 + 4I, I}
M (4 + I) = {4 + I, o, 1 + 4I}
M (1 + 2I) = {1 + 2I, 3I, 4 + 3I 2 I}
M (2 + 2I) = {2 + 2I, 4I, 3 + 3I, I}
M (3 + 2I) = {3 + 2I, 0, 2 + 3I, 3 I}
M (4 + 2I) = {4 + 2I, I, + 3I, 4I}
M (1 + 3I) = {1 + 3I, 4I, 4 + 2I, I}
M (2 + 3I) = {2 + 3I, 0, 3 + 2I}
M (3 + 3I) = {3 + 3I, I, 2 + 2I, 4I}
M (4 + 3I) = {4 + 3I, 2I, 1 + 2I, 3I}
M (1 + 4I) = {1 + 4I, 0, 4 + I}
M (2 + 4I) = {2 + 4I, I, 3 + I, 4I}
M (3 + 4I) = {3 + 4I, 2I, 2 + I, 3I}
M (4 + 4I) = { 4 + 4I, 3I, 1 + I, 2I}.

Thus we see the S-right cosets in general do not partition the S-neutrosophic group.

Interested reader can develop more results in this direction. Now we proceed on to define the notion of Smarandache neutrosophic conjugate subgroups of a S-neutrosophic group N(G).

**DEFINITION 2.1.9:** *Let N(G) be a S-neutrosophic group. Let P and K be any two S-neutrosophic subgroups of N(G), we say P and K are Smarandache conjugate if we can find x, y $\in$ N(G) such that Px = Ky (yK).*

*On similar lines if P' and K' are pseudo neutrosophic subgroups of N(G) we say they are pseudo Smarandache conjugate (pseudo S-conjugate) if P' x' = K' y' (or y' K') for some x', y' $\in$ N(G).*

The main interesting thing to note about these conjugate relation in S-neutrosophic groups is that if P and K are S-conjugate we need not have o(P) = o(K).

## 2.2 Smarandache neutrosophic Bigroups

In this section we introduce the notion of Smarandache neutrosophic bigroups. As in case of neutrosophic bigroups all substructures like S-Lagrange neutrosophic bigroup, S-Sylow neutrosophic bigroups and other results can be derived as a matter of routine exercise. Interested readers will certainly try their hand in this!

Now we proceed on to define Smarandache neutrosophic bigroup.

**DEFINITION 2.2.1:** *Let $B_N(G) = \{B(G_1) \cup B(G_2), *_1, *_2\}$, be a neutrosophic bigroup. $B_N(G)$ is defined to be a Smarandache neutrosophic bigroup (S-neutrosophic bigroup) if the neutrosophic group $B_1(G)$ is S-neutrosophic group (for S-group is not defined) and $B_2(G)$ is a group.*

Now we proceed on to define the notion of more generalized structures viz., Smarandache neutrosophic bigroup I and II.

**DEFINITION 2.2.2:** *Let $B_N(G) = \{B(G_1) \cup B(G_2), *_1, *_2\}$ be a proper subset with two binary operations if $B(G_1)$ is a S-neutrosophic group and $B(G_2)$ a S-semigroup then we call $B_N(G)$ a Smarandache neutrosophic bigroup I (S-neutrosophic bigroup I).*

*If in $B_N(G)$; $B(G_1)$ is a neutrosophic group and $B(G_2)$ a S-neutrosophic semigroup then we call $B_N(G)$ a Smarandache neutrosophic bigroup II (S-neutrosophic bigroup II) .*

We see Smarandache neutrosophic bigroup I and II are not related in general.

***Example 2.2.1:*** Let $B_N(G) = \{B(G_1) \cup B(G_2), *_1, *_2\}$ where

$B(G_1)$ = {1, 2, 3, 4, I, 2I, 3I, 4I}, a S-neutrosophic group and
$B(G_2)$ = $\{Z_{12}$, S-semigroup under multiplication modulo 12}.

$B_N(G)$ is a S-neutrosophic bigroup I.

Now we give an example of a S-neutrosophic bigroup II.

***Example 2.2.2:*** Let $B_N(G) = \{B(G_1) \cup B(G_2), *_1, *_2\}$, where

$B(G_1)$ = {0, I+1, 1}, is a neutrosophic group.
$B(G_2)$ = {1, 2, 3, 4, 5, I, 2I 3I, 4I, 5I}, S-neutrosophic semigroup under multiplication modulo 6.

$B_N(G)$ is a S-neutrosophic bigroup II.

Now we proceed on to define the substructures on S-neutrosophic bigroups.

**DEFINITION 2.2.3:** *Let $B_N(G) = \{B(G_1) \cup B(G_2), *_1, *_2\}$ be a neutrosophic bigroup. A proper subset $P_N(G) = \{P(G_1) \cup P(G_2), *_1, *_2\}$ is said to be Smarandache neutrosophic* *subbigroup (S-neutrosophic subbigroup) of $B_N(G)$ if ($P_N(G)$, $*_1$, $*_2$) is itself a S-neutrosophic bigroup. From this we have to make the following observations.*

*We see if $B_N$ (G) is a neutrosophic bigroup and has a S-neutrosophic subbigroup then we see $B_N$ (G) is itself a S-neutrosophic bigroup.*

Now we illustrate these by the following examples.

***Example 2.2.3:*** Let $B(G_1)$ = {0, 1, 2, 3, 4, I, 2I, 3I, 4I, 1 + I, 2 + I, , …, 3 + 4I, 4 + 4I} be a neutrosophic group under multiplication modulo 5 and $B(G_2) = A_4$; $B_N(G) = B(G_1) \cup B(G_2)$ is a S-neutrosophic bigroup.

$P(G_1) = \{1, I, 4, 4I\} \subset B(G_1)$ and
$P(G_2)$ =

$$\left\{\begin{pmatrix}1&2&3&4\\1&2&3&4\end{pmatrix},\begin{pmatrix}1&2&3&4\\2&1&4&3\end{pmatrix},\begin{pmatrix}1&2&3&4\\3&4&1&2\end{pmatrix},\begin{pmatrix}1&2&3&4\\4&3&2&1\end{pmatrix}\right\}$$

$P_N(G) = P(G_1) \cup P(G_2)$ is a S-neutrosophic subbigroup of $B_N(G)$.

We see in general the order of the S-neutrosophic subbigroup need not divide the order of the S-neutrosophic bigroup.

***Example 2.2.4:*** Let $B_N(G) = \{B(G_1) \cup B(G_2), *_1, *_2\}$ where

$B(G_1)$ = {0, 1, 2, 1 + I, 2 + I, 1 + 2I, 2 + 2I, I, 2I} a neutrosophic group under multiplication modulo 3.
$B(G_2)$ = $\{g \mid g^7 = e\}$, a cyclic group of order 7.

$o(B_N(G)) = 16$.

Consider $P_N(G) = P(G_1) \cup P(G_2) = \{1, 2, I, 2I\} \cup \{e\}$; $o(P_N(G)) = 5$. (5, 16) = 1.

Take $T_N(G) = T(G_1) \cup T(G_2)$, $T(G_1)$ = {1, 2, I, 2I} and $T(G_2) = \{g \mid g^7 = 1\}$. $o(T_N(G)) = 11$, (11, 16) = 1.

Clearly they are neutrosophic subbigroups and their order does not divide the order of the neutrosophic bigroup.

Now we give yet another example of a S-neutrosophic group.

***Example 2.2.5:*** Let $B_N$ (G) = {B($G_1$) ∪ B($G_2$), $*_1$, $*_2$} be a S-neutrosophic bigroup where

B($G_1$) = {0, 1, I, 1+I} is a S-neutrosophic group and
B($G_2$) = {$S_3$},

o($B_N$(G)) = 10. Let $P_N$(G) = P($G_1$) ∪ P($G_2$) where

$$P(G_1) = \{1, I, 1+I\} \cup \left\{\begin{pmatrix}1 & 2 & 3\\ 1 & 2 & 3\end{pmatrix}, \begin{pmatrix}1 & 2 & 3\\ 2 & 1 & 3\end{pmatrix}\right\} = P(G_2).$$

o($P_N$ (G)) = 5, 5 / o ($B_N$ (G).

Now we proceed on to define the notion of Lagrange Smarandache neutrosophic bigroup.

**DEFINITION 2.2.4:** *Let $B_N(G) = \{B(G_1) \cup B(G_2), *_1, *_2\}$ be a S-neutrosophic bigroup. Let $P_N(G) = \{P(G_1) \cup P(G_2), *_1, *_2\}$ be a S-neutrosophic subbigroup if $o(P_N(G))$ / $o(B_N(G))$ then $P_N(G)$ is called a Lagrange Smarandache neutrosophic subbigroup (Lagrange S-neutrosophic subbigroup) .*

*If every S-neutrosophic subbigroup is Lagrange then we call $B_N(G)$ to be a Lagrange Smarandache neutrosophic bigroup (Lagrange S-neutrosophic bigroup). If $B_N(G)$ has atleast one Lagrange Smarandache neutrosophic subbigroup (S-Lagrange Smarandache neutrosophic subbigroup) then we call $B_N(G)$ a weakly Lagrange Smarandache neutrosophic bigroup (weakly Lagrange S-neutrosophic bigroup) . If $B_N(G)$ has no Lagrange S-neutrosophic subbigroup we call $B_N(G)$ a Lagrange free Smarandache neutrosophic bigroup (Lagrange free S-neutrosophic bigroup) .*

Interested reader can derive several inter related properties about these structures. Next we proceed on to recall the notion of Smarandache neutrosophic bigroup II. From now onwards we define notions and properties for S-neutrosophic bigroup II denoted by $\langle G \cup I\rangle$.

**DEFINITION 2.2.5:** *Let $(\langle G \cup I\rangle = \langle G_1 \cup I\rangle \cup \langle G_2 \cup I\rangle, *, o)$ be a neutrosophic bigroup $(\langle G \cup I\rangle, *, o)$ is defined as a Smarandache neutrosophic bigroup II (S-neutrosophic bigroup II) if*

i. *$(\langle G_1 \cup I\rangle, *)$ is a neutrosophic group.*
ii. *$(\langle G_2 \cup I\rangle, o)$ is a S-neutrosophic semigroup.*

***Example 2.2.6:*** Let $(\langle G \cup I\rangle, *, o)$ be a neutrosophic bigroup where $(\langle Z \cup I\rangle, +)$ is a neutrosophic group and {0, 1, 2, 3, 4, 5, I, 2I, 3I, 4I, 5I} is a S-neutrosophic semigroup under multiplication modulo 6. So $(\langle G \cup I\rangle, *, o)$ is a S-neutrosophic bigroup II.

We define the notion of Smarandache neutrosophic subbigroup.

**DEFINITION 2.2.6:** *Let $(\langle G \cup I\rangle = \langle G_1 \cup I\rangle \cup \langle G_2 \cup I\rangle, *, o)$ be a neutrosophic bigroup, $H = (H_1 \cup H_2, *, o)$ is called a Smarandache neutrosophic subbigroup II (S-neutrosophic subbigroup II) if $(\langle H_1 \cup I\rangle, *)$ is a neutrosophic subgroup of $(\langle G_1 \cup I\rangle, *)$ and $(\langle H_2 \cup I\rangle, o)$ is a S-neutrosophic subsemigroup of $(\langle G_2 \cup I\rangle, o)$.*

The following theorem is interesting and the reader to expected to prove it.

**THEOREM 2.2.1:** *Let $(\langle G \cup I\rangle = \langle G_1 \cup I\rangle \cup \langle G_2 \cup I\rangle, *, o)$ be a neutrosophic bigroup having a S-neutrosophic subbigroup, then $\langle G \cup I\rangle$ is a S-neutrosophic bigroup.*

The other result can also be stated as a theorem.

**THEOREM 2.2.2:** *Let $\langle G \cup I\rangle = (\langle G_1 \cup I\rangle \cup \langle G_2 \cup I\rangle, *, o)$ be any S-neutrosophic bigroup. Every neutrosophic subbigroup of $\langle G \cup I\rangle$ need not in general be a S-neutrosophic subbigroup.*

As in case of neutrosophic bigroup we can even in the case of S-neutrosophic bigroup define the order of the S-neutrosophic bigroup, the order is finite if the number of elements in $\langle G \cup I\rangle$ is finite; otherwise infinite.

We denote the order of the S-neutrosophic bigroup by $\langle G \cup I\rangle$ $o(\langle G \cup I\rangle)$.

Now we see in general if a S-neutrosophic bigroup is finite still the order of the S-neutrosophic subbigroup need not in general divide the order of the S-neutrosophic bigroup. So to characterize such S-neutrosophic subbigroups which divides the order of the S-neutrosophic bigroup we make the following definition.

**DEFINITION 2.2.7:** *Let $(\langle G \cup I\rangle = \langle G_1 \cup I\rangle \cup \langle G_2 \cup I\rangle, *, o)$ be a S-neutrosophic bigroup of finite order. Let $P = (\langle P_1 \cup I\rangle \cup \langle P_2 \cup I\rangle, *, o)$ be a S-neutrosophic subbigroup of $\langle G \cup I\rangle$. If $o(P) / o(\langle G \cup I\rangle)$ then we call P a Lagrange S-neutrosophic subbigroup II, if every S-neutrosophic subbigroup of $\langle G \cup I\rangle$ is Lagrange, then we call the S-neutrosophic bigroup $\langle G \cup I\rangle$ to be a Lagrange Smarandache neutrosophic bigroup II (Lagrange S-neutrosophic bigroup II).*

*If $\langle G \cup I\rangle$ has atleast one Lagrange S-neutrosophic subbigroup P then we call $\langle G \cup I\rangle$ to be a weakly Lagrange Smarandache neutrosophic bigroup II (weakly Lagrange S-neutrosophic bigroup II).*

*If $\langle G \cup I\rangle$ has no Lagrange S-neutrosophic subbigroup then we call $\langle G \cup I\rangle$ to be a Lagrange free S-neutrosophic bigroup II.*

*Now if $\langle G \cup I\rangle$ has a neutrosophic subbigroup T which is not a S-neutrosophic subbigroup II and if $o(T) / o(\langle G \cup I\rangle)$ then we call T to be a pseudo Lagrange neutrosophic subbigroup II. If every neutrosophic subbigroup T which is not a S-neutrosophic sub bigroup II is a pseudo Lagrange neutrosophic subbigroup then we call $\langle G \cup I\rangle$ to be pseudo Lagrange S-neutrosophic bigroup II.*

*If $\langle G \cup I\rangle$ is both a Lagrange S-neutrosophic bigroup II and a pseudo Lagrange S-neutrosophic bigroup then we call $\langle G \cup I\rangle$*

*a strong Lagrange S-neutrosophic bigroup. If the S-neutrosophic bigroup has atleast one pseudo Lagrange neutrosophic subbigroup then we call ⟨G ∪ I⟩ to be a weakly pseudo Lagrange S-neutrosophic bigroup II. If ⟨G ∪ I⟩ has no pseudo Lagrange neutrosophic subbigroup II then we call ⟨G ∪ I⟩ a pseudo Lagrange free S-neutrosophic bigroup II.*

Several interesting results in this direction can be obtained by an innovative reader. Now we proceed in to get some analogue for Sylow S-neutrosophic bigroups.

**DEFINITION 2.2.8:** *Let (⟨G ∪ I⟩ = ⟨$G_1$ ∪ I⟩ ∪ ⟨$G_2$ ∪ I⟩, *, o) be a S-neutrosophic bigroup of finite order. Suppose p is a prime such that $p^{\alpha}$ / o(⟨G ∪ I⟩) and $p^{\alpha+1}$ ∤ o(⟨G ∪ I⟩) and if ⟨G ∪ I⟩ has a S-neutrosophic subbigroup V of order $p^{\alpha}$, then we call V a p-Sylow S-neutrosophic subbigroup II. If for every prime p such that $p^{\alpha}$ / o(⟨G ∪ I⟩) and $p^{\alpha+1}$ ∤ o(⟨G ∪ I⟩) we have a p-Sylow S-neutrosophic subbigroup II then we call ⟨G ∪ I⟩ to be a Sylow S-neutrosophic bigroup II. If ⟨G ∪ I⟩ has atleast one p-Sylow S-neutrosophic subbigroup II then we call ⟨G ∪ I⟩ a weakly Sylow S-neutrosophic bigroup II. If ⟨G ∪ I⟩ has no p-Sylow S-neutrosophic subbigroup II then we call ⟨G ∪ I⟩ a Sylow free S-neutrosophic bigroup. Now for the finite S-neutrosophic bigroup ⟨G ∪ I⟩ if for p a prime with $p^{\alpha}$ / o(⟨G ∪ I⟩) and $p^{\alpha+1}$ ∤ o (⟨G ∪ I⟩) we have neutrosophic subbigroup W of order $p^{\alpha}$ which is not a S-neutrosophic subbigroup then we call W the pseudo p-Sylow neutrosophic subbigroup II of ⟨G ∪ I⟩.*

*If for every prime p with $p^{\alpha}$ / o(⟨G ∪ I⟩) and $p^{\alpha+1}$ ∤ o(⟨G ∪ I⟩) we have a pseudo p-Sylow neutrosophic subbigroup II then we call ⟨G ∪ I⟩ a pseudo p-Sylow S-neutrosophic bigroup II. If ⟨G ∪ I⟩ is both a pseudo Sylow S-neutrosophic bigroup II and Sylow S-neutrosophic bigroup II then we call ⟨G ∪ I⟩ a strong Sylow S-neutrosophic bigroup II.*

*If ⟨G ∪ I⟩ has atleast one pseudo p-Sylow neutrosophic subbigroup then we call ⟨G ∪ I⟩ a weakly pseudo Sylow S-neutrosophic bigroup. If ⟨G ∪ I⟩ is only a Sylow S-neutrosophic bigroup II or only a pseudo Sylow S-neutrosophic bigroup II then we call ⟨G ∪ I⟩ a semi Sylow S-neutrosophic bigroup II.*

Here also we leave it as an exercise for the reader to find the interrelations and illustrative examples. However we give a few examples so that it will help the reader to easily understand the problems.

***Example 2.2.7:*** Let (⟨G ∪ I⟩ = ⟨$G_1$ ∪ I⟩ ∪ ⟨$G_2$ ∪ I⟩, *, o), be a S-neutrosophic bigroup where

| | | |
|---|---|---|
| ⟨$G_1$ ∪ I⟩ | = | {1, 2, 3, 4, I, 2I, 3I, 4I}, a neutrosophic group under multiplication modulo 5 and |
| ⟨$G_2$ ∪ I⟩ | = | {0, 1, 2, 3, I, 2I, 3I}, is a S-neutrosophic semigroup under multiplication modulo 4. |

o(⟨G ∪ I⟩) = 15, 3/15 $3^2$ ╳ 15, 5/15 and $5^2$ ╳ 15.

Take H = {1, I} ∪ {0, 2, 2I} is a S-neutrosophic subbigroup o(H) / o(⟨G ∪ I⟩). Consider T = {1, I, 4, 4I} ∪ {1, I, 2, 2I, 0}, o(T) = 9 and 9 ╳ 15.

But ⟨G ∪ I⟩ has no S-neutrosophic subbigroup of order 3. Take V = {1, 2, 3, 4} ∪ {0, 1, 2, 3}, V is a subbigroup of ⟨G ∪ I⟩, o(V) ╳ o (⟨G ∪ I⟩).

Suppose B = {4, 4I} ∪ {0, 2, 2I}, B is a S-neutrosophic subbigroup of order 5. W = {1, 4, I, 4I} ∪ {0, 2}. W is a neutrosophic subbigroup and o(W) = 6 and 6 ╳ 15. Thus we can get many subgroups and interesting results. We however define another new notion.

**DEFINITION 2.2.9:** *Let (⟨G ∪ I⟩ = (⟨$G_1$ ∪ I⟩ ∪ ⟨$G_2$ ∪ I⟩, *, o) be a S-neutrosophic bigroup of finite order. Suppose ⟨G ∪ I⟩ is a Sylow S-neutrosophic bigroup II and if in addition we have for every prime p, $p^{\alpha}$ / o(⟨G ∪ I⟩) and $p^{\alpha+1}$ ╳ o (⟨G ∪ I⟩) their exists a S-neutrosophic subbigroup of order $p^{\alpha+t}$ ( $t \geq 1$) then we call ⟨G ∪ I⟩ a super Sylow S-neutrosophic bigroup II.*

It is interesting to see that all super Sylow S-neutrosophic bigroups are Sylow S-neutrosophic bigroups but a Sylow S-neutrosophic bigroup need not in general be a super Sylow S-neutrosophic bigroup.

Now we proceed on to define the notion of Cauchy S-neutrosophic group weakly Cauchy S-neutrosophic S-bigroup and so on.

**DEFINITION 2.2.10:** *Let* $(\langle G \cup I\rangle = (\langle G_1 \cup I\rangle \cup \langle G_2 \cup I\rangle, *, o)$ *be a finite S-neutrosophic bigroup. Let* $x \in \langle G \cup I\rangle$ *be such that* $x^m = 1$ *and if* $m / o\,(\langle G \cup I\,\rangle)$ *we call x a Cauchy element. If for* $y \in \langle G \cup I\rangle$, $y^t = I$ *and* $t / o(\langle G \cup I\rangle)$ *then we call y a Cauchy neutrosophic element of* $\langle G \cup I\rangle$*. If every element x such that* $x^m = 1$ *and y such that* $y^t = I$*, are Cauchy element or Cauchy neutrosophic element then we call* $\langle G \cup I\rangle$ *to be a Cauchy S-neutrosophic bigroup. If* $\langle G \cup I\rangle$ *has at least one Cauchy element and one Cauchy neutrosophic element then we call* $\langle G \cup I\rangle$ *a weakly Cauchy S-neutrosophic bigroup II. If* $\langle G \cup I\rangle$ *has no Cauchy element or Cauchy neutrosophic element then we call* $\langle G \cup I\rangle$ *a Cauchy free S-neutrosophic bigroup II.*

The following result is interesting.

**THEOREM 2.2.3:** *All S-neutrosophic bigroups of order p, p a prime are Cauchy free S-neutrosophic bigroups.*

*Proof:* Given $o(\langle G \cup I\rangle) = p$, p a prime so every element which are such that $x^m = 1$ or $y^t = I$ will have $m < p$ and $t < p$ so $(m, p) = 1$ and $(t, p) = 1$. Thus $\langle G \cup I\rangle$ is a Cauchy free S-neutrosophic bigroup.

Now we proceed on to give some of its properties about their structure.

**DEFINITION 2.2.11:** *Let* $\{\langle G \cup I\rangle = \langle G_1 \cup I\rangle \cup \langle G_2 \cup I\rangle, *, o\}$ *be a S-neutrosophic bigroup.* $\langle G \cup I\rangle$ *is said to be a Smarandache commutative bigroup (S-commutative bigroup) if* $\langle G_1 \cup I\rangle$ *which is the neutrosophic group is commutative and every proper subset of the S-neutrosophic semigroup* $\langle G_2 \cup I\rangle$ *which is a neutrosophic group is a commutative group. If both* $\langle G_1 \cup I\rangle$ *and* $\langle G_2 \cup I\rangle$ *happen to be commutative* $\langle G \cup I\rangle$ *is trivially a S-commutative bigroup.*

Now we proceed on to define the notion of Smarandache weakly commutative neutrosophic bigroup.

**DEFINITION 2.2.12:** *Let* $(\langle G \cup I\rangle = \langle G_1 \cup I\rangle \cup \langle G_2 \cup I\rangle, o, *)$ *be S-neutrosophic bigroup. If every neutrosophic subgroup of* $\langle G_1 \cup I\rangle$ *and atleast a proper subset which is a subgroup of the S-neutrosophic semigroup happen to be commutative we call* $\langle G \cup I\rangle$ *to be a Smarandache weakly neutrosophic bigroup (S-weakly neutrosophic bigroup).*

The following theorem is left as an exercise for the reader to prove.

**THEOREM 2.2.4:** *Let* $(\langle G \cup I\rangle = \langle G_1 \cup I\rangle \cup \langle G_2 \cup I\rangle, o, *)$*, be a S-commutative neutrosophic bigroup.* $\langle G \cup I\rangle$ *is a S-weakly commutative neutrosophic bigroup and not conversely.*

We now proceed on to define Smarandache hyper neutrosophic bigroups for this we basically need the notion of largest group in a bigroup.

**DEFINITION 2.2.13:** *Let* $(\langle G \cup I\rangle = \langle G_1 \cup I\rangle \cup \langle G_2 \cup I\rangle, o, *)$ *be a S-neutrosophic bigroup* $P = \langle P_1 \cup I\rangle \cup \langle P_2 \cup I\rangle, o, *)$ *is said to be a Smarandache largest neutrosophic bigroup (S- largest neutrosophic bigroup) if* $\langle P_1 \cup I\rangle$ *is the largest neutrosophic subgroup of* $\langle G_1 \cup I\rangle$ *(If* $\langle G_1 \cup I\rangle$ *has no proper neutrosophic subgroup take* $\langle G1 \cup I\rangle$ *itself to be trivially largest subgroup) and* $\langle P_2 \cup I\rangle$ *is largest proper subset of* $\langle G_2 \cup I\rangle$ *which is the largest neutrosophic group of the S-neutrosophic semigroup* $\langle G_2 \cup I\rangle$*. We call this largest Smarandache neutrosophic subbigroup to be the Smarandache hyper neutrosophic bigroup (S-hyper neutrosophic bigroup).*

**DEFINITION 2.2.14:** *If* $(\langle G \cup I\rangle = \langle G_1 \cup I\rangle \cup \langle G_2 \cup I\rangle, *, o)$ *is a S-neutrosophic bigroup having no S-hyper neutrosophic bigroup then we call* $\langle G \cup I\rangle$ *to be a Smarandache simple neutrosophic bigroup (S-simple neutrosophic bigroup).*

Many interesting properties about S-neutrosophic bigroups can be derived as in case of other S-bigroups, bigroups and neutrosophic bigroups. This can be taken up as simple exercises by any interesting reader.

## 2.3 Smarandache Neutrosophic N-groups

In this section for the time the definition of Smarandache neutrosophic N-groups are defined. They are more generalized structures for they involve S-neutrosophic semigroups and S-semigroups. These structures enjoy more properties than other structures for this class inducts the generalized structures in the very definition. Several properties are derived. Illustrative examples are provided for the reader.

Now we proceed on to define the notion of Smarandache neutrosophic N-groups. For more about neutrosophic N-groups, refer [51].

**DEFINITION 2.3.1:** *Let* $(\langle G \cup I\rangle = \langle G_1 \cup I\rangle \cup \ldots \cup \langle G_N \cup I\rangle, *_1, *_2, \ldots, *_N)$ *be a nonempty set with N-binary operations. We say* $\langle G \cup I\rangle$ *is a Smarandache neutrosophic N-group (S-neutrosophic N-group) if*

i. *Some of* $(G_i, *_i)$ *is a S-semigroup.*
ii. *Some of* $(G_j, *_j)$ *are neutrosophic groups or groups.*
iii. *Rest of* $(G_k, *_k)$ *are S-neutrosophic semigroups,* $1 \leq i, j, k \leq N$.

***Example 2.3.1:*** Let $(\langle G \cup I\rangle = \langle G_1 \cup I\rangle \cup \langle G_2 \cup I\rangle \cup \langle G_3 \cup I\rangle, *_1, *_2, *_3)$ where $\langle G_1 \cup I\rangle = \{1, 2, 3, 4, I, 2I, 3I, 4I\}$, a neutrosophic group under multiplication modulo 5. $\langle G_2 \cup I\rangle = \{0, 1, 2, 3, 4, 5, I, 2I, 3I, 4I, 5I\}$, S-neutrosophic semigroup under multiplication modulo 6 and $G_3 = A_4$. Clearly $\langle G \cup I\rangle$ is a S-neutrosophic 3-group.

Now we define order and substructures in them.

**DEFINITION 2.3.2:** *Let* $(\langle G \cup I\rangle = \langle G_1 \cup I\rangle \cup \langle G_2 \cup I\rangle \cup \ldots \cup \langle G_N \cup I\rangle, *_1, *_2, \ldots, *_N)$ *be a neutrosophic N-group. A proper subset P of* $\langle G \cup I\rangle$ *is said to be a Smarandache neutrosophic sub N-group (S-neutrosophic sub N-group) if* $P = (P_1 \cup P_2 \cup \ldots \cup P_N, *_1, *_2, \ldots, *_N)$ *is a S-neutrosophic N-group, where* $P_i \subset G_i$ *and* $P_i = P \cap G_i$, $i = 1, 2, \ldots, N$.

It is interesting to note that we have taken in the definition of the S-neutrosophic sub N-group the neutrosophic N-group and not a S-neutrosophic N-group because of the following theorem.

**THEOREM 2.3.1:** *Let* $(\langle G \cup I\rangle = \langle G_1 \cup I\rangle \cup \langle G_2 \cup I\rangle \cup \ldots \cup \langle G_N \cup I\rangle, *_1, *_2, \ldots, *_N)$ *be a neutrosophic N-group. If* $\langle G \cup I\rangle$ *has a S-neutrosophic sub N-group then* $\langle G \cup I\rangle$ *is itself a S-neutrosophic N-group.*

*Proof:* Given $\langle G \cup I\rangle$ is a neutrosophic N-group having a S-neutrosophic sub N-group. Let $P = \{P_1 \cup P_2 \cup \ldots \cup P_N, *_1, \ldots, *_N\}$ be a S-neutrosophic sub N-group of $\langle G \cup I\rangle$ so

(1) $P_i$'s are subgroups or neutrosophic subgroups.
(2) $P_j$'s are S-semigroups.
(3) $P_k$'s are S-neutrosophic subsemigroups.

$1 \le i, j, k \le N$. As $P \subset \langle G \cup I\rangle$ and each $P_i \subset G_i$ we see by the very definition $\langle G \cup I\rangle$ is a S-neutrosophic N-group.

The number of distinct elements in a S-neutrosophic N-group $\langle G \cup I\rangle$ gives the order of $(\langle G \cup I\rangle)$ denoted by $o(\langle G \cup I\rangle)$. If $\langle G \cup I\rangle$ has finite number of elements we call the S-neutrosophic N-group to be finite otherwise infinite.

***Example 2.3.2:*** Let $\{\langle G \cup I\rangle = \{\langle Z \cup I\rangle\} \cup S_3 \cup \langle 0, 1, 2, 3, I, 2I, 3I\rangle, *_1, *_2, *_3\}$ be a S-neutrosophic 3 group where $\langle Z \cup I\rangle$ is a neutrosophic group under addition. $S_3$ is the symmetric group of degree 3 and $\langle 0, 1, 2, \ldots, 3I\rangle$ is a neutrosophic semigroup under multiplication modulo 4. Clearly $\langle G \cup I\rangle$ is an infinite S-neutrosophic N-group.

Now we see even if one of the $(G_i, *_i)$ is infinite then $\langle G \cup I\rangle$ is infinite. Only if all the $(G_i, *_i)$, are finite we get $\langle G \cup I\rangle$ to be finite.

It is an interesting thing to note that in general the order of the S-neutrosophic sub N-group need not divide the order of the S-neutrosophic group. To get some conditions we define the following.

**DEFINITION 2.3.3:** *Let $(\langle G \cup I\rangle = \langle G_1 \cup I\rangle \cup \langle G_2 \cup I\rangle \cup \ldots \cup \langle G_N \cup I\rangle, *_1, *_2, \ldots, *_N)$ be a finite S-neutrosophic N-group. A S-neutrosophic sub N-group $P = \{P_1 \cup P_2 \cup \ldots \cup P_N, *_1, \ldots, *_N\}$ is said to be a Lagrange S-neutrosophic sub N-group if $o(P) / o(\langle G \cup I\rangle)$.*

*If every S-neutrosophic sub N-group of $\langle G \cup I\rangle$ is a Lagrange S-neutrosophic sub N-group then we call $\langle G \cup I\rangle$ to be a Lagrange S-neutrosophic N-group. If $\langle G \cup I\rangle$ has atleast one Lagrange S-neutrosophic sub N-group then we call $\langle G \cup I\rangle$ to be a weakly Lagrange S-neutrosophic N-group. If $\langle G \cup I\rangle$ has no Lagrange S-neutrosophic sub N-group then we call $\langle G \cup I\rangle$ a Lagrange free S-neutrosophic N-group.*

Now we will illustrate this by an example before we proceed on to define other types of Lagrange neutrosophic sub N-groups.

***Example 2.3.3:*** Let $(\langle G \cup I\rangle = \langle G_1 \cup I\rangle \cup G_2 \cup G_3, *_1, *_2, *_3)$ be a S-neutrosophic 3-group where

$\langle G_1 \cup I\rangle$ = {0, 1, 2, 3, 4, 5, I, 2I, 3I, 4I, 5I}, a S-neutrosophic semigroup under multiplication modulo 6.

$G_2$ = $S_3$ and

$G_3$ = {0, 1, 2, …, 12, 13, 14}, S-semigroup under multiplication modulo 15.

$o(\langle G \cup I\rangle) = 32$.

$P = \{P_1 \cup P_2 \cup P_3, *_1, *_2, *_3\}$ where $P_1$ = {0, 2, 4, 2I, 4I}, a S-neutrosophic subsemigroup of $\langle G_1 \cup I\rangle$,

$$P_2 = \left\{ \begin{pmatrix} 1 & 2 & 3 \\ 1 & 2 & 3 \end{pmatrix}, \begin{pmatrix} 1 & 2 & 3 \\ 2 & 1 & 3 \end{pmatrix} \right\} \subset S_3$$

and $P_3 = \{0, 3, 6, 9, 12\} \subset G_3$. P is a S-neutrosophic sub 3-group of $\langle G \cup I \rangle$. o(P) = 12, $12 \nmid 32$. So P is not a Lagrange S-neutrosophic sub 3-group of $\langle G \cup I \rangle$. $K = (K_1 \cup K_2 \cup K_3, *_1, *_2, *_3)$ where, $K_1 = \{1, 5, I, 5I\}$ a S-neutrosophic subsemigroup;

$$P_2 = \left\{ \begin{pmatrix} 1 & 2 & 3 \\ 1 & 2 & 3 \end{pmatrix}, \begin{pmatrix} 1 & 2 & 3 \\ 1 & 3 & 2 \end{pmatrix} \right\} \subset S_3 \text{ and}$$

$P_3 = \{1, 14\} \subset G_3$. Clearly o(P) = 8 and 8/32, so K is a Lagrange S-neutrosophic sub N-group of $\langle G \cup I \rangle$, hence $\langle G \cup I \rangle$ is a weakly Lagrange S-neutrosophic N-group. Thus we see when we say $\langle G \cup I \rangle$ has no Lagrange S-neutrosophic sub N-group we mean that even we can have many S-neutrosophic sub N-groups but the order of them will not divide the order of $\langle G \cup I \rangle$.

Now we proceed on to define two types of Cauchy elements.

**DEFINITION 2.3.4:** *Let $(\langle G \cup I \rangle = \langle G_1 \cup I \rangle \cup \langle G_2 \cup I \rangle \cup \ldots \cup \langle G_N \cup I \rangle, *_1, *_2, \ldots, *_N)$ be a finite S-neutrosophic N-group. An element x is Cauchy neutrosophic if $x^n = I$ and $n / o(\langle G \cup I \rangle)$. y is Cauchy if $y^m = 1$ and $m / o(\langle G \cup I \rangle)$. x is S-Cauchy neutrosophic if $x^n = I$ and $n / o\langle H \cup I \rangle$ where $x \in \langle H \cup I \rangle$ and $\langle H \cup I \rangle$ is a neutrosophic sub N-group of $\langle G \cup I \rangle$. $y \in \langle H \cup I \rangle$ is S-Cauchy if $y^m = 1$, $m / o\langle H \cup I \rangle$.*

It is important and interesting to note that a S-Cauchy neutrosophic element may not in general be a Cauchy neutrosophic element of $\langle G \cup I \rangle$. Similarly a S-Cauchy element also may not in general be a Cauchy element of $\langle G \cup I \rangle$.

The interested reader can find conditions under which a S-Cauchy neutrosophic element is a Cauchy neutrosophic element of $\langle G \cup I \rangle$.

Now we illustrate these by the following example.

***Example 2.3.4:*** Let $(\langle G \cup I\rangle = \langle G_1 \cup I\rangle \cup \langle G_2 \cup I\rangle \cup \langle G_3 \cup I\rangle$, $*_1$, $*_2$, $*_3)$ where $(\langle G_1 \cup I\rangle) = \{1, 2, 3, 4, I, 2I, 3I, 4I\}$, S-neutrosophic group under multiplication modulo 5. $G_2 = Z_{10}$, S-semigroup under multiplication modulo 10, $G_3 = \{0, 1, 2, 3, I, 2I, 3I\}$, S-neutrosophic group under multiplication modulo 4. $o(\langle G \cup I\rangle) = 25$.

Let $4I \in \langle G \cup I\rangle$ $(4I)^2 = I$ but $2 \nmid o\langle G \cup I\rangle$ so 4I is not a Cauchy neutrosophic element of $G \cup I$. Take $3 \in \langle G \cup I\rangle$, $3^4 = 1 \pmod{10}$ so 3 is not a Cauchy element.

Take $P = \{P_1 \cup P_2 \cup P_3, *_1, *_2, *_3\}$, a S-neutrosophic sub 3-group of $\langle G \cup I\rangle$. $P_1 = \{1, 4, I, 4I\} \subset (\langle G_1 \cup I\rangle)$, $P_2 = \{1, 9\} \subset G_2 = Z_{10}$, $G_3 = \{1, 3\} \subset G_3$, $o(P) = 8$. Every element is either a S-Cauchy element or a S-Cauchy neutrosophic element of $\langle G \cup I\rangle$; but none of them is a Cauchy element or a Cauchy neutrosophic element of $\langle G \cup I\rangle$.

Now we proceed on to define two types S-Sylow structures on Smarandache neutrosophic N-groups.

**DEFINITION 2.3.5:** *Let $(\langle G \cup I\rangle = \langle G_1 \cup I\rangle \cup \ldots \cup \langle G_N \cup I\rangle$, $*_1, \ldots, *_N)$ be a S-neutrosophic N-group of finite order. Let p be a prime such that $p^\alpha / o(\langle G \cup I\rangle)$ and $p^{\alpha+1} \nmid o\langle G \cup I\rangle$. If $\langle G \cup I\rangle$ has a S-neutrosophic sub N-group P of order $p^\alpha$ then we call P a p-Sylow S-neutrosophic sub N-group.*

*If for every prime p such that $p^\alpha / o(\langle G \cup I\rangle)$ and $p^{\alpha+1} \nmid o(\langle G \cup I\rangle)$ we have a p-Sylow S-neutrosophic sub-N-group then we call $\langle G \cup I\rangle$ a Sylow S-neutrosophic N-group.*

*Note:* We will have for S-neutrosophic N-groups of finite order p-Sylow neutrosophic N-group and also the notion of Sylow neutrosophic N-group. Our main concern and interest is can a S-neutrosophic N-group be both Sylow-S-neutrosophic N-group and Sylow neutrosophic N-group. We just give some examples.

***Example 2.3.5:*** Let $(\langle G \cup I\rangle = \langle G_1 \cup I\rangle \cup \langle G_2 \cup I\rangle \cup G_3$, $*_1$, $*_2$, $*_3)$ where

| | | |
|---|---|---|
| $\langle G_1 \cup I\rangle$ | = | {1, 2, 3, 4, I, 2I, 3I, 4I}, S-neutrosophic group. |
| $\langle G_2 \cup I\rangle$ | = | {0, 1, 2, 3, I, 2I, 3I 1 + I, 2 + I, 3 + I, 1 + 2I, 2 + 2I, 3 + 2I, 3 + I, 3 + 2I, 3 + 3I}, and |
| $G_3$ | = | {$Z_{12}$, S-semigroup under multiplication modulo 12}, |

o($\langle G \cup I\rangle$) = 36. Now 2/36, $2^2$ / 36 and $2^3 \not| $ 36. Also 3/36, $3^2$/36 and $3^3 \not|$ 36. Take P = {1, I} ∪ {1, 3, 0} ∪ {0, 3, 6, 9}. P is a S-neutrosophic sub 3-group and infact P is a 3-Sylow S-neutrosophic sub 3-group. Clearly $\langle G \cup I\rangle$ has no 2-Sylow S-neutrosophic sub 3-group.

Take T = {1 + 3I, 1} ∪ {0} ∪ {4}; T is a neutrosophic sub 3-group. o(T) / o($\langle G \cup I\rangle$). So T is a 2-Sylow neutrosophic sub 3-group. However $\langle G \cup I\rangle$ has no 2 Sylow S-neutrosophic sub 3-group.

Take V = {$V_1 \cup V_2 \cup V_3$, $*_1$, $*_2$, $*_3$} where $V_1$ = {1, I, 4, 4I}, $V_2$ = {0, 2, 1} and $V_3$ = {1, 11}. V is a 3-Sylow neutrosophic sub 3-group as o(V) = 9. This $\langle G \cup I\rangle$ is only a Sylow neutrosophic 3-group and not a Sylow S-neutrosophic 3-group.

However it is important to note $\langle G \cup I\rangle$ has several neutrosophic sub N-groups which are not Lagrange. It is to be noted that a S-neutrosophic N-group need not in general be a neutrosophic N-group. So while defining ($p_1$, …, $p_N$)-Sylow subgroup of ($\langle G \cup I\rangle = \langle G_1 \cup I\rangle \cup \ldots \cup \langle G_N \cup I\rangle$, $*_1$, …, $*_N$) we have to be very careful. Infact we can more so define only ($p_1$, …, $p_N$)-Sylow group for the proper subset of $\langle G \cup I\rangle$ which is a neutrosophic N-group.

So we proceed on to define keeping in mind these observations.

**DEFINITION 2.3.6:** *Let ($\langle G \cup I\rangle = \langle G_1 \cup I\rangle \cup \ldots \cup \langle G_N \cup I\rangle$, $*_1, \ldots, *_N$) be a S-neutrosophic N-group of finite order. Suppose P = {$P_1 \cup P_2 \cup \ldots \cup P_N$, $*_1, \ldots, *_N$} is a proper subset of $\langle G \cup I\rangle$ which is a neutrosophic N-group. We define ($p_1$, …,$p_N$) S-Sylow neutrosophic sub N-group to be a $p_i$-Sylow subgroup of $P_i \subset G_i$ , i = 1, 2, …, N.*

Interested reader can develop results in this direction, however we illustrate this by the following example.

***Example 2.3.6:*** Let $\langle G \cup I\rangle = \{G_1 \cup G_2 \cup G_3, *_1, *_2, *_3\}$ be a S-neutrosophic group where $G_1 = S_4$, $\langle G_2 \cup I\rangle = \{1, 2, 3, 4, I, 2I, 3I, 4I\}$ and $G_3 = \{(a, b)| a, b \in \{0, 1, 2, I, 2I, 1 + I, 2 + 2I, 2+I, 2I + 2\}$ under component wise multiplication modulo 3.

Take $T = A_4 \cup \{1, 4, I, 4I\} \cup \{(1, 1), (2, 2), (1, 2), (2, 1)\}$, T has a (3, 2, 2) Sylow S-neutrosophic 3-group and (2, 2, 2) – Sylow S-neutrosophic group.

Several examples can be constructed.

Now we proceed on to define the concept of Smarandache neutrosophic homomorphism of N-groups.

**DEFINITION 2.3.7:** *Let $(\langle G \cup I\rangle = \langle G_1 \cup I\rangle \cup \langle G_2 \cup I\rangle \cup \ldots \cup \langle G_N \cup I\rangle, *_1, *_2, \ldots, *_N)$ be a S-neutrosophic N-group and $\langle H \cup I\rangle = (\langle H_1 \cup I\rangle \cup \ldots \cup \langle H_N \cup I\rangle, *_1, \ldots, *_N)$ be another S-neutrosophic N-group. We say a map $\phi$ from $\langle G \cup I\rangle$ to $\langle H \cup I\rangle$ is a S-neutrosophic homomorphism if $\phi(I) = I$ and $\phi$ is a neutrosophic N-group homomorphism from neutrosophic sub N-group $(T = T_1 \cup \ldots \cup T_N, *_1, \ldots, *_N)$ in $\langle G \cup I\rangle$ to $V = V_1 \cup V_2 \cup V_3 \cup \ldots \cup V_N, *_1, \ldots, *_N)$ in $(H \cup I)$.*

It is important to note that in general $\phi$ need not be even defined on the whole of $\langle G \cup I\rangle$.

Now $\phi \mid G_i : T_i \to V_i$ is either a group homomorphism or a neutrosophic group homomorphism for i = 1, 2, …, N. We denote $\phi / G_i$ need be defined on $T_i$ and need not be defined on the whole of $G_i$, which we denote by $\phi_i$, i = 1, 2, …, N i.e. $\phi = \phi_1 \cup \ldots \cup \phi_N$; $\phi_i : T_i \to V_i$, $1 \leq i \leq N$.

Next we define the notion of weak Smarandache neutrosophic N-groups.

**DEFINITION 2.3.8:** *A proper subset $(\langle G \cup I\rangle = \langle G_1 \cup I\rangle \cup \ldots \cup \langle G_N \cup I\rangle, *_1, \ldots, *_N)$ is said to be a weak Smarandache neutrosophic N-group (weak S-neutrosophic N-group) , if the* *following conditions are satisfied.* $\langle G \cup I\rangle = \langle G_1 \cup I\rangle \cup \ldots \cup G_i \cup \ldots \cup \langle G_N \cup I\rangle$ *in which* $\langle G_i \cup I\rangle$ *is a neutrosophic group or* $G_j$ *is S-semigroup. (Here also* $\langle G_i \cup I\rangle$ *and* $G_j$ *are proper subsets of* $\langle G \cup I\rangle$*).*

All S-neutrosophic N-groups are weak S-neutrosophic N-group and not conversely. We first illustrate this by the following example.

***Example 2.3.7:*** Let $\langle G \cup I\rangle = \{\langle G_1 \cup I\rangle \cup G_2 \cup G_3\}$ where $G_1 \cup I = \{1, 2, 3, 4 I, 2I, 3I, 4I\}$, neutrosophic group under multiplication modulo 5. $G_2 = \{Z_{12}$, S-semigroup under multiplication modulo 12} and $G_3$ = S (3) the symmetric group of all mappings of the set (1 2 3) to itself.

Clearly $Z_{12}$ and S (3) are S-semigroups as T = {4, 8} is group under multiplication modulo 12 and $S_3 \subset S(3)$, is the symmetric group of all one to one mappings of (1 2 3) to itself. Thus $\langle G \cup I\rangle$ is a weak S-neutrosophic group.

The subsets $\{\langle G_1 \cup I\rangle \cup T \cup S_3\}$ which are neutrosophic N-groups will be known as neutrosophic sub N-groups of the S-neutrosophic N-group $\langle G \cup I\rangle$.

As in case of all other algebraic structures we in case of weak S-neutrosophic N-groups $\langle G \cup I\rangle$, define the order of $\langle G \cup I\rangle$ to be the number of distinct elements in $\langle G \cup I\rangle$. If the number of elements in $\langle G \cup I\rangle$ is finite we call $\langle G \cup I\rangle$ a finite S-neutrosophic N-group otherwise we call them as infinite S-neutrosophic N-group and we denote the order of ($\langle G \cup I\rangle$ ) by o($\langle G \cup I\rangle$). Even if one of the neutrosophic groups ($\langle G_i \cup I\rangle$, $*_i$) or S-semigroups ($G_j$, $*_j$) is infinite then $\langle G \cup I\rangle$ will be infinite. If the order of every ($G_j$, $*_j$) is finite, then $\langle G \cup I\rangle$ is finite.

Now we proceed on to define the notion of substructures.

**DEFINITION 2.3.9:** *Let* $\langle G \cup I\rangle = (\langle G_1 \cup I\rangle \cup \ldots \cup \langle G_N \cup I\rangle, *_1, \ldots, *_N)$ *be a weak S-neutrosophic N-group. A proper subset P of* $\langle G \cup I\rangle$ *is said to be a weak Smarandache neutrosophic sub N-group (weak S-neutrosophic sub N-group) if P satisfies the following conditions*

i. *$P = P_1 \cup \ldots \cup P_N$ is such that each $P_i = P \cap G_i$ and $P_i$ are neutrosophic subgroups or S-subsemigroups, $1 \leq i \leq N$*

ii. *($P = P_1 \cup \ldots \cup P_N$, $*_1$, …, $*_N$) is itself a weak S-neutrosophic N-group.*

*(Weak S-neutrosophic N-groups in general are not neutrosophic N-groups).*

Now as in case of all other algebraic structures by the order of the weak S-neutrosophic N-group we mean the number of distinct elements in them. It is finite if the number of distinct elements in $\langle G \cup I \rangle$ is finite i.e. all $\langle G_i \cup I \rangle$ are finite. The weak S-neutrosophic N-group is infinite if its order is infinite i.e. it has atleast one of $(G_i, *_i)$ to be infinite.

Now we just give an example, which will lead us to define several types of algebraic properties.

***Example 2.3.8:*** Let $\{\langle G \cup I \rangle = \langle G_1 \cup I \rangle \cup G_2 \cup G_3 \cup G_4\}$ be a weak S-neutrosophic 4-group, where

| | | |
|---|---|---|
| $\langle G_1 \cup I \rangle$ | = | {1, 2, 3, 4, I, 2I, 3I, 4I} neutrosophic group under multiplication modulo 5. |
| $G_2$ | = | S(3), |
| $G_3$ | = | {$Z_{12}$, semigroup under multiplication modulo 12} and |
| $G_4$ | = | $\{(a, b) \mid a \in Z_2 \text{ and } b \in Z_4\}$; |

$o(\langle G \cup I \rangle) = 55$.

Take $P = \{P_1 \cup P_2 \cup P_3 \cup P_4\}$ where $P_1 = \{1, 4, I, 4I\}$,

$$P_2 = \left\{ \begin{pmatrix} 1 & 2 & 3 \\ 1 & 2 & 3 \end{pmatrix}, \begin{pmatrix} 1 & 2 & 3 \\ 2 & 3 & 1 \end{pmatrix}, \begin{pmatrix} 1 & 2 & 3 \\ 3 & 1 & 2 \end{pmatrix} \right\},$$

$P_3 = \{1, g^3, g^6, g^4\}$ and $P_4 = \{(1, 1), (1, 3)\}$.
$o(P) = 13$, $(13, 55) = 1$.

Now take $T = T_1 \cup T_2 \cup T_3 \cup T_4$ where $T_1 = P_1$, $T_2 = P_2$, $T_3 = \{1, g^6\}$ and $T_4 = P_4$. $o(T) = 11$ and $o(T) / o(\langle G \cup I\rangle)$. Thus we see the order of weak S-neutrosophic sub N-group need not in general divide the order of $\langle G \cup I\rangle$.

To this end we have many definitions.

**DEFINITION 2.3.10:** *Let ($\langle G \cup I\rangle = \langle G_1 \cup I\rangle \cup \langle G_2 \cup I\rangle \cup \ldots \cup \langle G_N \cup I\rangle$, $*_1$, $*_2$, …, $*_N$) be a weakly S-neutrosophic N-group of finite order. Let $P = \{P_1 \cup P_2 \cup P_3, \cup \ldots \cup *_1, \ldots, *_N\}$ be a weakly S-neutrosophic sub N-group of $\langle G \cup I\rangle$. We say P is a Lagrange weakly S-neutrosophic sub N-group if $o(P) / o(\langle G \cup I\rangle)$ otherwise P is called non Lagrange weakly S-neutrosophic sub N-group.*

*If every weakly S-neutrosophic sub N-group is Lagrange then we call $\langle G \cup I\rangle$ a Lagrange weakly S-neutrosophic N-group. If $\langle G \cup I\rangle$ has atleast one Lagrange weakly S-neutrosophic sub N-group we call $\langle G \cup I\rangle$ a weakly Lagrange weakly S-neutrosophic N-group.*

Clearly all Lagrange weakly S-neutrosophic N-groups are weakly Lagrange weakly S-neutrosophic N-group, however the converse is not true. If $\langle G \cup I\rangle$ has no Lagrange S-neutrosophic sub N-group then we call $\langle G \cup I\rangle$, a Lagrange free weakly S-neutrosophic N-group.

We can study these structures. It is interesting to note all weakly S-neutrosophic N-groups of order p, p a prime are Lagrange free weakly S-neutrosophic N-groups. However we have weakly S- neutrosophic N-groups of finite order say a composite number can also be a Lagrange free weakly S-neutrosophic N-groups.

Now we can define Cauchy element and Cauchy neutrosophic element in an analogous way as in case of S-neutrosophic N-group. All the notions of p-Sylow weakly S-neutrosophic sub N-groups, Sylow weakly S-neutrosophic N-groups Sylow free weakly S-neutrosophic N-groups can also be defined in an analogous way. Now we define semi Smarandache neutrosophic N-groups.

**DEFINITION 2.3.11:** *A non empty set* $\{\langle G \cup I\rangle, *_1, *_2,\ldots, *_N\}$ *with N-binary operations is said to be a semi Smarandache neutrosophic N-group (semi S-neutrosophic N-group) if the following conditions are satisfied.*

i. $\langle G \cup I\rangle = G_1 \cup G_2 \cup \langle G_3 \cup I\rangle \cup \ldots \cup G_N$ *is such that each* $G_i$ *or* $\langle G_j \cup I\rangle$ *is a proper subset of* $\langle G \cup I\rangle$.
ii. *Some of* $(G_i, *_i)$ *are groups*
iii. *Rest of* $(\langle G_j \cup I\rangle, *_j)$ *are S-neutrosophic semigroups* $(1 \leq i, j \leq N)$.

***Example 2.3.9:*** Consider $\langle G \cup I\rangle = \{\langle G_1 \cup G_2 \cup G_3 \cup \langle G_4 \cup I\rangle, *_1, \ldots, *_4\}$ where $G_1 = S_3$, $G_2 = \langle g \mid g^{12} = e\rangle$, $G_3 = D_{2.7}$ and $\langle G_4 \cup I\rangle = \{0, 1, 2, 3, 4, 5, I, 2I, 3I, 4I, 5I\}$ is a S-neutrosophic semigroup. Clearly $\langle G \cup I\rangle$ is a semi S-neutrosophic 4-group.

As in case of other algebraic structures we can define order and semi S-neutrosophic sub N-group, the Lagrange analogue, the Sylow analogue, the Cauchy analogue and homomorphism and so on.

## Chapter three

# SMARANDACHE NEUTROSOPHIC SEMIGROUPS AND THEIR GENERALIZATIONS

This chapter has three sections. Section one introduces for the first time the notion of Smarandache neutrosophic semigroups, bisemigroups and N-semigroups. We also give examples of them. In section 1 S-neutrosophic semigroups are introduced. Their substructures like S-neutrosophic subsemigroups, S-neutrosophic ideals, S-neutrosophic Lagrange subsemigroups, p-Sylow S-neutrosophic subsemigroups etc. are introduced and examples are provided. Section two deals with Smarandache neutrosophic bisemigroups. Special properties about S-neutrosophic bisemigroups are studied. The notion of Smarandache neutrosophic N-semigroups are introduced in section three. Definition of S-Lagrange neutrosophic sub N semigroup, S-maximal neutrosophic N-ideals, S-minimal neutrosophic N ideals, S-quasi minimal neutrosophic N-ideals and other interesting substructures are introduced and illustrated with examples.

## 3.1 Smarandache Neutrosophic semigroup

Now we proceed on to define the notion of Smarandache neutrosophic semigroup and define several interesting substructures. The notion of S-conjugate neutrosophic subsemigroups is introduced and explained by examples. Strong

S-neutrosophic ideals happen to be nice algebraic substructure of S-neutrosophic N-semigroups.

**DEFINITION 3.1.1:** *Let (S, o) be a neutrosophic semigroup. S is said to be a Smarandache neutrosophic semigroup (S-neutrosophic semigroup) if S contains a proper subset P such that (P, o) is a group. If (P, o) is a neutrosophic group we call (S, o) a Smarandache strong neutrosophic semigroup.*

It is important to note that in general all neutrosophic semigroups need not be S-neutrosophic semigroups. We illustrate our observations by an example.

***Example 3.1.1:*** Let $\langle Z^+ \cup I\rangle$ be a neutrosophic semigroup under multiplication. Clearly $(\langle Z^+ \cup I\rangle, \times)$ has no proper subset P such that $(P, \times)$ is a group. So $(\langle Z^+ \cup I\rangle, *)$ is not a S-neutrosophic semigroup.

***Note:*** The identity element alone is never taken as a proper subgroup. That is why in the above example 3.1.1 {1} is not a proper group under ×.

***Example 3.1.2:*** Let $\langle Z \cup I\rangle$ be a neutrosophic semigroup under multiplication ×, P = {1, –1; ×} is a S-neutrosophic semigroup as P is a group under multiplication.

***Example 3.1.3:*** Let $\langle Z_6 \cup I\rangle$ = {0, 1, 2, 3, 4, 5, I, 2I, 3I, 4I, 5I} be a neutrosophic semigroup under multiplication modulo 6. Take P = {1, 5} $\subset \{\langle Z_6 \cup I\rangle\}$. P is a group under multiplication modulo 6. So $\{\langle Z_6 \cup I\rangle\}$ is a S-neutrosophic semigroup.

Now we proceed on to define neutrosophic subsemigroup.

**DEFINITION 3.1.2:** *Let (S, o) be a neutrosophic semigroup. A proper subset P of (S, o) is said to be a neutrosophic subsemigroup if (P, o) is a neutrosophic semigroup.*

***Example 3.1.4:*** Let $\langle Z^+ \cup I\rangle$ be a neutrosophic semigroup under multiplication. Let P = $\{\langle 2Z^+ \cup I\rangle\}$ be a proper subset of $\langle Z^+ \cup$

I⟩. Clearly P is a neutrosophic subsemigroup under multiplication. It may also happen a neutrosophic semigroup can have proper subsets which are just subsemigroups.

To this end we have the following example.

***Example 3.1.5:*** Let $\langle Z^+ \cup I\rangle$ be a neutrosophic semigroup under multiplication. Clearly $P = 2Z^+$ is a subsemigroup under multiplication.

Thus a neutrosophic semigroup can have two substructures one is a neutrosophic substructure other is just a substructure.

Now we proceed on to define the notion of Smarandache neutrosophic subsemigroup.

**DEFINITION 3.1.3:** *Let (S, o) be a neutrosophic semigroup. A proper subset P of S is said to be Smarandache neutrosophic subsemigroup (S-neutrosophic subsemigroup) if P has a proper subset U such that (U, o) is group then we call (P, o) to be a S-neutrosophic subsemigroup.*

*Thus we can define (P, o) to be a S-neutrosophic subsemigroup, if P is itself a S- neutrosophic semigroup under the operation of S.*

**THEOREM 3.1.1:** *Let (S, o) be a neutrosophic semigroup. If S has a proper subset P such that (P, o) is a S-neutrosophic subsemigroup then we call (S, o) to be a S-neutrosophic semigroup.*

*Proof:* Given (S, o) is a neutrosophic semigroup which has a proper subset P such that (P, o) is a S-neutrosophic subsemigroup i.e. P contains a proper subset X such that (X, o) is a group as $(X, o) \subset (S, o)$ so (S, o) is a S-neutrosophic semigroup.

It is interesting and important to note that in general all neutrosophic semigroups need not in general have a proper S-neutrosophic subsemigroup.

We prove this by the following example.

***Example 3.1.6:*** Let $\langle Z^+ \cup I\rangle$ be a neutrosophic semigroup under multiplication. Clearly $\langle Z^+ \cup I\rangle$ has no S-neutrosophic subsemigroup.

Now we proceed on to define the notion of S-neutrosophic ideals of a neutrosophic semigroup. It is important to note that these ideals are in the first place neutrosophic.

**DEFINITION 3.1.4:** *Let (S, o) be a neutrosophic semigroup. A proper subset I of S is said to be Smarandache neutrosophic ideal (S-neutrosophic ideal) of (S, o) if the following conditions are satisfied.*

*(1) (I, o) is a S-neutrosophic subsemigroup.*
*(2) For every s $\in$ S and i $\in$ I, si and is are in I.*

***Example 3.1.7:*** Let $\langle Z \cup I\rangle$ be a neutrosophic semigroup under multiplication, I = $\{\langle 3Z \cup I\rangle\}$ is a neutrosophic ideal of $\langle Z \cup I\rangle$. It can be easily verified; $\langle Z \cup I\rangle$ has no S-neutrosophic ideals.

Now we give yet another example.

***Example 3.1.8:*** Let (S, $\times$) = {0, 1, 2, 3, 4, 5, I, 2I, 3I, 4I, 5I} be a neutrosophic semigroup under multiplication modulo 6. Take J = {0, 2, 4, 2I, 4I} $\subset$ S; Clearly J is a S- neutrosophic ideal of S.

It is still interesting to note that in general the order of a S-neutrosophic subsemigroup need not divide the order of the neutrosophic semigroup.

***Example 3.1.9:*** Consider the S-neutrosophic semigroup (S, o) given in the above example 3.1.8 i.e. S = {0, 1, 2, 3, 4, 5, I, 2I, 3I, 4I, 5I} a neutrosophic semigroup under multiplication modulo 6. P = {0, 2, 4, 2I, 4I} is a S-neutrosophic subsemigroup of S and o(P) $\nmid$ (o (S)) as o(P) = 5 and o(S) = 11 and 5 $\nmid$ 11.

Thus as in case of other algebraic structures we define the order of the S-neutrosophic semigroup.

**DEFINITION 3.1.5:** *Let (S, o) be a S-neutrosophic semigroup, the number of distinct elements in S is called the order of the S-neutrosophic semigroup and is denoted by o(S) or |S|. If the number of elements in S is finite we call (S, o) to be a finite S-neutrosophic semigroup. It the number of elements in (S, o) is infinite then we call (S, o) be an infinite S-neutrosophic semigroup i.e. o(S) = ∞.*

It is very surprising to see that when S is a neutrosophic semigroup of finite order in general the order of the S-neutrosophic subsemigroup does not divide the order of S. To this end we define S-Lagrange neutrosophic semigroup.

**DEFINITION 3.1.6:** *Let (S, o) be a finite neutrosophic semigroup. If the order of every S-neutrosophic subsemigroup divides the order of S then we call S to be a Smarandache Lagrange neutrosophic semigroup (S-Lagrange neutrosophic semigroup).*

*If on the other hand the order of atleast one S-neutrosophic subsemigroup divides order of S we call S a Smarandache weakly Lagrange neutrosophic semigroup (S-weakly Lagrange neutrosophic semigroup) (A S-neutrosophic subsemigroup P is Lagrange if o(P) / o(S)).*

*If the order of no S-neutrosophic subsemigroup divides the order of the neutrosophic semigroup S then we call S to be a Smarandache Lagrange free neutrosophic semigroup (S-Lagrange free neutrosophic semigroup).*

We now illustrate this with example before we prove a theorem.

***Example 3.1.10:*** Let (S, ×) = {$Z_9$ = {0, 1, 2, 3, 4, 5, 6, 7, 8, I, 2I, 3I, 4I, 5I, 6I, 7I, 8I}, be a neutrosophic semigroup of finite order under multiplication modulo 9. o(S) = 17, a prime number. Clearly T = {0, 1, I, 8, 8I} is a S-neutrosophic subsemigroup for X = {1, 8} is a group under multiplication modulo 9. But o(T) cannot divide the prime 17. So S is a S-Lagrange free neutrosophic semigroup.

**THEOREM 3.1.2:** *Let (S, o) be a finite neutrosophic semigroup of prime order. Then (S, o) is S-Lagrange free neutrosophic semigroup.*

*Proof:* If (S, o) has no S-neutrosophic subsemigroups we have nothing to prove. So first we assume o(S) = p, p a prime has a S-neutrosophic subsemigroup P with o(P) = n since n < p and p a prime (n, p) = 1 i.e. n $\not|$ p. Hence S is a S-Lagrange free neutrosophic semigroup.

Now we proceed on to define Smarandache Cauchy elements and Smarandache Sylow neutrosophic semigroup. Throughout we assume S is always a finite neutrosophic semigroup, unless we make a special mention of it.

**DEFINITION 3.1.7**: *Let (S, o) be a finite neutrosophic semigroup. (say o(S) = n). Suppose for every proper subset T of S such that (T, o) is a group; we have for every element x in T, $x^m$ = e (since S is finite) and if m / n then we call S to be Smarandache Cauchy neutrosophic semigroup (S-Cauchy neutrosophic semigroup).*

*If the condition is true atleast for a proper subset which is a group in S then we call S to be a Smarandache weakly Cauchy neutrosophic semigroup (S-weakly Cauchy neutrosophic semigroup). If for no proper set which is a group under 'o' we have $x^m$ = e and m $\not|$ n then we call S to be Smarandache Cauchy free neutrosophic semigroup (S-Cauchy free neutrosophic semigroup).*

We see the class of S-Cauchy free neutrosophic semigroup is non empty. This is evident by the following theorem.

**THEOREM 3.1.3:** *Let (S, o) be a S-neutrosophic semigroup of finite order say p, p a prime. Then (S, o) is a S-Cauchy neutrosophic free semigroup.*

*Proof:* Given o(S) = p, a prime. Let G ⊂ S be a group under 'o'. Clearly o(G) < p. so for any x ∈ G we have $x^m$ = e, m ≤ o(G) < p, but (m, p) = 1 so m $\not|$ p. Hence the claim.

Thus we have the powerful corollary.

**COROLLARY:** *All S-neutrosophic semigroups of prime order are S-neutrosophic Cauchy free semigroup.*

Now we illustrate this by the following example.

***Example 3.1.11:*** Let (S, o) = {⟨$Z_6 \cup I$⟩ = (0, 1, 2, 3, 4, 5, I, 2I, 3I, 4I, 5I) under multiplication modulo 5}. o(S) = 11, a prime. G = {1, 5} is a group and G ⊂ S now $5^2 = 1$ (mod 6) but 2 ∤ 11. Hence the claim.

**DEFINITION 3.1.8:** *Let (S, o) be a S-neutrosophic semigroup of finite order. If for every prime p such that $p^\alpha$ | o(S) and $p^{\alpha+1}$ ∤ o(S) we have a S-neutrosophic subsemigroup of order $p^\alpha$ then we call S to be a Smarandache Sylow neutrosophic semigroup (S-Sylow neutrosophic semigroup).*

*If S has atleast for a prime p with $p^\alpha$ | o(S) and $p^{\alpha+1}$ ∤ o(S) a S-neutrosophic subsemigroup of order $p^\alpha$ then we call S a Smarandache weak Sylow neutrosophic semigroup (S-weak Sylow neutrosophic semigroup).*

*If for every prime p such that $p^\alpha$ | o(S) and $p^{\alpha+1}$ ∤ o(S) we don't have a S-neutrosophic subsemigroup then we call S a Smarandache neutrosophic Sylow free semigroup (S-neutrosophic Sylow free semigroup).We call the S-neutrosophic subsemigroup P of order $p^\alpha$ with o(P) / o(S) and $p^{\alpha+1}$ ∤ o(S) to be a p-Sylow S-neutrosophic subsemigroup.*

***Example 3.1.12:*** Let (S, o) = {0, 1, 2, 3, 4, 5, 6, 7, I, 2I, 3I, 4I, 5I, 6I, 7I} be a neutrosophic semigroup under multiplication modulo 8. Clearly (S, o) is a S-neutrosophic semigroup. o(S) = 15 we have 3 / 15 and $3^2$ ∤ 15, 5 / 15 and $5^2$ ∤ 15.

We have to find out whether (S, o) has S-neutrosophic subsemigroups of order 3 and 5. Take P = {0, 1, 7, I, 7I} ⊂ S. P is a S-neutrosophic subsemigroup of order 5. But S has no S-neutrosophic subsemigroup of order 3. Thus S is only a S-weak Sylow neutrosophic semigroup.

Let us define two new notions called Smarandache neutrosophic super Sylow semigroup and Smarandache neutrosophic semi Sylow semigroup.

**DEFINITION 3.1.9:** *Let (S, o) be a S-neutrosophic semigroup of finite order. Suppose S is a S-Sylow neutrosophic semigroup and if in addition to this for every prime p such that* $p^{\alpha} \mid o(S)$ *and* $p^{\alpha+1} \nmid o(S)$ *we have in S a S-neutrosophic subsemigroup of order* $p^{\alpha+t}$, $t \geq 1$ *then we call S to be a Smarandache neutrosophic super Sylow semigroup (S-neutrosophic super Sylow semigroup).*

**DEFINITION 3.1.10:** *Let (S, o) be a finite S-neutrosophic semigroup. If for every prime p such that* $p^{\alpha} \mid o(S)$ *and* $p^{\alpha+1} \nmid o(S)$ *we have a S-neutrosophic subsemigroup of order t,* $t < \alpha$ *then we call S to be a Smarandache neutrosophic semi Sylow semigroup (S-neutrosophic semi Sylow semigroup).*

Interested reader is requested to construct examples for these. It is also easy to verify that every S-neutrosophic super Sylow semigroup is a S-neutrosophic semigroup and not conversely.

Now we proceed on to define the notion of Smarandache conjugate neutrosophic semigroups.

**DEFINITION 3.1.11:** *Let (S, o) be a finite S-neutrosophic semigroup. Let T and R be any two S-neutrosophic subsemigroups of S. We say T and R are Smarandache neutrosophic conjugate subsemigroups (S-neutrosophic conjugate subsemigroups) if we have non empty subsets* $T_1$ *in T and* $R_1$ *in R,* $T_1 \neq R_1$ *but* $T_1$ *and* $R_1$ *are subgroups under 'o' and* $T_1$ *and* $R_1$ *as groups are conjugate i.e. there exists y, x in S with* $xT_1 = R_1x$ *or* $(xR_1 = T_1x)$ *and* $yR_1 = T_1 y$ $(R_1y = y T_1$ *for* $y \in S)$.

**Remark:** Clearly $|R_1| = |T_1|$ but $|R|$ need not be equal to order of T. We may not have any other order relation between them.

***Example 3.1.13:*** Consider (S, o) a S-neutrosophic semigroup of finite order given by S = {0, 1, 2, …, 14, I, 2I, 3I, …, 14I} a semigroup under multiplication modulo 15.

Take P = {0, 1, 4, I, 4I} and T = {1, 14, I, 14I} two proper subsets of S which are S-neutrosophic subsemigroups. $P_1$ = {1, 4} and $T_1$ = {1, 14} we have $3 \in S$ with 3 (1, 4) = (3, 12) and (1, 14) 3 = (3, 42 (mod 15)) = (3, 12). So $3T_1 = P_1\, 3$. 6 (1, 4) = (6, 9) and (1, 14) 6 = (6, 9). Thus $6\, P_1 = T_1\, 6$.

It is interesting to note that we may have more than one x in S which satisfies the condition. We say P and T are S-conjugate neutrosophic semigroups. Further 9 (1, 4) = (9, 6) and (1, 14) 9 = (9, 6). Thus $9\, T_1 = P_1\, 9$. It is easily verified $12\, T_1 = P_1\, 12$.

Thus we see {3, 6, 9, 12, 3I, 6I, 9I, 12I} are the elements in S which make P and T S-conjugate. 0 trivially makes P and T S-conjugate so {0, 3, 6, 9, 12, 3I, 6I, 9I, 12I} serves as a S-neutrosophic semigroup.

Several interesting problems can be raised at this juncture and we would be defining a new notion Smarandache conjugating subset of S.

**DEFINITION 3.1.12:** *Let (S, o) be a S-neutrosophic semigroup. Suppose P and T be any two S-conjugate neutrosophic subsemigroups of S. i.e. $P_1 \subset P$ and $T_1 \subset T$ are such that $xP_1 = T_1x$ (or $P_1x = xT_1$).*

*Now* $V = \left\{ x \in S \,\middle|\, \begin{array}{l} xP_1 = T_1x \quad or \\ P_1x = xT_1 \end{array} \right\}$ *is a proper subset of S and it is defined as the Smarandache conjugating subset (S-conjugating subset) of S for (P, T).*

We saw in the example 3.1.13 the Smarandache conjugating subset happens to be S-neutrosophic subsemigroup of S.

**DEFINITION 3.1.13:** *Let (S, o) be a S-neutrosophic semigroup we say a pair of elements x, y $\in$ S is said to be Smarandache conjugate pair (S-conjugate pair) if we can find a pair of elements a, b $\in$ S that ax = yb (or by) (or xa = by or yb).*

***Example 3.1.14:*** Let (S, o) be a S-neutrosophic semigroup where S = {0, 1, 2, 3, 4, 5, I, 2I, 3I, 4I, 5I} the operation 'o' just multiplication modulo 6.

Take the pair (3, 5) $\in$ S we see (1, 3) $\in$ S is such that 3.5 $\equiv$ 1.3 (mod 6). So the given pair is a S-conjugate pair of the S-neutrosophic semigroup.

One can develop several other interesting properties about them, as we are more interested in introducing Smarandache neutrosophic N-semigroups we now proceed on to define the notion of Smarandache neutrosophic bisemigroups. All properties can be derived in case of S-strong neutrosophic semigroups without any difficulty.

## 3.2 Smarandache neutrosophic bisemigroup

In this section the notion of Smarandache neutrosophic bisemigroups are introduced. Several of the analogous definitions about S-neutrosophic bisemigroups can be derived as in case of neutrosophic bisemigroups as a simple exercise.

**DEFINITION 3.2.1:** *Let {$B_S = B_1 \cup B_2$, $\times$, o} be a nonempty set with two binary operations. $B_S$ is said to be a Smarandache neutrosophic bisemigroup (S-neutrosophic bisemigroup) if*

- *i.* *($B_1$, $\times$) is a S-neutrosophic semigroup*
- *ii.* *($B_2$, o) is a S-neutrosophic semigroup or neutrosophic semigroup or S-semigroup or a semigroup.*

*If both $B_1$ and $B_2$ are S-neutrosophic semigroups. We call $B_S$ a Smarandache neutrosophic strong bisemigroup (S-neutrosophic strong bisemigroup). If ($B_1$, $\times$) is a S-neutrosophic semigroup and ($B_2$, o) neutrosophic semigroup then $B_S$ is a Smarandache neutrosophic bisemigroup (S- neutrosophic bisemigroup).*

If $(B_1, \times)$ as a S-neutrosophic semigroup and $(B_2, o)$ a S-semigroup then we call B a strong S-neutrosophic bisemigroup.

Interested reader can find the relations between them. Now we proceed onto define the substructures of $B_S$.

**DEFINITION 3.2.2:** *Let $B_S = \{B_1 \cup B_2, \times, o\}$ be a S-neutrosophic bisemigroup. A proper subset $V_S$ of $B_S$ is said to be Smarandache neutrosophic subbisemigroup (S-neutrosophic subbisemigroup) if $V_S = V_1 \cup V_2$ with $V_1 = V_S \cap B_1$ and $V_2 = V_S \cap B_2$ and $(V_S, o, *)$ is a S-neutrosophic bisemigroup.*

***Example 3.2.1:*** Let $B = (B_1 \cup B_2, *, o)$ be a strong S-neutrosophic bisemigroup where $B_1 = \langle Z \cup I\rangle$ S-neutrosophic semigroup under multiplication and $B_2 = \{0\ 1\ 2\ 3\ 4\ 5, I, 2I, 3I, 4I, 5I\}$; S-neutrosophic semigroup multiplication modulo 6, so B is a S-neutrosophic bisemigroup.

***Example 3.2.2:*** Let $B = (B_1 \cup B_2, \times, o)$ be a S-neutrosophic bisemigroup where $B_1 = \{\langle Z \cup I\rangle\}$, a semigroup under multiplication and $B_2 = \{0, 1, 2, 3, 4, 5\}$ semigroup under multiplication modulo 6. B is a S-neutrosophic bisemigroup.

Now $T = T_1 \cup T_2$ be a proper subset of B where $T_1 = \{\langle 3Z \cup I\rangle\}$ and $T_2 = \{0, 2, 4\}$ a semigroup under multiplication modulo 6. T is a S-neutrosophic subbisemigroup. One can construct any number of such examples.

**DEFINITION 3.2.3:** *Let $B_S = (B_1 \cup B_2, *, o)$ be any S-neutrosophic bisemigroup. A proper subset $J = I_1 \cup I_2$ of $B_S$ is said to be a Smarandache neutrosophic biideal (S-neutrosophic biideal) of $B_S$ if*

- *i.* *$J = I_1 \cup I_2$ is a S-neutrosophic subsemigroup, where $I_1 = B_1 \cap J$ and $I_2 = B_2 \cap J$.*
- *ii.* *Each of $I_1$ and $I_2$ are S-ideals of $B_1$ and $B_2$ respectively.*

*Now we say the neutrosophic biideal J of $B_S$ is Smarandache maximal neutrosophic biideal (S-maximal neutrosophic biideal) if $I_1$ is a maximal S-ideal in $B_1$ and $I_2$ is a maximal S-ideal in $B_2$.*

*Similarly we say J of $B_S$ is Smarandache neutrosophic minimal biideal (S-neutrosophic minimal biideal) of $B_S$ if $I_1$ is a minimal S-ideal of $B_1$ and $I_2$ is a minimal S-ideal of $B_2$. Now if one of $I_1$ is maximal S-ideal and $I_2$ is not maximal S-ideal, we call the S-neutrosophic biideal to be a Smarandache neutrosophic quasi maximal biideal (S-neutrosophic quasi maximal biideal). In similar manner Smarandache neutrosophic quasi minimal biideal (S-neutrosophic quasi minimal biideal) is defined.*

Several other properties related to ideals can be extended in case of S-neutrosophic biideals with proper and appropriate modifications.

Now we will proceed onto define neutrosophic N-semigroups (N a positive integer greater than or equal to two) and S-neutrosophic N-semigroups.

## 3.3 Smarandache Neutrosophic N-semigroup

In this section we introduce the new notion of Smarandache neutrosophic N-semigroups their substructures.

**DEFINITION 3.3.1:** *Let $\{S_s(N), *_1, *_2, \ldots, *_N\}$ be a nonempty set with N-binary operations. $S_s(N)$ is said to be a Smarandache neutrosophic N-semigroup (S-neutrosophic N-semigroup) if the following conditions are satisfied.*

i. *$S_s(N) = S_1 \cup \ldots \cup S_N$ where each $S_i$ is a proper subset of $S_s(N)$.*
ii. *Each of the set $(S_i, *_i)$ is either a S-semigroup or a S-neutrosophic semigroup.*

We first illustrate these by the following examples.

***Example 3.3.1:*** Let $S_s(N) = (S_1 \cup S_2 \cup S_3, *_1, *_2, *_3)$ where $S_1 = \{\langle Z \cup I \rangle\}$, neutrosophic semigroup under multiplication, $S_2 =$

{0, 1, 2, 3, 4, 5} semigroup under multiplication modulo 6 and $S_3$ = {0, 1, 2, 3, I, 2I, 3I} neutrosophic semigroup under multiplication modulo 4. $S_s(N)$ is a S-neutrosophic 3-semigroup.

Now we see if each of the semigroups ($S_i$, $*_i$) are S-neutrosophic semigroups i.e. for i = 1, 2, …, N then we call $S_s(N)$ to be strong Smarandache neutrosophic N-semigroup (Strong S-neutrosophic N-semigroup).

It is clear that all strong S-neutrosophic N-semigroup is always a S-neutrosophic N-semigroup but a S-neutrosophic N-semigroup in general is not a strong S-neutrosophic N-semigroup evident by the very definition and example 3.3.1.

**DEFINITION 3.3.2:** *Let $S_s(N) = (S_1 \cup S_2 \cup \ldots \cup S_N, *_1, \ldots, *_N)$ be a S-neutrosophic N semigroup. If the number of distinct elements in $S_s(S)$ is finite we call $S_s(S)$ a finite S-neutrosophic N semigroup if the number of elements in $S_s(S)$ is infinite we say $S_s(S)$ is an infinite S-neutrosophic N-semigroup and denote the order by $o(S_N (S))$.*

We are more interested in working mainly with finite S-neutrosophic N-semigroups.

***Example 3.3.2:*** Let $S_N (S) = (S_1 \cup S_2 \cup S_3, *_1, *_2, *_3)$ be a S-neutrosophic 3 semigroup where $S_1$ = {0, 1, 2, …, 9}, semigroup under multiplication modulo 10, $S_2$ = {0, 1, 2, 3, I, 2I, 3I} neutrosophic semigroup under multiplication modulo 4 and $S_3 = \{Z_3 \times Z_2$ = {(0, 0), (1, 0), (0, 1), (1, 1), (2, 0), (2, 1)}. $S_N (S)$ is a S-neutrosophic 3-semigroup and $o(S_N (S)) = 22$.

Now we can give yet another definition S-neutrosophic N-semigroup.

**DEFINITION 3.3.3:** *Let $S_N (S) = (S_1 \cup S_2 \cup \ldots \cup S_N, *_1, \ldots, *_N)$ be a neutrosophic N-semigroup. $S_N (S)$ is said to be a Smarandache neutrosophic N-semigroup (S-neutrosophic N-semigroup) if $S_N (S)$ has a proper subset P where P is a neutrosophic N-group.*

Now we proceed on to define substructures in them.

**DEFINITION 3.3.4:** *Let $S_N(S) = (S_1 \cup S_2 \cup \ldots \cup S_N, *_1, \ldots, *_N)$ be a neutrosophic N-semigroup. A proper subset $P = \{P_1 \cup P_2 \cup \ldots \cup P_N, *_1, \ldots, *_N\}$ of $S_N(S)$ is said to be a Smarandache neutrosophic sub N-semigroup (S-neutrosophic sub N-semigroup) of $S_N(S)$ if and only if P with the binary operations $*_1, \ldots, *_N$ is a S-neutrosophic N-semigroup.*

**THEOREM 3.3.1:** *If $S_N(S) = (S_1 \cup S_2 \cup \ldots \cup S_N, *_1, \ldots, *_N)$ is a neutrosophic N-semigroup having a S-neutrosophic sub N-semigroup then $S_N(S)$ is a S-neutrosophic N-semigroup.*

*Proof:* Given $S_N(S) = (S_1 \cup S_2 \cup \ldots \cup S_N, *_1, \ldots, *_N)$ is a neutrosophic N-semigroup having $P = \{P_1 \cup P_2 \cup \ldots \cup P_N, *_1, \ldots, *_N\}$ to be S-neutrosophic sub N-semigroup of $S_N(S)$. Since P is a S-neutrosophic sub N-semigroup, P has a proper subset $T = (T_1 \cup T_2 \cup \ldots \cup T_N, *_1, *_2, \ldots, *_N)$ which is a neutrosophic N-group. Now $T \subset P \subset S_N(S)$, so $T \subset S_N(S)$, thus $S_N(S)$ is a S-neutrosophic N-semigroup.

Now we illustrate this by an example.

***Example 3.3.3:*** Let $S_N(S) = (S_1 \cup S_2 \cup S_3, *_1, *_2, *_3)$ be a neutrosophic 3-semigroup where

$S_1$ = $\{Z_{12}$, semigroup under multiplication modulo 12$\}$,
$S_2$ = {0, 1, 2, 3, 4, 5, I, 2I, 3I, 4I, 5I}, neutrosophic semigroup under multiplication modulo 6 and

$$S_3 = \left\{\begin{pmatrix} a & b \\ c & d \end{pmatrix} \mid a, b, c, d \in \{0, 1, I, 1 + I\}\right\}.$$

Take $P = (P_1 \cup P_2 \cup P_3, *_1, *_2, *_3)$ in $S_N(S)$ where $P_1 = \{0, 2, 4, 6, 8, 10\} \subset S_1$, $P_2 = \{0, 2I, 1, I, 4, 2, 4I\} \subset S_2$ and $P_3 = \left\{\begin{pmatrix} a & b \\ c & d \end{pmatrix} \mid |ad - bc| \neq 0\right\} \cup \left\{\begin{pmatrix} 0 & 0 \\ 0 & 0 \end{pmatrix}\right\}$. P is a S-neutrosophic sub 3-semigroup of $S_N(S)$.

Now we proceed on to define the notion of Smarandache neutrosophic N-ideal of $S_N(S)$.

**DEFINITION 3.3.5:** *Let $S_N (S) = (S_1 \cup S_2 \cup ... \cup S_N, *_1, ..., *_N)$ be a S-neutrosophic N-semigroup. A proper subset $I = (I_1 \cup I_2 \cup ... \cup I_N, *_1, ..., *_N)$ of $S_N (S)$ is called the Smarandache neutrosophic N-ideal (S- neutrosophic N-ideal) of $S_N (S)$ if the following conditions are satisfied.*

*i.* *I is a S-neutrosophic sub N-semigroup.*
*ii.* *Each $(I_i, *_i)$ is a S-ideal of $(S_i, *_i)$ i.e. for $1 \leq i \leq N$.*

We illustrate this by the following example.

***Example 3.3.4:*** Let $S_N(S) = (S_1 \cup S_2 \cup S_3, *_1, *_2, *_3)$ be a S-neutrosophic 3-semigroup; where

$S_1$ = $\langle Z \cup I \rangle$, neutrosophic semigroup under multiplication.
$S_2$ = $\{Z_{12}$, semigroup under multiplication modulo 12$\}$ and
$S_3$ = $\{0, 1, 2, 3, 4, 5, I, 2I, 3I, 4I, 5I\}$, neutrosophic semigroup under multiplication modulo 6.

Take $I = (I_1 \cup I_2 \cup I_3, *_1, *_2, *_3)$ where $I_1 = \{\langle 3Z \cup I \rangle\} \subset \langle Z \cup I \rangle$, $I_2 = \{0, 2, 4, 6, 8, 10\} \subset Z_{12} \subseteq S_2$ and $I_3 = \{0, 2, 4, 2I, 4I\} \subset S_3$, I is a S-neutrosophic N-ideal of $S_N(S)$.

*Note:* We can as in case of N-semigroup define maximal Smarandache neutrosophic N-ideal, minimal Smarandache neutrosophic N-ideal, quasi maximal Smarandache neutrosophic N-ideal and quasi minimal Smarandache neutrosophic N-ideal. The concept of Lagrange Smarandache neutrosophic N-semigroup and their extensions viz. Sylow Smarandache neutrosophic N-semigroup and Cauchy Smarandache neutrosophic N-semigroup which is different from the usual definitions is given in case of Smarandache neutrosophic N-semigroup.

**DEFINITION 3.3.6:** *Let $S_N (S) = (S_1 \cup S_2 \cup ... \cup S_N, *_1, ..., *_N)$ be a S-neutrosophic N-semigroup of finite order. A proper* *subset P = {$P_1 \cup P_2 \cup \ldots \cup P_N$, $*_1$, …, $*_N$} of $S_N$ (S) which is a S-neutrosophic sub N-semigroup is said to be a Lagrange Smarandache neutrosophic sub N-semigroup (Lagrange S-neutrosophic sub N-semigroup) if the o(P) / o($S_N$(S)).*

*If the order of every S-neutrosophic sub N-semigroup divides the order of $S_N(S)$ then we call $S_N(S)$ a Lagrange Smarandache neutrosophic N-semigroup (Lagrange S-neutrosophic N-semigroup) i.e. if every S-neutrosophic sub N-semigroup of $S_N(S)$ is a Lagrange S-neutrosophic sub N-semigroup.*

*If $S_N(S)$ has atleast one Lagrange S-neutrosophic sub N-semigroup then we call $S_N(S)$ a weakly Lagrange Smarandache neutrosophic N-semigroup (Weakly Lagrange S-neutrosophic N-semigroup). If $S_N(S)$ has no Lagrange S-neutrosophic sub N-semigroup then we call $S_N(S)$ to be a Lagrange free Smarandache neutrosophic N-semigroup (Lagrange free S-neutrosophic N-semigroup).*

We illustrate this by the following example.

***Example 3.3.5:*** Let $S_N(S) = (S_1 \cup S_2 \cup S_3, *_1, *_2, *_3)$ be a S-neutrosophic N-semigroup of finite order, where

$S_1$ = {$Z_{12}$, a semigroup under multiplication modulo 12},
$S_2$ = {0, 1, 2, 3, 4, 5, I, 2I, 3I, 4I, 5I}, neutrosophic semigroup under multiplication modulo 6 and
$S_3$ = {(x, y) | x, y ∈ {0, 1, 2, I, 2I}}, neutrosophic semigroup under component-wise multiplication modulo 3.

o(S) = 12 + 11 + 25 = 48, a composite number.

Let $P = (P_1 \cup P_2 \cup P_3, *_1, *_2, *_3)$ where $P_1 = S_1$, $P_2$ = {0, 2, 4, 2I, 4I, I} and $P_3$ = {(0, 1) (1, 0) (2, 0) (0, 0), (I, 0) (2 I, 0)}. Clearly P is a S-neutrosophic sub N-semigroup, o(P) = 12 + 6 + 6 = 24 and o(P) / o($S_N$ (S)) , i.e. 24 / 48.

Consider T = {$T_1 \cup T_2 \cup T_3$, $*_1$, $*_2$, $*_3$} where $T_1$ = (0, 2, 4, 6, 8, 10) ⊂ $Z_{12}$, $T_2$ = {0, 3, 3I} ⊂ $S_2$ and $T_3$ = {(I, 0) (2 I, 0) (0, 2) (1, 0) (2, 0) (0, 0)} ⊂ $S_3$. Now o(T) = 15, 15 $\not|$ 48. Thus we can say $S_N$ (S) in this example is only a weakly Lagrange S-neutrosophic N-semigroup.

Now we give an example of a Lagrange free S-neutrosophic N-semigroup.

***Example 3.3.6:*** Let $S_N$ (S) = ($S_1 \cup S_2 \cup S_3$, $*_1$, $*_2$, $*_3$) where

$S_1$ = {$Z_{10}$, semigroup under multiplication modulo 10},
$S_2$ = {0, 1, 2, 3, 4, 5, I, 2I, 3I, 4I, 5I}, neutrosophic semigroup under multiplication modulo 6 and
$S_3$ = {$Z_2 \times Z_5$ | {(a, b)} such a $\in Z_2$ and b $\in Z_5$}, a semigroup under component wise multiplication.

o($S_N$ (S)) = 31, a prime number.

Take T = {$T_1 \cup T_2 \cup T_3$, $*_1$, $*_2$, $*_3$}, a proper subset of $S_N$(S) where $T_1$ = {0, 1, 9} $\subset S_1$, $T_2$ = {0, 2, 2I, 4, 4I} $\subset S_2$ and $T_3$ = {(1, 1) (1, 4) (1, 3) (1, 2) (0, 0)} $\subset S_3$. Clearly T is a S-neutrosophic 3-semigroup for $T_1 \supset$ {1, 9} a group under multiplication modulo 10. $T_2 \supset$ (2, 4) a group under multiplication modulo 6 and $T_3 \supset$ {(1 1), (1, 4)} a group.

Now o(T) = 6, (6, 31) = 1. Thus $S_N$(S) has non trivial S-neutrosophic sub N-semigroups but since the order of $S_N$(S) is a prime 31, the order of no S-neutrosophic sub N-semigroup will divide the order of $S_N$(S).

Hence $S_N$(S) in this example is a Lagrange free S-neutrosophic N-semigroup.

In view of the above example we have the following theorem.

**THEOREM 3.3.2:** *All S-neutrosophic N -semigroups, $S_N$ (S) = ($S_1 \cup S_2 \cup ... \cup S_N$, $*_1$, ..., $*_N$), of prime order p are Lagrange free S- neutrosophic N-semigroups.*

*Proof:* If $S_N$(S) has no proper S-neutrosophic sub N-semigroup nothing to prove. Given o($S_N$(S)) = a prime, p, so even if $S_N$(S) has non trivial proper S-neutrosophic sub N-semigroup we see their order cannot divide a prime. Hence the claim.

Now we proceed on to define the notion of strong Lagrange Smarandache neutrosophic N-semigroup.

**DEFINITION 3.3.7:** *Let $S_N(S) = (S_1 \cup S_2 \cup \ldots \cup S_N, *_1, \ldots, *_N)$ be a S-neutrosophic N-semigroup of finite order. If every proper subset $\{P_1 \cup P_2 \cup \ldots \cup P_N, *_1, \ldots, *_N\}$ of $S_N(S)$ which are N-groups are such that $o(P) / o(S_N(S))$ then we define $S_N(S)$ to be a strong Lagrange Smarandache neutrosophic N-semigroup (strong Lagrange S-neutrosophic N-semigroup).*

We do not see any relation between Lagrange S-neutrosophic N-semigroup or strong Lagrange S-neutrosophic N-semigroup. The only common feature between these two definitions is that when the S-neutrosophic N-semigroup is of prime order they are Lagrange free.

Strongly free Lagrange Smarandache neutrosophic N-semigroup (defined analogous to free Lagrange Smarandache neutrosophic N-semigroups by replacing sub N-semigroups by N-groups). All prime order S-neutrosophic N-semigroups are also strongly free Lagrange S-neutrosophic N-semigroups.

Now we proceed on to define the notion of Sylow Smarandache neutrosophic N-semigroups.

**DEFINITION 3.3.8:** *Let $S_N(S)$ be a finite S-neutrosophic N-semigroup. If for each prime p such that that $p^\alpha / o(S_N(S))$ and $p^{\alpha+1} \nmid o(S_N(S))$, we have a proper S-neutrosophic sub N-semigroup, $P = \{P_1 \cup P_2 \cup \ldots \cup P_N, *_1, \ldots, *_N\}$ with $o(P) = p^\alpha$ then we call $S_N(S)$ a Sylow Smarandache neutrosophic N-semigroup (Sylow S-neutrosophic N-semigroup).*

*We call each P as the p-Sylow Smarandache neutrosophic sub N-semigroup (p-Sylow S-neutrosophic sub N-semigroup). If $S_N(S)$ has atleast one p-Sylow S-neutrosophic sub N-semigroup P then we call $S_N(S)$ to be a weakly Sylow Smarandache neutrosophic N-semigroup (weakly Sylow S-neutrosophic N-semigroup). If $S_N(S)$ has no p-Sylow S-neutrosophic sub N-semigroup then we call $S_N(S)$ to be a Sylow free Smarandache neutrosophic N-semigroup(Sylow free S-neutrosophic N-semigroup) .*

***Example 3.3.7:*** Let $S_N(S) = (S_1 \cup S_2 \cup S_3, *_1, *_2, *_3)$ be a S-neutrosophic N-semigroup where

$S_1$ = {$Z_9$, a semigroup under multiplication modulo 9},
$S_2$ = {0, 1, 2, 3, I, 2I, 3I}, a neutrosophic semigroup under multiplication modulo 4 and
$S_3$ = {$Z_4 \times Z_2$ = (x, y) | x $\in Z_4$ and y $\in Z_2$}, a semigroup under component wise multiplication.

o($S_N$ (S) ) = 24, $2^3$ / 24 and $2^4 \not|$ 24, 3 / 24, $3^2 \not|$ 24.

We see $S_N$ (S) cannot have a S-neutrosophic sub N-semigroup of order 3. Further it cannot have a S-neutrosophic sub N-semigroup of order 8, Thus ($S_N$ (S) is Sylow free. In view of this we make a little modified definition.

**DEFINITION 3.3.9:** *Let $S_N$ (S) = ($S_1 \cup S_2 \cup ... \cup S_N$, $*_1$, ..., $*_N$) be a finite S-neutrosophic N-semigroup. If P = $P_1 \cup ... \cup P_N$ is any proper subset of $S_N$ (S) which is a N-group such that if for every prime p / o($S_N$ (S)) we have N-group of order $p^t$ ($t \geq 1$ ) then we call $S_N(S)$ to be a strong Sylow Smarandache neutrosophic N-semigroup (strong Sylow S-neutrosophic N-semigroup).*

Interested reader can construct examples.

Now we define the notion of Cauchy Smarandache neutrosophic N-semigroup and special Cauchy Smarandache neutrosophic N-semigroup.

**DEFINITION 3.3.10:** *Let $S_N$ (S) = ($S_1 \cup S_2 \cup ... \cup S_N$, $*_1$, ..., $*_N$) be a S-neutrosophic N-semigroup of finite order say n. If for all $x_i \in S_N$ (S) such that $x_i^{m_i} = e_i$ we have $m_i$ / n then we call $S_N$ (S) to be a Cauchy Smarandache neutrosophic N-semigroup (Cauchy S-neutrosophic N-semigroup).*

*We call $x_i$ the Cauchy Smarandache elements (Cauchy S-elements) of the S-neutrosophic N-semigroup. If only the elements of the N-groups are alone Cauchy S-elements of $S_N$ (S) then we call $S_N$ (S) to be a special Cauchy Smarandache neutrosophic N-semigroup (special Cauchy S-neutrosophic N-semigroup).*

It is clear from the definition that Cauchy S-neutrosophic N-semigroup by the very definition is a special Cauchy S-neutrosophic N-semigroup and not conversely. Now we give a simple example to illustrate our definitions.

***Example 3.3.8:*** Let $S_N(S) = (S_1 \cup S_2 \cup S_3, *_1, *_2, *_3)$ be a S-neutrosophic 3-semigroup. Let

$S_1$ = $\{Z_8$, semigroup under multiplication modulo 8$\}$,
$S_2$ = $\{0, 1, 2, 3, I, 2I, 3I\}$, semigroup under multiplication modulo 4 and
$S_3$ = $\{Z_3 \times Z_3\}$, semigroup under component wise multiplication.

$o(S_N(S)) = 24$. Now $S_N(S)$ is a Cauchy S-neutrosophic N-semigroup for $3^2 \equiv 1 \pmod 8$, $5^2 \equiv 1 \pmod 8$. $(2, 1)^2 = (1, 1)$; $(2, 2)^2 = (1\ 1)$, thus 2 / 24.

Now we proceed on to define the notion of Smarandache conjugate neutrosophic sub N-semigroups.

**DEFINITION 3.3.11:** *Let $S_N(S) = (S_1 \cup S_2 \cup S_3 \cup \ldots \cup S_N, *_1, \ldots, *_N)$ be a S-neutrosophic N-semigroup. Let $P = \{P_1 \cup P_2 \cup \ldots \cup P_N, *_1, \ldots, *_N\}$ and $M = \{M_1 \cup M_2 \cup \ldots \cup M_N, *_1, \ldots, *_N\}$ be any two S-neutrosophic sub N-semigroups of $S_N(S)$.*

*Let $T = (T_1 \cup T_2 \cup \ldots \cup T_N, *_1, *_2, \ldots, *_N)$ and $L = (L_1 \cup L_2 \cup \ldots \cup L_N, *_1, *_2, \ldots, *_N)$ be the proper subset of P and M respectively which are N-groups.*

*If T and L are Smarandache conjugate as N-groups (S-conjugate as N-groups) i.e. each $T_i$ and $L_i$ are S-conjugate as S-semigroups then we say P and M are Smarandache conjugate neutrosophic sub N-semigroups (S-conjugate neutrosophic sub N semigroups) of $S_N(S)$.*

***Remark:*** For P and M to be S-conjugate neutrosophic sub N semigroups we do not demand $o(P) = o(M)$. Neither do we demand any such traits among the N-groups.

**Chapter Four**

# SMARANDACHE NEUTROSOPHIC LOOPS AND THEIR GENERALIZATION

This chapter introduces the notion of Smarandache neutrosophic loops (S- neutrosophic loops) and their generalizations like Smarandache neutrosophic biloops, Smarandache neutrosophic N-loops and the new class of neutrosophic N-loops. The new class of neutrosophic loops introduced in [51] are all S-neutrosophic loops of order 4t or 2(n + 1), n odd n ≥ 3. This chapter has three sections. In section 1 S-neutrosophic loops are introduced and several interesting properties are derived. In section 2 Smarandache neutrosophic biloops are introduced and some of its interesting properties derived. Section three just introduces the notion of S- neutrosophic N-loops and suggests the reader to develop all properties.

## 4.1 S-Neutrosophic Loops and their Properties

In this section we newly define the notion of S-neutrosophic loops and show a class of loops of order 4n are S-neutrosophic loop. For more about neutrosophic loops please refer [51]. We infact prove the existence of S-neutrosophic loop which are right or left alternative, Moufang, Bol, Bruck and WIP.

Now we proceed on to define the notion of Smarandache neutrosophic loops.

**DEFINITION 4.1.1:** *Let $\langle L \cup I\rangle$ be a neutrosophic loop $\langle L \cup I\rangle$ is said to be a Smarandache neutrosophic loop (S- neutrosophic loop) if $\langle L \cup I\rangle$ has a proper subset T which is a neutrosophic group.*

***Example 4.1.1:*** Let $\langle L \cup I\rangle = \{\langle L_5(3) \cup I\rangle = \{e, 1, 2, 3, 4, 5, eI, 2I, 3I, 4I, 5I\}$ be a neutrosophic loop. $\{e, 2, I, 2I\} \subset \langle L_5(3) \cup I\rangle$ is a neutrosophic group under the operations of $\langle L_5(3) \cup I\rangle$. So $\langle L_5(3) \cup I\rangle$ is a S-neutrosophic loop.

We have a very important result to state [51].

**THEOREM 4.1.1:** *Every neutrosophic loop in the new class of neutrosophic loops $\langle L_n \cup I\rangle$ is a S-neutrosophic loop.*

*Proof:* We see every neutrosophic loop $\langle L_n(m) \cup I\rangle$ in $\langle L_n \cup I\rangle$ is a S-neutrosophic loop. This is evident from the fact that $\langle L_n(m) \cup I\rangle$ contains a subset {tI, e, eI, t} which is given by the following table under the operations in $\langle L_n(m) \cup I\rangle$

| * | e | t | eI | tI |
|---|---|---|---|---|
| e | e | t | eI | tI |
| t | t | e | tI | eI |
| eI | eI | tI | eI | tI |
| tI | tI | eI | tI | eI |

This is a neutrosophic loop and $\{\langle e, t\rangle \cup I\}$ is neutrosophic group as $\langle e, t\rangle$ is a group. Hence the claim.

Find examples of neutrosophic loops which are not S-neutrosophic loops. Now we proceed on to define Smarandache neutrosophic subloop.

**DEFINITION 4.1.2:** *Let $\langle L \cup I\rangle$ be a neutrosophic loop. A proper subset P of $\langle L \cup I\rangle$ is said to be a Smarandache neutrosophic subloop (S-neutrosophic subloop) if P satisfies the following conditions.*

(1) *P is a neutrosophic subloop of $\langle L \cup I\rangle$.*
(2) *P contains a proper subset H such that H is a neutrosophic group or equivalently we can define a Smarandache neutrosophic subloop P, to be a Smarandache neutrosophic loop under the operations of $\langle L \cup I\rangle$.*

The following result is very interesting.

**THEOREM 4.1.2:** *If $\langle L \cup I\rangle$ is a neutrosophic loop having a proper subset P which is a S-neutrosophic subloop then $\langle L \cup I\rangle$ is itself a S-neutrosophic loop.*

*Proof:* Given $\langle L \cup I\rangle$ is a neutrosophic loop and $P \subset \langle L \cup I\rangle$ is a S-neutrosophic subloop of $\langle L \cup I\rangle$ i.e. P has a proper subset T, such that T is a neutrosophic group and is contained in P.

But $T \subset \langle L \cup I\rangle$ as $T \subset P \subset \langle L \cup I\rangle$, Thus $\langle L \cup I\rangle$ has a proper subset T, where T is a neutrosophic group, so $\langle L \cup I\rangle$ is a S-neutrosophic loop. Hence the claim.

Now it may happen that the S-neutrosophic loop $\langle L \cup I\rangle$ may contain finite number of elements then we say $\langle L \cup I\rangle$ is a finite S-neutrosophic loop.

If the number of elements in a S-neutrosophic loop is infinite then we call $\langle L \cup I\rangle$ to be an infinite S-neutrosophic loop.

Now in general the order of a S-neutrosophic subloop need not divide the order of the finite S-neutrosophic loop. To tackle this problem we define the notion of Lagrange S-neutrosophic subloop and their other Lagrange concepts.

**DEFINITION 4.1.3:** *Let $\langle L \cup I\rangle$ be a S-neutrosophic loop of finite order. P be a proper subset of $\langle L \cup I\rangle$; we say P is a S-Lagrange neutrosophic subloop if P is a S-neutrosophic subloop and $o(P) / o(\langle L \cup I\rangle)$. If every S-neutrosophic subloop of $\langle L \cup I\rangle$ happens to be S-Lagrange neutrosophic subloop then we call $\langle L \cup I\rangle$ to be a Smarandache Lagrange neutrosophic loop (S-Lagrange neutrosophic loop). If $\langle L \cup I\rangle$ has atleast one S-Lagrange neutrosophic subloop then we call $\langle L \cup I\rangle$ a weak*

*Smarandache Lagrange neutrosophic loop (Weak S-Lagrange neutrosophic loop). If ⟨L ∪ I⟩ has no S-Lagrange neutrosophic subloop then we call ⟨L ∪ I⟩ to be a free S-Lagrange neutrosophic loop.*

*It is clear every Lagrange S-neutrosophic loop is weakly S-neutrosophic loop and converse is not true.*

*The reader is expected to construct concrete examples of*

i. *Weakly Lagrange S-neutrosophic loop,*
ii. *Lagrange S-neutrosophic loop*
iii. *Lagrange free S-neutrosophic loop.*

*Characterize those S-neutrosophic loops, which are Lagrange S-neutrosophic loops. However we prove the following theorem.*

**THEOREM 4.1.3:** *Let ⟨L ∪ I⟩ be S-neutrosophic loop of finite order n and n a prime. Then ⟨L ∪ I⟩ is a Lagrange free S-neutrosophic loop.*

*Proof:* Given ⟨L ∪ I⟩ is a S-neutrosophic loop of finite order n, n a prime. So every S-neutrosophic subloop P of ⟨L ∪ I⟩ is such that (n, o(P)) = 1 i.e. o(P) $\not|$ n, so ⟨L ∪ I⟩ has no Lagrange S-neutrosophic subloop. Hence ⟨L ∪ I⟩ is a Lagrange free S-neutrosophic loop.

We have also given a class of S-neutrosophic loops. Now we proceed on to define the notion of p-Sylow S neutrosophic subloops and their extended properties.

**DEFINITION 4.1.4:** *Let ⟨L ∪ I⟩ be a finite S-neutrosophic loop. If p is a prime such that $p^{\alpha}$ / o (⟨L ∪ I⟩) and $p^{\alpha+1}$ $\not|$ o (⟨L ∪ I⟩) and if ⟨L ∪ I⟩ has a S-neutrosophic subloop P of order $p^{\alpha}$ then we call P a p-Sylow S-neutrosophic subloop of ⟨L ∪ I⟩. If ⟨L ∪ I⟩ has atleast one p-Sylow S-neutrosophic subloop then we call ⟨L ∪ I⟩ to be a weakly Sylow S-neutrosophic loop. If for every prime p such that $p^{\alpha}$ / o (⟨L ∪ I⟩) and $p^{\alpha+1}$ $\not|$ o (⟨L ∪ I⟩) we* *have a p-Sylow S-neutrosophic subloop (or order $p^{\alpha}$) then we call $\langle L \cup I\rangle$ a Sylow S-neutrosophic loop. If $\langle L \cup I\rangle$ has no p-Sylow S-neutrosophic subloop then we call $\langle L \cup I\rangle$ a Sylow free S-neutrosophic loop.*

We have the following facts from the definition. Every Sylow S-neutrosophic loop is weakly Sylow S-neutrosophic loop. Clearly the converse is not true. The reader is expected to construct examples of those new classes of Sylow-S-neutrosophic loops. Now we define S-homomorphism of S-neutrosophic loops.

**DEFINITION 4.1.5:** *Let $\langle L \cup I\rangle$ and $\langle V \cup I\rangle$ be any two S-neutrosophic loops. A map $\phi$ from $\langle L \cup I\rangle$ to $\langle V \cup I\rangle$ is said to be Smarandache neutrosophic homomorphism (S- neutrosophic homomorphism) if the following conditions are true.*

i. $\phi(I) = I$.
ii. $\phi/P : P \to W$.

*(where P is a S-neutrosophic subgroup in $\langle L \cup I\rangle$ and W is a S-neutrosophic subgroup of $\langle V \cup I\rangle$) is a neutrosophic group homomorphism.*

*It is important to note $\phi$ need not be defined on all the elements of $\langle L \cup I\rangle$ only for a choosen $P \subset \langle L \cup I\rangle$ where P is a neutrosophic group $\phi$ is defined. Also for varying P we can have varying S neutrosophic homomorphism of $\langle L \cup I\rangle$ to $\langle V \cup I\rangle$. Kernel of $\phi$ in general may not be a normal neutrosophic subgroup of P.*

Now we illustrate this situation by the following example.

***Example 4.1.2:*** Let $\langle L \cup I\rangle$ and $\langle P \cup I\rangle$ be two S-neutrosophic loops, where $\langle L \cup I\rangle = \langle L_5(3) \cup I\rangle$ and $\langle P \cup I\rangle = \langle L_7(2) \cup I\rangle$. Let S = {e, 3, 3I, eI} $\subset \langle L_5(3) \cup I\rangle$ be a neutrosophic group and T = {e, 5, 5I, eI} $\subset \langle L_7(2) \cup I\rangle$ be a neutrosophic group.

Define $\phi_s = \phi/S : S \to T$ by

$\phi_S (e) = e \quad \phi_s (eI) = eI$
$\phi_s (3) = 5 \quad \phi_s (3I) = 5I.$

It is easily verified $\phi_s$ is a S-neutrosophic group homomorphism. Clearly $\phi$ need not be defined on the totality of $\langle L \cup I \rangle$.

Now we proceed on to define the notion of S-neutrosophic Moufang loop, S-neutrosophic Bol loop, S-neutrosophic Bruck loop and so on.

**DEFINITION 4.1.6:** *Let $\langle L \cup I \rangle$ be a S neutrosophic loop, $\langle L \cup I \rangle$ is defined to be a Smarandache neutrosophic Bol loop (S-neutrosophic Bol loop) if every proper subset A of $\langle L \cup I \rangle$, which is a neutrosophic subloop, is a Bol loop, i.e. elements of A satisfy the Bol identity.*

The following theorem is left as an exercise for the reader to prove.

**THEOREM 4.1.4:** *Every neutrosophic Bol loop is a S-Bol loop but not conversely.*

*We give an example of a S-neutrosophic loop, which is not a S-Bol loop.*

***Example 4.1.3:*** *Let $\langle L_5(3) \cup I \rangle$ be a S-neutrosophic loop. Clearly $\langle L_5(3) \cup I \rangle$ is not a neutrosophic Bol loop.*

**DEFINITION 4.1.7:** *Let $\langle L \cup I \rangle$ be a S-neutrosophic loop. $\langle L \cup I \rangle$ is defined to be Smarandache neutrosophic Moufang loop (S-neutrosophic Moufang loop) if every proper subset A of $\langle L \cup I \rangle$ which is a neutrosophic subloop is a Moufang loop i.e., A satisfies the Moufang identity.*

It is left as an exercise for the reader to prove the following.

Prove the S-neutrosophic loop $\langle L_{45}(8) \cup I\rangle$ is not a S-neutrosophic Moufang loop.

Now we can in the similar way define S-neutrosophic Bruck loop, S-neutrosophic alternative loop and so on, using the notion of neutrosophic loops defined in [51].

***Example 4.1.4:*** Let $L = \langle L_n(m) \cup I\rangle$ be any S neutrosophic loop where n is a prime, then $\langle L_n(m) \cup I\rangle$ is of order $2(n + 1)$. It is easily verified;

i. $\langle L_n(m) \cup I\rangle$ is a S-Lagrange neutrosophic loop.
ii. $\langle L_n(m) \cup I\rangle$ Lagrange neutrosophic loop.

Thus we have a nontrivial class of infinite number of S-neutrosophic loops which are S-Lagrange neutrosophic loop and at the same time Lagrange neutrosophic loop.

It is left for the reader to find S-neutrosophic loops which are S-Lagrange neutrosophic loops.

***Example 4.1.5:*** Let $L_{15}(2) \in L_{15}$ be a loop of order 16. $\langle L_{15}(2) \cup I\rangle$ is a S-neutrosophic loop of order 32. Consider the neutrosophic subloop $\{H_1(3) \cup I\}$ = {e, 1, 4, 7, 10, 13, eI, I, 4I, 7I, 10I, 13I}. Clearly $o(\langle H_1(3) \cup I\rangle) \mid o\langle L_{15}(2) \cup I\rangle$ i.e. $12 \nmid 32$. So $\langle H_1(3) \cup I\rangle$ is not a S-Lagrange neutrosophic subloop of $\langle L_{15}(2) \cup I\rangle$.

Now consider P = {e, 3, eI, 3I}. P is a S-Lagrange neutrosophic subloop of $\langle L_{15}(2) \cup I\rangle$; i.e. o (P) / 32. Hence P is a weak S-Lagrange neutrosophic loop.

Let $P_1 = \{e, 3\} \subset \langle L_{15}(2) \cup I\rangle$. $P_1$ is a Lagrange subloop. Take $P_2 = \{e, 1, 4, 7, 10, 13\} \subset \langle L_{15}(2) \cup I\rangle$, $P_2$ is also a Lagrange subloop of $\langle L_{15}(2) \cup I\rangle$.

Now $o(P_1) / o(\langle L_{15}(2) \cup I\rangle)$ and $o(P_2) \nmid o(\langle L_{15}(2) \cup I\rangle$.

Clearly $P_1$ and $P_2$ are not neutrosophic loops.

*Note:* The study of these properties especially in S-loops and S-neutrosophic loops will certainly help researches to know more about the loops in general for we have very natural classes of loops.

Now we illustrate this by an example before we proceed on to define the notion of Cauchy elements in S-neutrosophic loops of finite order.

***Example 4.1.6:*** Let $\langle L_n(m) \cup I\rangle$ be a S-neutrosophic loop of finite order. Here n is taken as a prime. Then every S-neutrosophic subloop and every neutrosophic subloop are p-Sylow neutrosophic subloops. Thus $\langle L_n(m) \cup I\rangle$ is a Sylow S-neutrosophic loop.

It still interesting to note that we have infinite number of classes of Sylow S-neutrosophic loops given by $\langle L_n \cup I\rangle\}$, n varies only over primes and each $\langle L_n(m) \cup I\rangle \in \{\langle L_n \cup I\rangle\}$ is a Sylow S-neutrosophic loop. It is still important to note that all loops in the new class of S-neutrosophic loops in general are not Sylow S-neutrosophic loops.

The reader is expected to give an example of the same.

**DEFINITION 4.1.8:** *Let $\langle L \cup I\rangle$ be a S-neutrosophic loop of finite order.*

*Let $x \in \langle L \cup I\rangle$ if $x^t = 1$ and $t / o \langle L \cup I\rangle$ we call x a Cauchy element if $x^t = 1$ and $t / o (P)$ where $x \in P$ and P is a S-neutrosophic subloop of $\langle L \cup I\rangle$ then we call x to be a S-Cauchy element of $\langle L \cup I\rangle$.*

*Let $y \in \langle L \cup I\rangle$ if, $y^t = I$ for some t and if $t / o\langle L \cup I\rangle$ we call y to be a Cauchy neutrosophic element of $\langle L \cup I\rangle$. If $y \in P \subset \langle L \cup I\rangle$ and P a S-neutrosophic subloop of $\langle L \cup I\rangle$ and if $y^s = I$ and $s / o(P)$ we call y a S-Cauchy neutrosophic element of $L \cup I$. If every element in P is a S-Cauchy element or a S-Cauchy neutrosophic element for at least one S-neutrosophic subloop P then we call $\langle L \cup I \rangle$ to be a S-Cauchy neutrosophic loop.*

On similar lines we can define Cauchy neutrosophic loop $\langle L \cup I\rangle$.

***Example 4.1.7:*** Let $\langle L_n(m) \cup I\rangle$ be a S-neutrosophic loop of finite order n, n a prime. Every element in $\langle L_n(m) \cup I\rangle$ is either a S-Cauchy element or a S-Cauchy neutrosophic element. It is easily verified every element x in $\langle L_n(m) \cup I\rangle$ is such that $x^2 = I$ or $x^2 = I$ and 2 / 2 (n + 1). Hence the claim.

Thus we have a new class of S-neutrosophic loops which are both S-Cauchy as well as S-neutrosophic Cauchy i.e. $\{\langle L_n \cup I\rangle\}$, n a prime has infinite number of finite classes (i.e. $\{\langle L_n(m) \cup I\rangle\} \in \{\langle L_n \cup I\rangle\}$) which are both S-Cauchy and S-Cauchy neutrosophic S-loops.

It is left for the reader to analyze those loops not in the new class of S-neutrosophic loops to be S-Cauchy. It is still important to note even the condition n is a prime, is not essential as every element in the new class of S-neutrosophic loops happens to be either a S-Cauchy element or a S-Cauchy neutrosophic element. Now we proceed on to define the notion of strong S-neutrosophic loops which are Moufang, Bruck, Bol or alternative and so on. For this we define it only locally and not globally.

**DEFINITION 4.1.9:** *Let $\langle L \cup I\rangle$ be a S-neutrosophic loop. We say $\langle L \cup I\rangle$ is a strong Moufang S-neutrosophic loop if every proper S-neutrosophic subloop P of $\langle L \cup I\rangle$ is a Moufang loop i.e. every P satisfies the Moufang identity.*

*Note:* It is interesting and important to note that we don't demand the whole of $\langle L \cup I\rangle$ to satisfy the Moufang identity. If every S-neutrosophic subloop of $\langle L \cup I\rangle$ satisfies the Moufang identity it is enough and we call $\langle L \cup I\rangle$ to be a S-neutrosophic Moufang loop.

We do not know whether there exist a relation between strong Moufang S-neutrosophic loops and Moufang S-neutrosophic loops. Interested reader is requested to analyze this problem.

**DEFINITION 4.1.10:** *Let $\langle L \cup I\rangle$ be a S-neutrosophic loop. Let P be a proper subset of $\langle L \cup I\rangle$ if P is just a neutrosophic subloop* *of L and satisfies the Moufang identity then we call P to be a Moufang neutrosophic subloop.*

*If every neutrosophic subloop of ⟨L ∪ I⟩ happen to satisfy the Moufang identity then we call ⟨L ∪ I⟩ to be a neutrosophic Moufang loop.*

Here also we do not demand the totality of the neutrosophic loop to satisfy the Moufang identity.

***Example 4.1.8:*** Let {⟨$L_n$(m) ∪ I⟩} be a S-neutrosophic loop where n is a prime. Every S-neutrosophic subloop P of ⟨$L_n$(m) ∪ I⟩ satisfies the Moufang identity. Thus we see ⟨$L_n$(m) ∪ I⟩ is a strong S-neutrosophic Moufang loop.

It can be easily verified that all the elements of the S-neutrosophic loop ⟨$L_n$(m) ∪ I⟩ does not satisfy the Moufang identity.

It is also interesting to note that every neutrosophic subloops of ⟨$L_n$(m) ∪ I⟩, n a prime also satisfy the Moufang identity.

Thus ⟨$L_n$(m) ∪ I⟩ is both a S-neutrosophic Moufang loop as well as the neutrosophic Moufang loop. Thus we have a class of S-neutrosophic Moufang loops as well as neutrosophic Moufang loops.

The cardinality of this class is clearly infinite {⟨$L_n$ ∪ I⟩ | n a prime and each ⟨$L_n$ ∪ I⟩ contains ⟨$L_n$(m) ∪ I⟩ for varying m}. Like wise we can define S-neutrosophic Bruck loop, S neutrosophic Bol loop and so on.

**DEFINITION 4.1.11:** *Let ⟨L ∪ I ⟩ be a S-neutrosophic loop. We say ⟨L ∪ I⟩ is a strong Weak Inverse Properly (WIP) S-neutrosophic loop, if every proper S-neutrosophic subloop of ⟨L ∪ I⟩ satisfies the identity.*

$$(xy)z = e \text{ imply } x (yz) = e \text{ for all } x, y , z \in \langle L \cup I\rangle.$$

*We have a new class of S-neutrosophic WIP loop.*

***Example 4.1.9:*** Let $\langle L_n(m) \cup I\rangle$ be a S-neutrosophic loop. $\langle L_n(m) \cup I\rangle$ is a S-neutrosophic WIP loop if and only if $(m^2 - m +1) \equiv 0 \pmod{n}$. This has been proved in case of new class of neutrosophic loops [51].

*Note:* It is important to note that every triple in $\langle L_n(m) \cup I\rangle$ such that $(m^2 - m +1) \equiv 0$ (modn) satisfies the WIP. So the whole S-neutrosophic loop $\langle L_n(m) \cup I\rangle$ is trivially a WIP loop.

Now does the new class of S-neutrosophic loops $\langle L_n \cup I \rangle$ have any other S-WI P-loops.

***Example 4.1.10:*** Let $\langle L_n(m) \cup I \rangle$ be a S-neutrosophic loop of finite order n, a prime. $\langle L_n(m) \cup I\rangle$ is a S-WIP loop. Clearly every triple of $\langle L_n(m) \cup I\rangle$ does not satisfy the WIP-identity.

Next we proceed on to define the concept of S-left and S-right alternative neutrosophic loops.

We will define for right or left the other can be done by the reader.

**DEFINITION 4.1.12:** *Let $\langle L \cup I \rangle$ be a S-neutrosophic loop $\langle L \cup I\rangle$ is said to be a strong S-right alternative neutrosophic loop if $(xy)\, y = x\, (yy)$ for all $x, y \in P$ and for every S-neutrosophic subloop P of $\langle L \cup I\rangle$. Similarly we can define S-left alternative neutrosophic loop.*

***Example 4.1.11:*** Let $\langle L_n (m) \cup I \rangle$ be a S-neutrosophic loop n a prime. Clearly $\langle L_n (m) \cup I\rangle$ is a S-neutrosophic loop but every element of $(\langle L_n (m) \cup I\rangle)$ need not in general satisfy the right (left) alternative identity. We have an infinite class of S-neutrosophic loops, which are left and right alternative i.e. alternative.

Several other properties enjoyed by loops S-loops or neutrosophic loops can also be easily extended to the case of S-neutrosophic loops.

## 4.2 Smarandache neutrosophic biloop

In this section the new notion of Smarandache neutrosophic biloops are introduced and a few of its properties are defined. Now we proceed on to define the notion of Smarandache neutrosophic biloop.

**DEFINITION 4.2.1:** *Let $S(B), *_1, *_2)$ be a non empty set on which two binary operations $*_1$ and $*_2$ are defined. We say $(S(B), *_1, *_2)$ is a Smarandache neutrosophic biloop (S- neutrosophic biloop) if the following conditions are satisfied*

*(1) $S(B) = S(B_1) \cup S(B_2)$ where $S(B_1)$ and $S(B_2)$ are proper subsets of $S(B)$.*
*(2) $(S(B_1), *_1)$ is a S-neutrosophic loop.*
*(3) $(S(B_2), *)$ is a S-loop or a group.*

*Thus all S-neutrosophic biloops are essentially neutrosophic biloops or neutrosophic subbiloops of the S-neutrosophic biloop.*

We illustrate this definition by the following example.

***Example 4.2.1:*** Let $S(B) = \{S(B_1) \cup D(B_2), *_1, *_2\}$ where

$S(B_1) = \langle L_5(2) \cup I\rangle$ and
$S(B_2) = A_4$;
S(B) is a S-neutrosophic biloop.

***Example 4.2.2:*** Let $S(B) = \{S(B_1) \cup S(B_2), *_1, *_2\}$ where

$S(B_1) = \langle L_5(2) \cup I\rangle$ and
$S(B_2) = (L_5(3))$
S(B) is a S-neutrosophic biloop.

We now proceed on to define the notion of S-neutrosophic subbiloop of a S-neutrosophic biloop.

**DEFINITION 4.2.2:** *Let S (B) = {S($B_1$) ∪ S($B_2$), $*_1$, $*_2$} be any neutrosophic biloop. We call a proper subset P of S(B) to be S-neutrosophic subbiloop, if P is a S-neutrosophic biloop under the operations $*_1$ and $*_2$.*

*Note:* P need not be even closed under the operation $*_1$ or $*_2$. So before we make a further definition we give an example.

***Example 4.2.3:*** Let S(B) = {S($B_1$) ∪ S($B_2$), $*_1$, $*_2$} be a neutrosophic biloop where S($B_1$) = $\langle L_7(3) \cup I \rangle$ and S($B_2$) = $A_4$.

Take P

$$= \left\{e, eI, 2, 2I, \begin{pmatrix} 1 & 2 & 3 & 4 \\ 1 & 2 & 3 & 4 \end{pmatrix}, \begin{pmatrix} 1 & 2 & 3 & 4 \\ 2 & 3 & 1 & 4 \end{pmatrix}, \begin{pmatrix} 1 & 2 & 3 & 4 \\ 3 & 1 & 2 & 4 \end{pmatrix}\right\}$$

be a proper subset of S(B).

Clearly (P, $*_1$) is not a loop or a group or a neutrosophic loop; likewise ($P_1$, $*_2$) is not a loop or a group or a neutrosophic loop. In fact P is not even closed under the binary operations $*_1$ or $*_2$. But take P =

$$\{e, eI, 2, 2I\} \cup \left\{\begin{pmatrix} 1 & 2 & 3 & 4 \\ 1 & 2 & 3 & 4 \end{pmatrix}, \begin{pmatrix} 1 & 2 & 3 & 4 \\ 2 & 3 & 1 & 4 \end{pmatrix}, \begin{pmatrix} 1 & 2 & 3 & 4 \\ 3 & 1 & 2 & 4 \end{pmatrix}\right\}$$

$= P_1 \cup P_2$;
where

$P_1$ = S($B_1$) ∩ P and
$P_2$ = S($B_2$) ∩ P; clearly

($P_1$, $*_1$) is a neutrosophic loop i.e. a neutrosophic subloop of S($B_1$) and ($P_2$, $*_2$) is a subgroup of S($B_2$) or a group under $*_2$. Thus P is a S-neutrosophic biloop.

In view of this we state the following theorem, the proof of which is left as an exercise for the reader to prove.

**THEOREM 4.2.1:** *Let S (B) = {S ($B_1$) $\cup$ S ($B_2$), $*_1$, $*_2$} be a neutrosophic biloop, P a proper subset of S (B) is a S-neutrosophic biloop; if P = $P_1 \cup P_2$ and $P_i = P \cap S(B_i)$ are neutrosophic subloops or groups of S($B_i$), ( i = 1, 2).*

It is interesting to note that if a neutrosophic biloop S(B) has a S-neutrosophic subbiloop then S(B) itself is a S-neutrosophic biloop. We make a little deviation in the definition of Lagrange S-neutrosophic biloops. It is also important to mention that all neutrosophic subbiloops in general are not S-neutrosophic subbiloops; however all S-neutrosophic subbiloops are neutrosophic subbiloops.

Now we give the definition of the order of a S-neutrosophic biloop. The order of the S-neutrosophic biloop B(S) is the number of distinct elements in B(S). If B(S) has finite number of elements then we call B(S) a finite S neutrosophic biloop. If the number of elements in B(S) is infinite, we say B(S) is an infinite neutrosophic biloop and denote the order by o(B(S)).

**DEFINITION 4.2.3:** *Let B(S) = {B($S_1$) $\cup$ B($S_2$), $*_1$, $*_2$} be a S-neutrosophic biloop of finite order. If P(S) = {$P_1 \cup P_2$, $*_1$, $*_2$} be a S-neutrosophic subbiloop of B(S) and if o(P(S)) / o(B(S)) then we call P(S) a Lagrange S-neutrosophic subbiloop.*

*If every S-neutrosophic subbiloop is Lagrange then we say S(B) is a Lagrange S-neutrosophic biloop. If S(B) has atleast one Lagrange S-neutrosophic subbiloop P then we say S(B) is a weakly S-neutrosophic biloop. If S(B) has no Lagrange S-neutrosophic subbiloop then we say S(B) is a Lagrange free S-neutrosophic biloop. Clearly all Lagrange S-neutrosophic biloops are weakly Lagrange S-neutrosophic biloops. But a weakly Lagrange S-neutrosophic biloop is not a Lagrange S-neutrosophic biloop.*

It is pertinent to mention here that we do not know whether the order of neutrosophic subbiloops of S(B) divide the order of S(B) and we are not interested about them when we study S-neutrosophic biloops.

However we define a new notion called strong Lagrange S-neutrosophic biloops of finite order.

**DEFINITION 4.2.4:** *Let $B(S) = \{B(S_1) \cup B(S_2), *_1, *_2\}$ be a S-neutrosophic biloop of finite order. If every S-neutrosophic subbiloop is Lagrange and if every neutrosophic subbiloop is Lagrange then we call B(S) to be a strong Lagrange S-neutrosophic biloop. All strong Lagrange S-neutrosophic biloops are Lagrange S-neutrosophic biloops but not conversely.*

Now we illustrate by an example.

***Example 4.2.4:*** Let $B(S) = \{B(S_1) \cup B(S_2), *_1, *_2\}$ be a S-neutrosophic biloop where $B(S_1) = \langle L_5(2) \cup I\rangle$ and $B(S_2) = \{g \mid g^6 = 1\}$. Clearly o[B(S)] = 18. Take $P = P_1 \cup P_2$ where $P_1 = \{e, eI, 3, 3I\}$ and $P_2 = \{g^3, 1\}$. Now o(P) = 6 and 6/18 so P is a Lagrange S-neutrosophic subbiloop.

Take $V = V_1 \cup V_2$ where $V_1 = \{e, 1, 2, 3, 4, 5\}$ and $V_2 = \{1, g^3\}$, o(V) = 8 and $o(V) \nmid o(S(B))$ so V is not a S Lagrange neutrosophic subbiloop.

Take $T = T_1 \cup T_2$ where $T_1 = \{1, e, 2, 3, 4, 5\}$ and $T_2 = \{1, g^2, g^4\}$, o(T) = 9, 9 / 18 so T is a Lagrange subbiloop, which is not even a subbiloop.

Suppose take $L = L_1 \cup L_2$ where $L_1 = \{e, eI, 3, 3I\}$ and $L_2 = \{1, g^2, g^4\}$, o(L) = 7 and L is a S-neutrosophic subbiloop which is not Lagrange as $7 \nmid 18$. Thus we see in the S-neutrosophic biloop we have several types of subbiloops. Interested reader can study and analyze them. But what we can say easily is, "Every S-neutrosophic biloop of order p, p a prime is a Lagrange free S-neutrosophic biloop".

Now we just give a simple example before we define Cauchy property for these S-neutrosophic biloops.

***Example 4.2.5:*** Let $S(B) = \{S(B_1) \cup S(B_2), *_1, *_2\}$ be a S-neutrosophic biloop where $S(B_1) = (\langle L_5(3) \cup I\rangle)$ and $S(B_2) = \{g \mid g^7 = 1\}$, o(S(B)) = 19. We have several S-neutrosophic subbiloops viz. $P = P_1 \cup P_2$, $P_1 = \{e, 3, eI, 3I\}$ and $P_2 = S(B_2)$ $o(P) \nmid 19$. Further we can have only trivial subgroups from $S(B_2)$, as it has no subgroups as $o(S(B_2)) = 7$, a prime order.

**DEFINITION 4.2.5:** *Let $S(B) = \{S(B_1) \cup S(B_2), *_1, *_2\}$ be a S-neutrosophic biloop of finite order. Suppose $x \in S(B)$ and $x^n = I$ and if $n / S(B)$ then x is called a Cauchy neutrosophic element of S(B). If $y \in S(B)$, is such that $y^m = 1$ and if $m / o(S(B))$, then y is called the Cauchy element of S(B). If every element of finite order or finite neutrosophic order happens to be a Cauchy element of S(B) then we call S(B) a Cauchy S-neutrosophic biloop.*

We now illustrate this situation by a simple example.

***Example 4.2.6:*** Let $S(B) = \{S(B_1) \cup (B_2), *_1, *_2\}$ where $S(B_1) = \{\langle L_7(3) \cup I\rangle\}$ and $S(B_2) = D_{2.4} = \{a, b, \mid a^2 = b^4 = 1, bab = a\}$. S(B) is a finite S-neutrosophic biloop. It is easily verified every element in S(B) is of order 2 or 4 only and o (S(B)) = 24, 2 and 4 divide 24. Infact 3/24 but S(B) has no element of order 3 or 6 that is why we have to define Cauchy element and Cauchy neutrosophic element in S- neutrosophic biloops. But our usual Cauchy theorem for finite groups is entirely different.

Now we see when we define p-Sylow structure also we do analyze in a very different manner.

**DEFINITION 4.2.6:** *Let $S(B) = \{S(B_1) \cup (B_2), *_1, *_2\}$ be a S-neutrosophic biloop of finite order. Let P be a prime such that $p^\alpha / o(S(B))$ and $p^{\alpha+1} \nmid o(S(B))$, if S(B) has a S-neutrosophic subbiloop P of order $p^\alpha$ then we call P a p-Sylow S-neutrosophic subbiloop of S (B). If for every prime p such that $p^\alpha / o(S(B))$ and $p^{\alpha+1} \nmid o(S(B))$ we have a p-Sylow S-neutrosophic subbiloop then we call S(B) to be a Sylow S-neutrosophic biloop. If S(B) has atleast one p-Sylow S-neutrosophic subbiloop then we call S(B) to be a weakly Sylow S-neutrosophic biloop. If S(B) has no p-Sylow S-neutrosophic subbiloops then we call S(B) a Sylow free S-neutrosophic* *biloop. It is a matter of routine to note that every Sylow S-neutrosophic biloop is a weakly Sylow S-neutrosophic biloop.*

However the converse is not true in general.

Several interesting results as in case of neutrosophic biloops can also be derived in case of S-neutrosophic biloops. Now we illustrate by two examples.

***Example 4.2.7:*** Let S(B) = {S($B_1$) ∪ ($B_2$), $*_1$, $*_2$} be a S-neutrosophic biloop where S($B_1$) = {⟨$L_7(3)$ ∪ I ⟩} and S($B_2$) = {g | $g^{10}$ = 1}.Clearly S(B) is a S-neutrosophic biloop of order 36. Now 3/36, $3^2$ / 36 but $3^3$ ∤ 36. 2 / 36, $2^2$ / 36 and $2^3$ ∤ 36. Now we study the p-Sylow neutrosophic subbiloops of S(B) {e, eI, 3I, 3} is the neutrosophic subloop. So S(B) has 2-Sylow neutrosophic subbiloops. Consider P = $P_1$ ∪ $P_2$ where $P_1$ = {e, eI, 3, 3I} and $P_2$ = {$g^2$, $g^4$, $g^6$, $g^8$, 1}. o(P) = 9 thus S (B) has a 3-Sylow neutrosophic subbiloop. Thus S (B) is only a weakly p-Sylow S-neutrosophic biloop.

***Example 4.2.8:*** Let S(B) = {S ($B_1$) ∪ ($B_2$), $*_1$, $*_2$}, be a S-neutrosophic biloop where S ($B_1$) = {⟨$L_{15}(7)$ ∪ I ⟩} and S ($B_2$) = $D_{2.20}$ = {a, b | $a^2$ = $b^{20}$ = 1, bab = a}. o (S(B)) = 72. Clearly $3^2$ / 72 and $3^3$ ∤ 72. We see $2^3$ / 72 and $2^4$ ∤ 72. The p-Sylow subloop of S(B) are given below. $P_1$ = {e, eI, 12, 12I} and = {$b^4$, $b^8$, $b^{12}$, $b^{16}$, 1}. o($P_1$ ∪ $P_2$) = 9, $P_1$ ∪ $P_2$ is a 3-Sylow neutrosophic subbiloop. Let V = $V_1$ ∪ $V_2$ ⊂ S(B) where $V_1$ = {e, eI, 3, 3I} and $V_2$ = {$b^5$, $b^{10}$, $b^{15}$, 1}. o ($V_1$ ∪ $V_2$) = 8. Thus S(B) is a p-Sylow neutrosophic biloop.

Thus we now proceed on to define the notion of Smarandache neutrosophic biloops of level II.

**DEFINITION 4.2.7:** *Let {⟨L ∪ I⟩ = $L_1$ ∪ $L_2$ $*_1$, $*_2$} be non empty set on which are defined two binary operations $*_1$ and $*_2$ such that*

*i.* *L = $L_1$ ∪ $L_2$ , $L_1$ and $L_2$ are is the proper subsets of L.*

*ii.* *$(L_1, *_1)$ is a S-neutrosophic loop.*
*iii.* *$(L_2, *_2)$ is a S-semigroup.*

*Then we call $\{\langle L \cup I\rangle = L_1 \cup L_2\ *_1, *_2\}$ to be a Smarandache neutrosophic biloop of level II (S-neutrosophic biloop of level II). In S-neutrosophic biloop of level II we take S-semigroup in place of a S-loop or group. If it so happens that $(L_2, *_2)$ is a neutrosophic group or a Smarandache neutrosophic semigroup then we call $\langle L \cup I\rangle$ to be Smarandache neutrosophic strong biloop (S- neutrosophic strong biloop).*

*Remark:* It is an easy consequence from the definition that every S-neutrosophic strong biloop is always a S-neutrosophic biloop but the converse is not true in general. We now illustrate this by the following example.

***Example 4.2.9:*** Let $\{\langle L \cup I\rangle = L_1 \cup L_2, *_1, *_2\}$ where $L_1 = \langle L_5(3) \cup I\rangle$, and $L_2 = S_4$. Clearly $\langle L \cup I\rangle$ is a S-neutrosophic biloop. Clearly $\langle L \cup I\rangle$ is not a S-neutrosophic strong biloop for $L_2$ is not a neutrosophic group or a S-neutrosophic semigroup.

***Example 4.2.10:*** Let $\{\langle L \cup I\rangle = L_1 \cup L_2, *_1, *_2\}$ where $L_1 = \langle L_7(3) \cup I\rangle$ a S-neutrosophic loop of order 16 and $L_2 = \langle Z_{10} \cup I\rangle$ a S-neutrosophic semigroup. Clearly L is a S-neutrosophic strong biloop.

Now as in case of loops or biloops we define the order of the S-neutrosophic biloop $\langle L \cup I\rangle$, to be the number of distinct elements in $\langle L \cup I\rangle$. If $\langle L \cup I\rangle$ has finite number of elements we call $\langle L \cup I\rangle$ to be a finite S-neutrosophic biloop, if $\langle L \cup I\rangle$ has infinite number of elements then we call $\langle L \cup I\rangle$ to be a infinite S-neutrosophic biloop. Both the examples given are S-neutrosophic biloops of finite order. We just give an example of a S-neutrosophic biloop of infinite order before which we define a substructure.

**DEFINITION 4.2.8:** *Let $\{\langle L \cup I\rangle = L_1 \cup L_2\ *_1, *_2\}$ be S-neutrosophic biloop of level II. Let P be a proper subset of $\langle L \cup$*

*I⟩ P is said to be Smarandache neutrosophic subbiloop (S-neutrosophic subbiloop) of ⟨L ∪ I⟩; if (P, $*_1$, $*_2$) is itself a S-neutrosophic biloop of level II under the operations of ⟨L ∪ I⟩.*

*Note:* P is not a loop under $*_1$ or $*_2$ or P need not be even closed under any one of the operations.

This is first illustrated by the following example.

***Example 4.2.11:*** Let {⟨L ∪ I⟩ = $L_1$ ∪ $L_2$ $*_1$, $*_2$} where $L_1$ = ⟨$L_5$ (3) ∪ I⟩ and $L_2$ = $S_3$ be a S-neutrosophic biloop. Take P = {e, eI, 3, 3I} ∪ {$A_3$}. P is a S-neutrosophic subbiloop. If

$$P = \{e, eI, 3, 3I, \begin{pmatrix} 1 & 2 & 3 \\ 1 & 2 & 3 \end{pmatrix}, \begin{pmatrix} 1 & 2 & 3 \\ 2 & 3 & 1 \end{pmatrix}, \begin{pmatrix} 1 & 2 & 3 \\ 3 & 1 & 2 \end{pmatrix} *_1, *_2\},$$

is not a loop under $*_1$ and is not a loop under $*_2$. But only P = $P_1$ ∪ $P_2$ is a S-neutrosophic subbiloop of ⟨L ∪ I⟩.

***Example 4.2.12:*** Let {⟨L ∪ I⟩ = $L_1$ ∪ $L_2$ $*_1$, $*_2$}, where $L_1$ = ⟨$L_7$(3) ∪ I⟩ be a S-neutrosophic loop and $L_2$ = ⟨Z ∪ I⟩, a semigroup under multiplication. ⟨L ∪ I⟩ is an infinite S-neutrosophic biloop as $L_2$ is an infinite S-neutrosophic semigroup.

Now we have seen both examples of infinite and finite S-neutrosophic biloops. We study some properties related with finite S-neutrosophic biloops. The term S-neutrosophic means Smarandache neutrosophic.

**DEFINITION 4.2.9:** *Let {⟨L ∪ I⟩ = $L_1$ ∪ $L_2$ $*_1$, $*_2$} be a S-neutrosophic biloop of level II of finite order. Let P be a proper subset of ⟨L ∪ I⟩ where P is a S-neutrosophic subbiloop of level II of ⟨L ∪ I⟩.*

*If o(P) / o(⟨L ∪ I⟩) then P is called the Lagrange S-neutrosophic subbiloop of ⟨L ∪ I⟩. If every S-neutrosophic subbiloop of level II of ⟨L ∪ I⟩ is Lagrange, then we say ⟨L ∪ I⟩ to be a Lagrange S-neutrosophic biloop. If ⟨L ∪ I⟩ has at least*

*one Lagrange S-neutrosophic subbiloop then we say ⟨L ∪ I⟩ is a weak Lagrange S-neutrosophic loop. If ⟨L ∪ I⟩ has no Lagrange S-neutrosophic biloop then we say ⟨L ∪ I⟩ is a Lagrange free S-neutrosophic biloop.*

***Example 4.2.13:*** Let {⟨L ∪ I⟩ = $L_1 \cup L_2$, $*_1$, $*_2$} where $L_1$ = ⟨$L_5(3)$ ∪ I⟩ and $L_2 = S_3$. o(⟨L ∪ I⟩) = 18.
Consider

$$P = \{e, 2, eI, 2I, \begin{pmatrix} 1 & 2 & 3 \\ 2 & 1 & 3 \end{pmatrix}, \begin{pmatrix} 1 & 2 & 3 \\ 1 & 2 & 3 \end{pmatrix}\}.$$

P is a S-neutrosophic subbiloop. o(P) = 6 and o(P) / o(⟨L ∪ I⟩). So P is a Lagrange S-neutrosophic subbiloop. Let

$$P' = \{e, 3, eI, 3I, \begin{pmatrix} 1 & 2 & 3 \\ 1 & 2 & 3 \end{pmatrix}, \begin{pmatrix} 1 & 2 & 3 \\ 2 & 3 & 1 \end{pmatrix}, \begin{pmatrix} 1 & 2 & 3 \\ 3 & 1 & 2 \end{pmatrix}\}.$$

P' is a S-neutrosophic subbiloop of ⟨L ∪ I⟩ but o(P') ╳ o(⟨L ∪ I⟩). So P' is not a Lagrange S-neutrosophic subbiloop of ⟨L ∪ I⟩. Thus ⟨L ∪ I⟩ is only a weak Lagrange S-neutrosophic biloop.

One can find many more examples and counter examples for these definitions. We need to mention only one property. If ⟨L ∪ I⟩ is a S-neutrosophic biloop of finite order say n and n is a prime then ⟨L ∪ I⟩ is a Lagrange free S-neutrosophic biloop.

***Example 4.2.14:*** Let {⟨L ∪ I⟩ = $L_1 \cup L_2$, $*_1$, $*_2$} where $L_1$ = ⟨$L_7(3)$ ∪ I⟩ and $L_2$ = S(3). ⟨L ∪ I⟩ is a S-neutrosophic biloop of order 43. Clearly ⟨L ∪ I⟩ is a Lagrange free S-neutrosophic biloop.

Now we proceed on to define p-Sylow S-neutrosophic subbiloops.

**DEFINITION 4.2.10:** *Let {⟨L ∪ I⟩ = $L_1 \cup L_2$ $*_1$, $*_2$} be a S-neutrosophic biloop of finite order. Let p be a prime such that $p^\alpha$ / o(⟨L ∪ I⟩) and $p^{\alpha+1}$ ╳ o(⟨L ∪ I⟩). If ⟨L ∪ I⟩ has a proper S-neutrosophic subbiloop of level II; P of order $p^\alpha$ then we call P* *the p-Sylow S-neutrosophic subbiloop of ⟨L ∪ I⟩. If for every prime p such that $p^{\alpha}$ / o(⟨L ∪ I⟩) and $p^{\alpha+1} \nmid$ o(⟨L ∪ I⟩) we have a p-Sylow S-neutrosophic subbiloop of level II then we call ⟨L ∪ I⟩ to be a Sylow S-neutrosophic biloop. If ⟨L ∪ I⟩ has atleast one p-Sylow S-neutrosophic subbiloop then we call ⟨L ∪ I⟩ to be a weak Sylow S-neutrosophic biloop. If ⟨L ∪ I⟩ has no p-Sylow S-neutrosophic subbiloop then we call ⟨L ∪ I⟩ to be a Sylow free S-neutrosophic biloop.*

***Example 4.2.15:*** Let {⟨L ∪ I⟩ = $L_1 \cup L_2$, $*_1$, $*_2$} be a finite S-neutrosophic biloop of level II where $L_1$ = ⟨$L_5(2)$ ∪ I⟩ and $L_2$ = {0, 1, 2, …, 6}, S-semigroup under multiplication modulo 7. For {1, 6} is a group in $L_2$. Clearly o(⟨L ∪ I⟩) = 19, a prime number. So ⟨L ∪ I⟩ has no p-Sylow S-neutrosophic subbiloop. Hence ⟨L ∪ I⟩ is a Sylow free S-neutrosophic biloop.

***Example 4.2.16:*** Let {⟨L ∪ I⟩ = $L_1 \cup L_2$, $*_1$, $*_2$} where

$L_1$ = ⟨$L_5(3)$ ∪ I⟩ and
$L_2$ = ⟨g | $g^6 = 1$⟩ be a S neutrosophic biloop.

o(⟨L ∪ I⟩) = 18, 3/ o(⟨L ∪ I⟩) and $3^2$ / o(⟨L ∪ I⟩) and $3^3 \nmid$ o(⟨L ∪ I⟩). 2 / 18 but $2^2 \nmid$ o(⟨L ∪ I⟩).

Clearly ⟨L ∪ I⟩ has no 2-Sylow S-neutrosophic subbiloop. It is easily verified that ⟨$L_5(3)$ ∪ I⟩ has only S-neutrosophic subloops of minimal and maximal order 4.

Thus ⟨L ∪ I⟩ has no S-neutrosophic subbiloop of order 9. So ⟨L ∪ I⟩ is a Sylow free S-neutrosophic biloop.

Thus even if the order of the S-neutrosophic biloop is n, n a composite number still we may not have a p-Sylow S-neutrosophic subbiloop i.e. ⟨L ∪ I⟩ can be Sylow free S-neutrosophic biloop.

Interested reader can derive results in this direction. Next we proceed on to define S-Cauchy element and S-Cauchy neutrosophic element of a S-neutrosophic biloop, ⟨L ∪ I⟩ = $L_1 \cup L_2$ $*_1$, $*_2$}.

**DEFINITION 4.2.11:** *Let $\{\langle L \cup I\rangle = L_1 \cup L_2 *_1, *_2\}$ be a S-neutrosophic biloop of level II of finite order. An element $x \in \langle L \cup I\rangle$ is said to be a Cauchy element of $\langle L \cup I\rangle$ if $x^t = 1$ and $t / o(\langle L \cup I\rangle)$. An element $y \in \langle L \cup I\rangle$ is said to be a Cauchy neutrosophic element of $\langle L \cup I\rangle$ if $y^m = I$ and $m / o(\langle L \cup I\rangle)$. If every element of $\langle L \cup I\rangle$ is either a Cauchy element or a Cauchy neutrosophic element of $\langle L \cup I\rangle$ then we call $\langle L \cup I\rangle$ to be Cauchy S-neutrosophic biloop.*

*If $\langle L \cup I\rangle$ has at least one Cauchy element and one Cauchy neutrosophic element then we call $\langle L \cup I\rangle$ to be a weak Cauchy S-neutrosophic biloop. If $\langle L \cup I\rangle$ has no Cauchy element or a Cauchy neutrosophic element then we call $\langle L \cup I\rangle$ to be Cauchy free S-neutrosophic biloop. If $\langle L \cup I\rangle$ has only Cauchy elements 'or' only Cauchy neutrosophic element then we call $\langle L \cup I\rangle$ to be a quasi Cauchy S-neutrosophic biloop. The term 'or' is used in the mutually exclusive sense.*

Relationship between these classes of S-neutrosophic biloop of level II is left for an innovative researcher.

Now we give some examples of these Cauchy S-neutrosophic biloops.

***Example 4.2.17:*** Let $\{\langle L \cup I\rangle = L_1 \cup L_2 *_1, *_2\}$, where

$L_1 = \langle L_5(3) \cup I\rangle$, and
$L_2 = \{0, 1, 2, 3\}$ S-semigroup under multiplication modulo 4.

$\langle L \cup I\rangle$ is a biloop of order 16.

Clearly $\langle L \cup I\rangle$ is a Cauchy S-neutrosophic biloop of level II, for every element in $\langle L \cup I\rangle$ is of either order 2.

Now we give an example of a Cauchy free S neutrosophic biloop.

***Example 4.2.18:*** Let us take $\{\langle L \cup I\rangle = L_1 \cup L_2, *_1, *_2\}$ where

$L_1 = \langle L_7(3) \cup I \rangle$ a neutrosophic loop of order 16 and
$L_2 = \{0, 1, 2, 3, 4$; S-semigroup under multiplication modulo $5\}$;

a cyclic group of order 5. Clearly $o(\langle L \cup I \rangle) = 21$ and $\langle L \cup I \rangle$ is Cauchy free S-neutrosophic biloop of level II .

Next we consider the following example.

***Example 4.2.19:*** Let $\{\langle L \cup I \rangle = L_1 \cup L_2, *_1, *_2\}$ where

$L_1 = \langle L_5 (3) \cup I \rangle$ and
$L_2 = \{0, 1, 2, 3, \ldots, 14$; S-semigroup under multiplication modulo $15\}$.

Clearly $\langle L \cup I \rangle$ is a S-neutrosophic biloop of level II of order 27. We see $\langle L \cup I \rangle$ has no Cauchy neutrosophic element but has only some Cauchy elements so $\langle L \cup I \rangle$ is a weak quasi Cauchy S-neutrosophic biloop.

Now we can define S-neutrosophic biloop which are Moufang, Bruck, Bol etc and give some illustrative examples of them.

**DEFINITION 4.2.12:** *Let $\{\langle L \cup I \rangle = L_1 \cup L_2, *_1, *_2\}$ be a S-neutrosophic biloop. We say $\langle L \cup I \rangle$ is a S-neutrosophic Moufang biloop, if every S-neutrosophic subbiloop of level II of $\langle L \cup I \rangle$ satisfies the Moufang identity. If $\langle L \cup I \rangle$ has atleast one S-neutrosophic subbiloop which satisfies the Moufang identity then we call $\langle L \cup I \rangle$ to be a weak S-neutrosophic Moufang biloop.*

*If $\langle L \cup I \rangle$ has no S-neutrosophic subbiloop which satisfies the Moufang identity then we say $\langle L \cup I \rangle$ is not a S-neutrosophic Moufang biloop.*

It is important to note that in general the whole loop $\langle L \cup I \rangle$ need not satisfy the Moufang identity. We try to illustrate them with examples.

***Example 4.2.20:*** Let $\{\langle L \cup I\rangle = L_1 \cup L_2, *_1, *_2\}$, where

$L_1 = \langle L_7 (3) \cup I\rangle$ and
$L_2 = \{0, 1, 2, 3, \ldots, 14$; S-semigroup under multiplication modulo $15\}$.

$\langle L \cup I\rangle$ is a S-neutrosophic biloop of finite order.

It is easily verified that every S-neutrosophic subbiloop of $\langle L \cup I\rangle$ satisfies the Moufang identity. Thus $\langle L \cup I\rangle$ is a S-neutrosophic Moufang biloop. But every element of $\langle L \cup I\rangle$ does not satisfy the Moufang identity.

Now we can on similar lines define S-neutrosophic Bol biloop, S-neutrosophic Bruck biloop, S-neutrosophic WIP biloop, S-neutrosophic alternative biloop, S-neutrosophic right (left) alternative biloops.

However here we give only examples of them.

***Example 4.2.21:*** Let $\{\langle L \cup I\rangle = L_1 \cup L_2, *_1, *_2\}$, where

$L_1 = \langle L_7(3) \cup I\rangle$ and
$L_2 = \{0, 1, 2, 3, \ldots, 11$; S-semigroup under multiplication modulo $12\}$.

$\langle L \cup I\rangle$ is a S-neutrosophic Bruck biloop.

***Example 4.2.22:*** Let $\{\langle L \cup I\rangle = L_1 \cup L_2, *_1, *_2\}$, where

$L_1 = \langle L_{13}(7) \cup I\rangle$ and
$L_2 = \{0, 1, 2, 3, 4, 5$; S-semigroup under multiplication modulo $6\}$.

$\langle L \cup I\rangle$ is a S-neutrosophic Bol biloop.

Before we give examples of S-neutrosophic alternative, (left) (right) biloop we make the following definition.

**DEFINITION 4.2.13:** *Let $\{\langle L \cup I\rangle = L_1 \cup L_2, *_1, *_2\}$, be a S-neutrosophic biloop of level II, $\langle L \cup I\rangle$ is said to be S-neutrosophic strong Moufang biloop if every element of $\langle L \cup I\rangle$ satisfies the Moufang identities.*

It is clear that all S-neutrosophic Moufang biloops in general need not be S-neutrosophic strong Moufang biloop. This is evident from the following example.

***Example 4.2.23:*** Let $\{\langle L \cup I\rangle = L_1 \cup L_2, *_1, *_2\}$ be a S-neutrosophic biloop where

$L_1 = \langle L_5(3) \cup I\rangle$ and
$L_2 = \{Z_{14}$, S-semigroup under multiplication modulo 14}.

Now $\langle L \cup I\rangle$ is a S-neutrosophic Moufang biloop, but it is easily verified $\langle L \cup I\rangle$ is not a S-neutrosophic strong Moufang biloop.

However the following theorem is easily verified.

**THEOREM 4.2.2:** *Every S-neutrosophic strong Moufang biloop $\langle L \cup I\rangle$ is a S-neutrosophic Moufang biloop. But the converse is not true.*

*Proof:* One side of the proof follows directly from the definition. To prove the other part the above 4.2.33 example will suffice.
We give example of S-neutrosophic strong left alternative biloops, S-neutrosophic strong right alternative biloops and S-neutrosophic strong WIP-biloop.

***Example 4.2.24:*** Let $\{\langle L \cup I\rangle = L_1 \cup L_2; *_1, *_2\}$, where

$L_1 = \langle L_n(m) \cup I\rangle$ with $(m^2 - m + 1) \equiv 0 (\text{mod } n)$, a S-neutrosophic loop and
$L_2 = \{Z_{12}$, S-semigroup under multiplication modulo 12}.

Clearly $\langle L \cup I\rangle$ is a S-neutrosophic strong WIP-biloop.

***Example 4.2.25:*** Let $\{\langle L \cup I\rangle = L_1 \cup L_2; *_1, *_2\}$; where

$L_1 = \langle L_n(2) \cup I\rangle$ and
$L_2 = \{Z_{24}$; S-semigroup under multiplication modulo 24$\}$.

$\langle L \cup I\rangle$ is a S-neutrosophic strong right alternative biloop.

***Example 4.2.26:*** Also $\{\langle L \cup I\rangle = L_1 \cup L_2; *_1, *_2\}$; where

$L_1 = \langle L_n(n - 1) \cup I\rangle$ and
$L_2 = G = \langle g \mid g^{12} = 1\rangle$.

$\langle L \cup I\rangle$ is a S-neutrosophic strong left alternative biloop.

***Example 4.2.27:*** Let $\langle L \cup I\rangle = \{L_1 \cup L_2, *_1, *_2\}$; where $L_1 = \langle L_{13}(11) \cup I\rangle$ and $L_2 = A_7$. $\langle L \cup I\rangle$ is only a S-neutrosophic (right) left alternative biloop but not a S-neutrosophic strong (right) left alternative biloop but not a S-neutrosophic strong (right) left alternative biloop. This is also S-neutrosophic WIP-biloop but is not a S-neutrosophic strong WIP-loop.

Several analogous properties can be defined.

Now we proceed on to define neutrosophic N-loops and Smarandache neutrosophic N-loops, $N \geq 2$. When $N = 2$ we get the neutrosophic biloop and the S-neutrosophic biloop.

## 4.3 Smarandache neutrosophic N-loops

In this section we just define the notion of Smarandache neutrosophic N-loop and illustrate with an example and suggest the reader to define a list of notions as in case of neutrosophic N-loops. For more about N-loops and neutrosophic N-loops refer [50, 51].

**DEFINITION 4.3.1:** *Let $\{\langle L \cup I\rangle = L_1 \cup L_2 \cup \ldots \cup L_N, *_1, \ldots, *_N\}$ be a nonempty set with N-binary operations. $\langle L \cup I\rangle$ is said to* *be a Smarandache neutrosophic N-loop (S-neutrosophic N-loop) if the following conditions are satisfied.*

i. *$\langle L \cup I\rangle = L_1 \cup L_2 \cup \ldots \cup L_N$ is such that each $L_i$ is a proper subset of $\langle L \cup I\rangle$.*
ii. *Some of $(L_i, *_i)$ are S-neutrosophic loops.*
iii. *Some of $(L_j, *_j)$ are just neutrosophic loops (and) or just S-loops.*
iv. *Some of $(L_k, *_k)$ are groups (and) or S-semigroups.*

***Example 4.3.1:*** Let $\langle L \cup I\rangle = \{L_1 \cup L_2 \cup L_3 \cup L_4, *_1, *_2, *_3, *_4\}$ where

$L_1$ = $\langle L_5(3) \cup I\rangle$,
$L_2$ = $L_7(2)$,
$L_3$ = $A_4$ and
$L_4$ = $\{Z_{10}$, semigroup under multiplication modulo $10\}$.

Clearly $\langle L \cup I\rangle$ is a S-neutrosophic 4-loop.

We can define order of a S-neutrosophic N-loop as in case of a neutrosophic N-loop. Also we define a proper subset P of $\{\langle L \cup I\rangle = L_1 \cup L_2 \cup \ldots \cup L_N, *_1, \ldots, *_N\}$ where $\langle L \cup I\rangle$ is a S-neutrosophic N-loop, is a Smarandache neutrosophic sub N-loop if P itself is a S-neutrosophic N-loop under the binary operations of $\langle L \cup I\rangle$.

Further P is not even closed under any of the binary operations $*_1, \ldots, *_N$. If fact $P_i = P \cap L_i$ for $1 \leq i \leq N$ and $\{P = P_1 \cup \ldots \cup P_N, *_1, \ldots, *_N\}$ is S-neutrosophic N-loop with each $P_i$ – a S-loop, a S-neutrosophic loop or a group or a S-semigroup. Clearly every neutrosophic loop in general need not be a S-neutrosophic N-loop. Further when N= 2 we get the S-neutrosophic biloop.

Now for N-subloops we make the following conditions as essential one. For finite S-neutrosophic N-loops we define the concept of Lagrange S-neutrosophic sub N-loop if the order S-neutrosophic sub N-loop divides the order of $\langle L \cup I\rangle$. On similar lines for neutrosophic N-loops we define Lagrange S-neutrosophic N-loops, p-Sylow S-neutrosophic sub N-loop, Sylow S-neutrosophic N-loop, Cauchy element, Cauchy neutrosophic element, Moufang S-neutrosophic N-loop, Bruck S-neutrosophic N-loop, WIP S-neutrosophic N-loop right (left) S-neutrosophic N-loop.

We as in case of neutrosophic N-loops define (N – t) deficit S-neutrosophic sub N-loop, $1 \leq t \leq N$. We for the Smarandache (N – t) deficit neutrosophic N-subloop define Lagrange concept, Sylow concept and Cauchy concept. Also we define when they are Moufang, Bol, alternative, WIP and so on, only if all the Smarandache (N – t) deficit neutrosophic N-subloop happens to be Moufang, Bol, alternative, WIP and so on respectively. All these are left as simple exercise for the reader to do with appropriate modifications. Also the reader is expected to give examples and enumerate some of the properties of them.

**Chapter Five**

# SMARANDACHE NEUTROSOPHIC GROUPOIDS AND THEIR GENERALIZATIONS

In this chapter we define the notion of Smarandache neutrosophic groupoids, Smarandache neutrosophic bigroupoids and Smarandache neutrosophic N-groupoids and give several of its interesting properties.

This chapter has two sections in the first section we define the notion of S-neutrosophic groupoids and give several related definitions. Section two of the chapter introduces first the concept of Smarandache neutrosophic bigroupoids and indicates briefly how the Smarandache neutrosophic substructures can be defined and analyzed.

The later part of the section introduces for the first time the notion of Smarandache neutrosophic N-groupoids and just define more particularizations and generalizations like Smarandache neutrosophic N-quasi loop, Smarandache neutrosophic N-quasi semigroups and Smarandache neutrosophic N-quasi groupoids.

Since this class of groupoids happen to be a class containing the associative structures, semigroups on one side and the non associative structure loops on the other side, we have not dealt elaborately with Smarandache N-groupoids. For this book deals elaborately with Smarandache neutrosophic N-loops and Smarandache neutrosophic N-semigroups.

## 5.1 Smarandache neutrosophic groupoids and their properties

In this section for the first time we introduce the notion of Smarandache neutrosophic groupoid (S- neutrosophic groupoids). We define special properties like Smarandache seminormal neutrosophic groupoid, Smarandache semi-conjugate neutrosophic subgroupoids and Smarandache neutrosophic P-groupoids and define several substructure as in case of S-neutrosophic loops.

Now we proceed on to define the notion of Smarandache neutrosophic groupoid.

**DEFINITION 5.1.1:** *A Smarandache neutrosophic groupoid (S-neutrosophic groupoid) $\{\langle G \cup I\rangle, *\}$ is a neutrosophic groupoid which has a proper subset S, $S \subset \langle G \cup I\rangle$ such that (S, *) is a neutrosophic semigroup.*

*The number of elements in the S-neutrosophic groupoid is called the order of the S-neutrosophic groupoid. If $\langle G \cup I\rangle$ has finite number of elements, then we say $\langle G \cup I\rangle$ is a finite S-neutrosophic groupoid, otherwise we say $\langle G \cup I\rangle$ is an infinite S-neutrosophic groupoid.*

Now we proceed onto define the notion of Smarandache neutrosophic subgroupoid.

**DEFINITION 5.1.2:** *Let $\{\langle G \cup I\rangle, *\}$ be a S-neutrosophic groupoid. A non empty subset H of $\langle G \cup I\rangle$ is said to be a Smarandache neutrosophic subgroupoid (S-neutrosophic subgroupoid) if H contains a proper subset $K \subset H$ such that (K, *) is a neutrosophic semigroup.*

The following theorem can be taken as simple exercises.

**THEOREM 5.1.1:** *Every neutrosophic subgroupoid of a S-neutrosophic groupoid need not in general be a S-neutrosophic subgroupoid.*

**THEOREM 5.1.2:** *Let $\langle G \cup I\rangle$ be a neutrosophic groupoid having a S-neutrosophic subgroupoid then $\langle G \cup I\rangle$ is a S-neutrosophic groupoid.*

Now we proceed on to define the notion of Smarandache commutative neutrosophic groupoid.

**DEFINITION 5.1.3:** *Let $\{\langle G \cup I\rangle, *\}$ be a S-neutrosophic groupoid, if every proper subset of $\langle G \cup I\rangle$ which is a neutrosophic semigroup is commutative then we call $\{\langle G \cup I\rangle, *\}$ a Smarandache commutative neutrosophic groupoid (S-commutative neutrosophic groupoid).*

*(It is interesting to note that whole of $\{\langle G \cup I\rangle, *\}$ need not be commutative).*

*We say G is a Smarandache weakly commutative neutrosophic groupoid if $\{\langle G \cup I\rangle, *\}$ has atleast one proper subset which is a neutrosophic semigroup is commutative.*

It is left as an exercise for the reader to prove.

**THEOREM 5.1.3:** *Every S-neutrosophic commutative groupoid is a S-weakly commutative neutrosophic groupoid.*

Now we proceed on to define the notion of Smarandache neutrosophic ideal of a S-neutrosophic groupoid.

**DEFINITION 5.1.4:** *A Smarandache left neutrosophic ideal (S-left neutrosophic ideal) A of a S-neutrosophic groupoid $\langle G \cup I\rangle$ satisfies the following conditions*

i. *A is a S-neutrosophic subgroupoid.*
ii. *$x \in \langle G \cup I\rangle$ and $a \in A$ then $x\, a \in A$.*

Similarly we can define Smarandache right neutrosophic ideal. If A is both a S-right neutrosophic ideal and S-left neutrosophic ideal then we say A is a Smarandache neutrosophic ideal (S-neutrosophic ideal) of $\langle G \cup I\rangle$.

Now we introduce a new notion called Smarandache semi normal neutrosophic groupoid.

**DEFINITION 5.1.5:** *Let V be a S-neutrosophic subgroupoid of $\langle G \cup I \rangle$. We say V is a Smarandache semi normal neutrosophic groupoid (S-semi normal neutrosophic groupoid) if*

$$aV = X \text{ for all } a \in \langle G \cup I \rangle$$
$$Va = Y \text{ for all } a \in \langle G \cup I \rangle$$

*where either X or Y is a S-neutrosophic subgroupoid of G but X and Y are both neutrosophic subgroupoids. V is said to be Smarandache normal neutrosophic groupoid (S-normal neutrosophic groupoid) if a V = X and Va = Y for all $a \in \langle G \cup I \rangle$, where both X and Y are S-neutrosophic subgroupoids of $\langle G \cup I \rangle$.*

We leave the following theorem as an exercise.

**THEOREM 5.1.4:** *Every S-normal neutrosophic groupoid is a S-seminormal neutrosophic groupoid and not conversely.*

We proceed on to define the notion of Smarandache semi conjugate neutrosophic subgroupoids.

**DEFINITION 5.1.6:** *Let $\langle G \cup I \rangle$ be a S-neutrosophic groupoid. H and P be two neutrosophic subgroupoids of $\langle G \cup I \rangle$. We say H and P are Smarandache semiconjugate neutrosophic subgroupoids (S- semiconjugate neutrosophic subgroupoids) of $\langle G \cup I \rangle$ if*

i. *H and P are S-neutrosophic subgroupoids of $\langle G \cup I \rangle$.*
ii. *H = xP or Px or*
iii. *P = xH or Hx for some $x \in \langle G \cup I \rangle$.*

We call two neutrosophic subgroupoids H and P of a neutrosophic groupoid $\langle G \cup I \rangle$ to be Smarandache conjugate neutrosophic subgroupoids (S- conjugate neutrosophic subgroupoids) of $\langle G \cup I \rangle$ if

i. *H and P are S-neutrosophic subgroupoids of* $\langle G \cup I \rangle$.
ii. $H = xP$ *or* $Px$ *and*
iii. $P = xH$ *or* $Hx$.

The following Theorem is left as an exercise for the reader to prove.

**THEOREM 5.1.5:** *Let* $\langle G \cup I \rangle$ *be a S-neutrosophic groupoid. If P and K are two S-neutrosophic subgroupoids of G which are S-conjugate then they are S-semi conjugate and the converse in general is not true.*

Now we proceed on to define a new notion called Smarandache inner commutative neutrosophic groupoid.

**DEFINITION 5.1.7:** *Let* $\langle G \cup I \rangle$ *be a S-neutrosophic groupoid. We say* $\langle G \cup I \rangle$ *is a Smarandache inner commutative (S-inner commutative) if every S-neutrosophic subgroupoid of G is inner commutative.*

Now we proceed on to define the notion of Smarandache neutrosophic Moufang groupoids.

**DEFINITION 5.1.8:** *Let* $\langle G \cup I \rangle$ *be a neutrosophic groupoid.* $\langle G \cup I \rangle$ *is said to be Smarandache Moufang neutrosophic groupoid (S-Moufang neutrosophic groupoid) if there exists H in* $\langle G \cup I \rangle$ *such that H is a S-neutrosophic subgroupoid of* $\langle G \cup I \rangle$ *and* $(x * y) * (z * x) = (x * (y * z)) * x$ *for all* $x, y, z \in H$.

If every S-neutrosophic subgroupoid of the neutrosophic groupoid $\langle G \cup I \rangle$ satisfies the Moufang identity then we call $\langle G \cup I \rangle$ a Smarandache strong neutrosophic Moufang groupoid.

On similar lines one can define Smarandache Bol neutrosophic groupoid, Smarandache strong Bol neutrosophic groupoid, Smarandache alternative neutrosophic groupoid

Smarandache strong alternative neutrosophic groupoid and so on.

Now we define Smarandache P-neutrosophic groupoid and Smarandache strong neutrosophic P-groupoid.

**DEFINITION 5.1.9:** *Let $\langle G \cup I\rangle$ be a S-neutrosophic groupoid we say $(\langle G \cup I\rangle, *)$ is a Smarandache neutrosophic P-groupoid (S-neutrosophic P-groupoid) if $\langle G \cup I\rangle$ contains a proper S-neutrosophic subgroupoid A such that $(x * y) * x = x * (y * x)$ for all $x, y \in A$.*

*We say $(\langle G \cup I\rangle, *)$ is Smarandache strong neutrosophic P-groupoid (S-strong neutrosophic P-groupoid) if every S-neutrosophic groupoid of $\langle G \cup I\rangle$ is a S-neutrosophic P-groupoid of $\langle G \cup I\rangle$.*

We can define the new notion called Smarandache direct product of neutrosophic groupoids.

**DEFINITION 5.1.10:** *Let $\langle G_1 \cup I\rangle, \langle G_2 \cup I\rangle, \ldots, \langle G_n \cup I\rangle$ be n-neutrosophic groupoids. We say $\langle G \cup I\rangle = \langle G_1 \cup I\rangle \times \ldots \times \langle G_n \cup I\rangle$ to be a Smarandache direct product of neutrosophic groupoids (S-direct product of neutrosophic groupoids) if $\langle G \cup I\rangle$ has a proper subset H of $\langle G \cup I\rangle$ which is a neutrosophic semigroup under the operations of $\langle G \cup I\rangle$. It is interesting and important to note that $\langle G_i \cup I\rangle$ need not be S-neutrosophic groupoids.*

One can define Smarandache neutrosophic homomorphism of S-neutrosophic groupoids.

**DEFINITION 5.1.11:** *Let $(\langle G_1 \cup I\rangle, *_1)$ and $(\langle G_2 \cup I\rangle, *_2 \rangle$ be any two S-neutrosophic groupoids. A map $\phi$ from $\langle G_1 \cup I\rangle$ to $\langle G_2 \cup I\rangle$ is said to be a Smarandache neutrosophic homomorphism (S-neutrosophic homomorphism) if $\phi$: $A_1 \to A_2$ $A_1 \subset (\langle G_1 \cup I\rangle$ and $A_2 \subset \langle G_2 \cup I\rangle$ are semigroups of $\langle G_1 \cup I\rangle$ and $\langle G_2 \cup I\rangle$ receptivity i.e. $\phi (a *_1 b) = \phi (a) *_2 \phi (b)$ for all $a, b \in A_1$ and $\phi (I) = I$.*

It is surprising to note that $\phi$ need not be even defined on whole of $(\langle G_1 \cup I\rangle)$.

We can define one to one S-homomorphism to be a Smarandache neutrosophic isomorphism.

## 5.2 Smarandache neutrosophic bigroupoids and their generalizations

In the first part of this section we introduce the notion of Smarandache neutrosophic bigroupoids (S-neutrosophic bigroupoids) and indicate some of its properties the later part of this section defines the notion of Smarandache neutrosophic N-groupoids (S- neutrosophic N-groupoids) and suggest some of its properties and few particularization and generalizations.

However we do not deal deeply into the result for with the given definition the reader can easily develop a theory as in the case of S-Neutrosophic loops and S-neutrosophic semigroups. Now we proceed on to define the notion of S-neutrosophic bigroupoids.

**DEFINITION 5.2.1:** *Let $(S(G), *_1, *_2)$ be a non-empty set with two binary operations. S(G) is said to be a Smarandache neutrosophic bigroupoid (S-neutrosophic bigroupoid) if the following conditions are satisfied.*

*i.* *$S(G) = G_1 \cup G_2$ where $G_1$ and $G_2$ are proper subsets of S(G).*

*ii.* *At least one of $(G_1, *_1)$ or $(G_2, *_2)$ is a S-neutrosophic groupoid.*

***Example 5.2.1:*** Let $\{S(G) = G_1, G_2, *_1, *_2\}$ where $G_1 = \langle Z_{12} \cup I\rangle$ the S-neutrosophic groupoid under the operation $*_1$, $a *_1 b = 8a + 4b \pmod{12}$ for $a, b \in Z_{12}$ and $G_2 = \langle Z_3 \cup I\rangle$, the neutrosophic groupoid under the operation $*_2$; $a *_2 b = 2a + b \pmod 3$. S(G) is a S-neutrosophic bigroupoid.

The order of the S-neutrosophic bigroupoid as in case of other S-neutrosophic algebraic structures is finite if it has finite number of distinct elements. Otherwise it is infinite. We can define Smarandache strong neutrosophic bigroupoid.

**DEFINITION 5.2.2:** *Let $S(G) = \{G_1 \cup G_2, *_1, *_2\}$ we say $S(G)$ to be a Smarandache strong neutrosophic bigroupoid (S-strong neutrosophic bigroupoid) if $S(G) = G_1 \cup G_2$ is such that $G_1$ and $G_2$ are proper subsets of $S(G)$ and both $(G_1, *_1)$ and $(G_2, *_2)$ are S-neutrosophic groupoids. Clearly all S-strong neutrosophic bigroupoids are S-neutrosophic bigroupoids.*

While defining the sub-structures the following are more interesting. One S-strong neutrosophic bigroupoid can have 5 types of subbigroupoids viz.

1. S-strong neutrosophic subbigroupoids
2. S-neutrosophic bigroupoids
3. S-bigroupoids
4. Neutrosophic bigroupoids
5. Bigroupoids.

It is to be noted that a S-neutrosophic bigroupoid cannot have S-strong neutrosophic subbigroupoids. Likewise we cannot have several types of biideals.

A S-strong neutrosophic bigroupoid can have only S-strong neutrosophic bi-ideals and a S-neutrosophic bigroupoid can have only S-neutrosophic biideals.

We define Smarandache super strong neutrosophic bigroupoid as follows:

**DEFINITION 5.2.3:** *Let $S(B) = (G_1 \cup G_2, *_1, *_2)$ be a non-empty set with two binary operations. $S(B)$ is called a Smarandache super strong neutrosophic bigroupoid (S-super strong neutrosophic bigroupoid) if $S(B) = G_1 \cup G_2$ where $G_1$ and $G_2$ are proper subsets of $S(B)$ and $(G_1, *_1)$ is S-neutrosophic loop and $(G_2, *_2)$ is a S-neutrosophic groupoid.*

All properties like defining S-Lagrange neutrosophic bigroupoids, S-Sylow neutrosophic bigroupoid, S-Moufang neutrosophic bigroupoid, S-neutrosophic biideals, and so on can be defined.

Now we proceed on to define the notion of Smarandache neutrosophic N-groupoids.

**DEFINITION 5.2.4:** *Let $S(G) = \{G_1 \cup G_2 \cup \ldots \cup G_N, *_1, *_2, \ldots, *_N\}$ be a non-empty set with N-binary operations. S(G) is called a Smarandache neutrosophic N-groupoid (S-neutrosophic N-groupoid) if $S(G) = G_1 \cup G_2 \cup \ldots \cup G_N$ is such that each $G_i$ is a proper subset of S(G).*

i. *Some of $(G_i, *_I)$ are S-neutrosophic groupoids.*
ii. *Some of $(G_j, *_j)$ are neutrosophic groupoids*
iii. *Rest of $(G_k, *_k)$ are S-neutrosophic semigroups or neutrosophic semigroups, $1 \leq i, j, k \leq N$.*

**DEFINITION 5.2.5:** *Let $\langle G \cup I\rangle = \{G_1 \cup G_2 \cup G_3 \cup \ldots \cup G_N, *_1, \ldots, *_N\}$ be a non-empty set with N-binary operations. We call $\langle G \cup I\rangle$ a Smarandache neutrosophic N-quasi loop (S-neutrosophic N-quasi loop) if at least one of the $(G_j, *_j)$ are S-neutrosophic loops. So a S-neutrosophic sub-N-quasiloop will demand one of the subsets $P_i$ contained $G_i$ to be a S-neutrosophic subloop.*

All properties pertaining to the sub-structures can be derived as in the case of S-neutrosophic N-groupoids. We define Smarandache N-quasi semigroups.

**DEFINITION 5.2.6:** *Let $\langle G \cup I\rangle = \{G_1 \cup G_2 \cup G_3 \cup \ldots \cup G_N, *_1, \ldots, *_N\}$ be a non-empty set with N-binary operations with each $G_i$ a proper subset of $\langle G \cup I\rangle$, $i = 1, 2, \ldots, N$. $\langle G \cup I\rangle$ is a Smarandache neutrosophic N-quasi semigroup (S-neutrosophic N-quasi semigroup) if some of the $(G_i, *_i)$ are S-neutrosophic loops and the rest are S- neutrosophic semigroups.*

*If in this definition some of the ($G_i$, *i) are S-neutrosophic loops and the rest just neutrosophic semigroup then we call* $\langle G \cup I\rangle$ *a weak S-neutrosophic N-quasi semigroup.*

*Note:* We do not have in the collection any neutrosophic groupoid or groupoid. Likewise we define Smarandache neutrosophic N-quasi groupoid as a non-empty set with N-binary operations $*_1, \ldots, *_N$ on $G_1, G_2, \ldots, G_N$ where $\langle G \cup I\rangle = \{G_1 \cup G_2 \cup G_3 \cup \ldots \cup G_N, *_1, \ldots, *_N\}$ where $(G_i, *_i)$ are either S-neutrosophic groups or S-neutrosophic groupoid or used in the mutually exclusive sense.

Further each $G_i$ is a proper subset of $\langle G \cup I\rangle$, i = 1, 2, …, N. Now all notions defined for S-neutrosophic N-loops can be easily extended to the class of S-neutrosophic groupoids.
Interested reader can work in this direction. To help the reader in chapter 7 if this book 25 problems are suggested only on S-neutrosophic groupoids and S-neutrosophic N-groupoids. All identities studied in case of S-neutrosophic loops and S-neutrosophic N-loops can also be defined and analyzed in case of S-neutrosophic groupoids.

If in this definition we have some of the $(G_t, *_t)$ are just groupoids or semigroups, then we call $(S(G), *_1, \ldots, *_N)$ to be S-weak neutrosophic N-groupoid. All the properties defined in case of S-neutrosophic groupoids can be defined with appropriate modification.

Concepts like Lagrange S-neutrosophic sub N-groupoids, p-Sylow S-neutrosophic sub N-groupoids, Lagrange S-neutrosophic N-groupoids, Sylow S-neutrosophic N-groupoids, weak Sylow S-neutrosophic N-groupoids, weak Lagrange S-neutrosophic N-groupoids, Sylow free neutrosophic N-groupoids, Lagrange free S-neutrosophic N-groupoids, S-Moufang neutrosophic N-groupoids, S-Bol neutrosophic N-groupoids can be defined.

Also, the notion of S-N conjugate S-sub N-groupoids can be defined. As in case of neutrosophic N-structures we can also define S-neutrosophic N-quasi loops, S-neutrosophic N-quasi semigroups. We have given over 25 problems for the interested researcher in S-neutrosophic N-groupoids and its related concepts, in chapter 7 of this book.

**Chapter Six**

# MIXED SMARANDACHE NEUTROSOPHIC STRUCTURES

This chapter for the first time introduces the notion of mixed Smarandache neutrosophic N-structures (mixed S- neutrosophic N-structures). This algebraic structure will be handy when we work with practical problems which takes it values from different algebraic structures. In such situations this mixed structures will be much useful.

We just sketch the properties for the innovative reader can do the job without any difficulty. The only criteria should be they must be familiar with the notions introduced in the book [50-1]. This chapter has only one section.

Now we proceed on to define Smarandache mixed neutrosophic N-structures.

**DEFINITION 6.1.1:** *Let $\langle M \cup I\rangle = \{M_1 \cup M_2 \cup \ldots \cup M_N, *_1, \ldots, *_N\}$ be a nonempty set with N-binary operations ($N \geq 5$). $\langle M \cup I\rangle$ is called the Smarandache mixed neutrosophic N-structure (S-mixed neutrosophic N-structure) if it satisfies the following conditions.*

i. *$\langle M \cup I\rangle = M_1 \cup M_2 \cup \ldots \cup M_N$ is such that each $M_i$ is a proper subset of $\langle M \cup I\rangle$, $1 \leq i \leq N$.*
ii. *$(M_i, *_i)$ for some i are S-neutrosophic groups.*
iii. *$(M_j, *_j)$ for some j are S-neutrosophic loops.*
iv. *$(M_k, *_k)$ for some k are S-neutrosophic semigroups.*

v. *$(M_t, *_t)$ for some t are S-neutrosophic groupoids.*

vi. *The rest of $(M_r, *_r)$ are groups or loops or groupoids or semigroups 'or' not used in the mutually exclusive sense.*

*We just first illustrate this by an example.*

***Example 6.1.1:*** Let $\{\langle M \cup I\rangle = M_1 \cup M_2 \cup M_3 \cup M_4 \cup M_5 \cup M_6, *_1, *_2, *_3, *_4, *_5, *_6\}$ where $M_1 = \langle L_5(3) \cup I\rangle$, $M_2 = \{1, 2, 3, 4, I, 2I, 3I, 4I\}$, $M_3 = \{\langle Z_6 \cup I\rangle$, is the S-neutrosophic semigroup under multiplication modulo 6}, $M_4 = \{\langle Z_8 \cup I\rangle \mid a * b = 3a + 5b \pmod 8$ for $a, b \in \langle Z_8 \cup I\rangle\}$ is the S-neutrosophic groupoid, $M_5 = A_5$ and $M_6 = S(3)$. $\langle M \cup I\rangle$ is a S-mixed neutrosophic N-structure.

If the number of distinct elements in the S-mixed neutrosophic N-structure is finite we call $\langle M \cup I\rangle$ to be finite, otherwise infinite. The S-mixed neutrosophic N-structure given in the above example is finite. We are more interested in finite S-mixed neutrosophic N-structure for by using them we can derive several interesting properties.

Now we proceed on to define the three types of sub N-structure as in case of mixed neutrosophic N-structure.

**DEFINITION 6.1.2:** *Let $\{\langle M \cup I\rangle = M_1 \cup M_2 \cup \ldots \cup M_N, *_1, \ldots, *_N\}$, where $\langle M \cup I\rangle$ is a S-mixed neutrosophic N-structure. A proper subset P of $\langle M \cup I\rangle$ is said to be S-mixed neutrosophic N-substructure if $P = \{P_1 \cup P_2 \cup \ldots \cup P_N, *_1, \ldots, *_N\}$ such that $P_i = P \cap L_i$ $1 \leq i \leq N$ and P under the operations of $\langle M \cup I\rangle$ is a S-mixed neutrosophic N-structure. On similar lines we define V a non empty proper subset of $\langle M \cup I\rangle$ to be a S-deficit mixed neutrosophic substructure if $V = \{V_1 \cup V_2 \cup \ldots \cup V_m, *_1, \ldots, *_m\}$ with $m < N$ happens to be a S-mixed neutrosophic sub structure where $*_1, \ldots, *_m \in \{*_1, \ldots, *_N\}$.*

Now we define dual S-mixed neutrosophic N-structure.

**DEFINITION 6.1.3:** *Let $\{\langle D \cup I\rangle = D_1 \cup D_2 \cup \ldots \cup D_N, *_1, *_2, \ldots, *_N\}$ where $\langle D \cup I\rangle$ is non empty set on which is defined N-binary operations ($N \geq 5$). We say $\langle D \cup I\rangle$ is a dual Smarandache mixed neutrosophic N-structure (dual S-mixed neutrosophic N-structure) if $\langle D \cup I\rangle$ satisfies the following conditions.*

- i. *$\langle D \cup I\rangle = D_1 \cup D_2 \cup \ldots \cup D_N$ where $D_i$ s are proper subsets of $\langle D \cup I\rangle$; $1 \leq i \leq N$.*
- ii. *Some of $(D_i, *_i)$ are S-loops.*
- iii. *Some of $(D_j, *_j)$ are S-semigroups.*
- iv. *Some of $(D_k, *_k)$ are S-groupoids.*
- v. *Some of $(D_t, *_t)$ are groups.*
- vi. *The rest of $(D_p, *_p)$ are neutrosophic loops or neutrosophic groups, or neutrosophic semigroups or neutrosophic groupoids 'or' not used in the mutually exclusive sense.*

Now we proceed on to give an example of it.

***Example 6.1.2:*** Let $\{\langle D \cup I\rangle = D_1 \cup D_2 \cup D_3 \cup D_4 \cup D_5, \cup D_6, *_1, *_2, *_3, *_4, *_5, *_6\}$, where $D_1 = L_5(3)$, $D_2 = S(3)$, $D_3 = \{a, b \in Z_{12}$, such that $a *_3 b = 2a + 4b \pmod{12}\}$, $D_4 = A_4$, $D_5 = \langle L_7(2) \cup I\rangle$ and $D_6 = \{0, 1, 2, 3, 4, 5, I, 2I, 3I, 4I, 5I$, under multiplication modulo 6$\}$. $\langle D \cup I\rangle$ is a dual S-mixed neutrosophic 6-structure.

We just define the substructure in it.

**DEFINITION 6.1.4:** *Let $\{\langle D \cup I\rangle = D_1 \cup D_2 \cup \ldots \cup D_N, *_1, *_2, \ldots, *_N\}$ be a dual S-mixed neutrosophic N-structure, a proper subset P of $\langle D \cup I\rangle$ is said to be a dual S-mixed neutrosophic sub N-structure of $\langle D \cup I\rangle$ if $P = \{P_1 \cup P_2, \cup P_3 \cup \ldots \cup P_N, *_1, \ldots, *_N\}$ where $P_i = P \cap D_i$, $1 \leq i \leq N$ and $P_i$ a substructure of $D_i$ under the operations of $*_i$ and P is also a dual S-mixed neutrosophic N-structure under the binary operations $*_1, \ldots, *_N$.*

As in case of mixed neutrosophic N-structures we define the notion of weak S-mixed neutrosophic N-structure and dual S-mixed neutrosophic N-structure.

Also we can define S-deficit mixed neutrosophic N-structure and dual S-deficit mixed neutrosophic N-structure. A major difference between S-mixed neutrosophic N-structure and mixed neutrosophic N-structure is we can in case of S-mixed neutrosophic N-structure can have just proper subsets P which are mixed neutrosophic structures. Now we define the Lagrange properties for S-mixed neutrosophic N-Structures.

**DEFINITION 6.1.5:** *Let* $\{\langle M \cup I\rangle = M_1 \cup M_2 \cup \ldots \cup M_N, *_1, \ldots, *_N\}$ *be a S-mixed neutrosophic N-structure of finite order. Let P be a proper subset of* $\langle M \cup I\rangle$, *be a S-mixed neutrosophic sub N-structure of* $\langle M \cup I\rangle$. *If* $o(P) / o(\langle M \cup I\rangle)$ *then we call P to be a Lagrange Smarandache mixed neutrosophic sub N-structure (Lagrange S-mixed neutrosophic sub N-structure). If every proper S-mixed neutrosophic sub N-structure is Lagrange then we call* $\langle M \cup I\rangle$ *itself to be a Lagrange S-mixed neutrosophic N-structure.*

*If* $\langle M \cup I\rangle$ *has atleast one Lagrange S-mixed neutrosophic sub N-structure then we call* $\langle M \cup I\rangle$ *to be weak Lagrange S-mixed neutrosophic N-structure. If* $\langle M \cup I\rangle$ *has no Lagrange S-mixed neutrosophic sub N-structure then we call* $\langle M \cup I\rangle$ *to be Lagrange free S-mixed neutrosophic sub N-structure.*

On similar lines we can define these concepts for weak S-mixed neutrosophic sub N structures and S deficit mixed neutrosophic sub N structures.

*This is left for the reader to develop.*

Now we proceed on to show that we have a class of S-mixed neutrosophic N-structure which are Lagrange free.
All S-mixed neutrosophic N-structure of finite order n, n a prime are Lagrange free S-mixed neutrosophic N-structures.

Now we proceed on to define Sylow concept to S-mixed neutrosophic N-structures.

**DEFINITION 6.1.6:** *Let* $\{\langle M \cup I\rangle = M_1 \cup M_2 \cup \ldots \cup M_N, *_1, \ldots, *_N\}$ *be a S-mixed neutrosophic N-structure of finite order. Let P* *be a prime such that $p^{\alpha}$ / o(⟨M ∪ I⟩) and $p^{\alpha+1}$ ∤ o(⟨M ∪ I⟩). If ⟨M ∪ I⟩ has a S-mixed neutrosophic sub N-structure of order $p^{\alpha}$, P then we call P to be a p-Sylow S-mixed neutrosophic sub N-structure. If for every prime p such that $p^{\alpha}$ / o(⟨M∪ I⟩) and $p^{\alpha+1}$ ∤ o(⟨M∪ I⟩) we have a p-Sylow S- mixed neutrosophic sub N-structure then we call ⟨M ∪ I⟩ to be a Sylow S-mixed neutrosophic N-structure. If ⟨M ∪ I⟩ has atleast one p-Sylow S-mixed neutrosophic N-substructure then we call ⟨M ∪ I⟩ to be a weak Sylow S-mixed neutrosophic N-structure. If ⟨M ∪ I⟩ has no p-Sylow S-mixed neutrosophic sub N-structure then we call ⟨M ∪ I⟩ to be a Sylow free S-mixed neutrosophic N-structure.*

This class has infinitely many S-mixed neutrosophic N-structures as all S-mixed neutrosophic N-structure of finite order n, n a prime are Sylow free S-mixed neutrosophic N-structures. These concepts can be defined in case of weak S-mixed neutrosophic N-structure which has weak S-mixed neutrosophic sub N-structure and for S-mixed neutrosophic deficit sub N-structures.

Also for all these the S-mixed neutrosophic N-structure of finite order n, n a prime is such that they will be Sylow free Lagrange free and Cauchy free. However we give an example.

***Example 6.1.3:*** Let {⟨M ∪ I⟩ = $M_1 \cup M_2 \cup M_3 \cup M_4 \cup M_5 \cup M_6$, $*_1$, $*_2$, $*_3$, $*_4$, $*_5$, $*_6$} be a S-mixed neutrosophic N-structure of finite order, where $M_1$ = {⟨$L_5(3)$ ∪ I⟩}, $M_2$ = {⟨$Z_6$ ∪ I⟩, neutrosophic semigroup under multiplication modulo 6}. $M_3$ = {1, 2, 3, 4, I, 2I, 3I, 4I} a neutrosophic group under multiplication modulo 5. $M_4$ = {0, 1, 2, 3, I, 2I, 3I with a $*_4$ b = (2a + b) mod 4} a neutrosophic groupoid, $M_5$ = {g | $g^4$ = 1} and $M_6 = Z_8$, semigroup under multiplication modulo 8. o(⟨M ∪ I⟩) = 50. Take P = {$P_1 \cup P_2 \cup P_3 \cup P_4 \cup P_5 \cup P_6$, $*_1$, $*_2$, $*_3$, $*_4$, $*_5$, $*_6$} where $P_1$ = {e, eI, 2, 2I}, $P_2$ = {0, 2, 2I, 4, 4I}, $P_3$ = {4, 4I, I, 1}, $P_4$ = {2, 2I, 0}, $P_5$ = {1, $g^2$ | $g^4$ = 1} and $P_6$ = {0, 2, 4, 6}. P is a S-mixed neutrosophic sub 6 structure.

Take T = {$T_1 \cup T_2 \cup T_3 \cup T_4 \cup T_5$, ∪ $T_6$, $*_1$, $*_2$, $*_3$, $*_4$, $*_5$, $*_6$} where $T_1$ = {e, 2}, $T_2$ = {1, 3, 3I, I}, $T_3$ = {4, 4I, I, 1}, $T_4$ =

$\{0, 2\}$, $T_5 = \{1, g^2\}$ and $T_6 = \{0, 4\}$. T is a weak S-mixed neutrosophic sub N-structure.

Now we give an example of a S-mixed neutrosophic deficit sub N-structure. Let $W = \{W_1 \cup W_2 \cup W_3, *_1, *_2, *_3\}$ where $W_1 = \{e, eI, 2, 2I\}$, $W_2 = \{I, 3, 3I, 1\}$, $W_3 = \{0, 2, 4, 6\}$ semigroup under multiplication modulo 8. Clearly W is a S-mixed neutrosophic deficit sub N-structure.

Thus we have seen all the 3 types of sub N-structures of a S-mixed neutrosophic deficit N-structure. Like wise we can see this S-mixed neutrosophic N-structure has no 2-Sylow S-mixed neutrosophic sub 6-structure but has a 5-Sylow S-mixed neutrosophic sub 6-structure given by $\langle X \cup I \rangle = \{X_1 \cup X_2 \cup X_3 \cup X_4 \cup X_5, \cup X_6, *_1, *_2, *_3, *_4, *_5, *_6\}$ where $X_1 = \{e, eI, 2, 2I\}$, $X_2 = \{2, 4, 2I, 4I, I, 1\}$, $X_3 = \{1, 4, 4I, I\}$, $X_4 = \{1, 0, 2, 2I, I\}$, $X_5 = \{1, g^2\}$ and $X_6 = \{0, 2, 4, 6\}$. $o(\langle X \cup I \rangle) = 25$. 25 / 50 and $5^3 \not| \; 50$. $\langle X \cup I \rangle$ is a 5-Sylow S-mixed neutrosophic sub 6-structure. So $\langle M \cup I \rangle$ is a weak Sylow S-mixed neutrosophic 6-structure.

This has both Cauchy and Cauchy neutrosophic element but not all elements are Cauchy or Cauchy neutrosophic. We propose several interesting problem in the last chapter of this book for any innovative and interested reader.

**Chapter Seven**

# PROBLEMS

This chapter suggests some problems for the interested reader. Most of the problems are simple exercises. Topics which are not elaborately covered will find its place in this chapter as suggested problems.

Since the notion of S-neutrosophic groupoids have not been dealt deeply here we give 25 problems on S-neutrosophic groupoids, S-neutrosophic bigroupoids and S-neutrosophic N-groupoids.

1. Suppose (S, o) is a finite Smarandache neutrosophic semigroup. P and T be two Smarandache neutrosophic subsemigroups of S which are conjugate i.e., we have $P_1$ in P and $T_1$ in T such that $(P_1, o)$ and $(T_1, o)$ are groups which re conjugate. Let $(x_1, \ldots, x_n) \subset S$ be such that $T_1x_i = x_iP_1$ (or $P_1x_i = x_iT_1$) for i = 1, 2, 3, …, n.

   Does $V = \{x_1, \ldots, x_n\} \cup \{0\}$ have any algebraic structure? When is it a Smarandache neutrosophic subsemigroup? Is it possible to classify those Smarandache neutrosophic semigroup whose Smarandache conjugacy set happens to a Smarandache neutrosophic sub semigroup. When will the Smarandache conjugacy set be a group?

2. Give examples of S-neutrosophic of finite order loops which are:

   i. S-Lagrange neutrosophic loop and not Lagrange neutrosophic loops.
   ii. Lagrange neutrosophic loops and not S-Lagrange neutrosophic loops.
   iii. Weak Lagrange neutrosophic loops and not weak S-Lagrange neutrosophic loops and vice versa.

3. Does the new class of S-neutrosophic loops $\langle L_n \cup I \rangle$, n not a prime always have a Sylow S-neutrosophic loop or never have a Sylow S-neutrosophic loop? If we have $\langle L_n \cup I \rangle$ for a non prime n to be a Sylow S-neutrosophic loop does it imply $2(n + 1) = 2^t$ or obtain condition on n.

4. Find conditions on the order of the neutrosophic loop so that

   i. a Cauchy element is a S-Cauchy element.
   ii. a S-Cauchy element is not a Cauchy element.

   *Note:* Work this problem at least in the case of the new class of S- neutrosophic loops.

5. Give examples of S-neutrosophic loop which are not S-Cauchy or S- Cauchy neutrosophic.

   *Note:* One has to find loops from other classes for the class of S-neutrosophic loops $\langle L_n \cup I \rangle$ happen to be S-Cauchy and S-Cauchy neutrosophic.

6. Define S-neutrosophic Bruck loop and give examples. Does the class of neutrosophic loops $\langle L_n(m) \cup I \rangle$ have S-neutrosophic Bruck loop?

7. Define S-neutrosophic Bol loop and give examples! Does the class of new S-neutrosophic loop $\langle L_n(m) \cup I \rangle$ have any S-neutrosophic Bol loop?

8. Find conditions on m and n other than mentioned in the book so that $\langle L_n(m) \cup I \rangle$ is a S-neutrosophic loop.

9. Give examples of quasi Cauchy S-neutrosophic biloops.

10. Find conditions on the S-neutrosophic N-groups or on its order so that a S-Cauchy neutrosophic element is S-Cauchy element of $\langle G \cup I \rangle$.

11. Find conditions so that S Cauchy elements are Cauchy elements of $\langle G \cup I \rangle$.

12. Give examples of weak Cauchy N-loops.

13. Does their exist relation between Cauchy S-neutrosophic biloops and Quasi Cauchy S-neutrosophic biloops?

14. Give an example of a S-neutrosophic 5-group which is of composite order which is S-Lagrange neutrosophic 5-group.

15. Give an example of a Lagrange S-neutrosophic S-group.

16. Can a S-neutrosophic group be both Super Sylow S-neutrosophic group as well as Lagrange S-neutrosophic group? Justify your claim.

17. Give an example of a S-neutrosophic bigroup which is weakly S-Lagrange.

18. Does there exist any relation between a S-pseudo Lagrange neutrosophic N-group and a S-Lagrange neutrosophic N-group?

19. Give an example of a S-neutrosophic N-group which is both S-Lagrange and S-pseudo Lagrange.

20. Let $\langle Z_{20} \cup I\rangle$ be a S-neutrosophic semigroup under multiplication modulo 20. Find all its

   i. S-neutrosophic subsemigroups
   ii. S-neutrosophic ideals
   iii. Is $\langle Z_{10} \cup I\rangle$ a S-Lagrange neutrosophic semigroup?

21. Does $\langle Z_{20} \cup I\rangle$ given in problem 20 have Cauchy neutrosophic elements?

22. Can S-Cayley theorem for S-neutrosophic N-groups be proved?

23. Let $\{S(S) = S_1 \cup S_2 \cup S_3, *_1, *_2, *_3\}$ be a S-neutrosophic 3-semigroup where $S_1 = \langle Z_{10} \cup I\rangle$, S-semigroup under multiplication modulo 10, $S_2 = S(3)$ and $S_3 = Z_4$ under multiplication modulo 4. Is S(S) a S-Lagrange neutrosophic 3-group? Is S(S) a S-Sylow neutrosophic 3-group?

24. Does the problem 23 given above have S-neutrosophic quasi maximal 3-ideal? Or S-neutrosophic minimal 3-ideal or S-neutrosophic maximal 3-ideal?

25. Give an example of a S-neutrosophic N-semigroup which has no S-neutrosophic sub N-semigroup? ($N \geq 3$).

26. Illustrate by an example the concept of S-Lagrange neutrosophic N-semigroup.

27. Can any condition be derived so that the S-neutrosophic N-semigroup is always S-free Lagrange neutrosophic N-semigroup?

28. Give an example of S-neutrosophic N-semigroup which is a special Cauchy S-neutrosophic N-semigroup.

29. Illustrate by an example S-neutrosophic 5-semigroup which has non-trivial S-conjugate neutrosophic sub-5-semigroups.

30. Is $\{\langle L \cup I\rangle = L_1 \cup L_2, *_1, *_2\}$ where $L_1 = \langle L_5(3) \cup I\rangle$ and $L_2 = \langle L_{15}(2) \cup I\rangle$ a S-neutrosophic biloop?

31. Find all S-neutrosophic subbiloops in a S-neutrosophic biloop given in problem 30.

32. Is the biloop given in problem 30 S-neutrosophic right alternative?

33. Illustrate by an example a S-neutrosophic Bruck biloop.

34. Define

    i. S-Sylow neutrosophic N-loop
    ii. S-Lagrange N-loop
    iii. S-weak Lagrange N-loop

    and illustrate them by examples.

35. Define S-(N – t) deficit neutrosophic sub N-loop and find all S-(N – t) deficit neutrosophic sub N-loop for N = 7. How many such sub N-loops exist for N = 7 ($1 \leq t < 7$)?

36. Give an example of S-WIP-neutrosophic-N-loop?

37. Define for a S-neutrosophic N-loop the notion of Cauchy element and Cauchy neutrosophic element.

38. Define a S-neutrosophic homomorphism between any two S-neutrosophic N-loops. Illustrate it by an example.

39. Define S-alternative neutrosophic groupoid and illustrate it by an example.

40. Is the S-neutrosophic groupoid $(\langle Z_{16} \cup I\rangle, *)$ where $a * b = 12a + 4b \pmod{16}$ a S-Lagrange neutrosophic groupoid?

41. Find all S-neutrosophic subgroupoids and S-subgroupoids of the S-neutrosophic groupoid given in problem 40.

42. Give an example of a S-Lagrange neutrosophic bigroupoid.

43. Define weakly S-Lagrange neutrosophic N-groupoid and illustrate it by an example.

44. Does there exist a S-neutrosophic N-groupoid which is both S-Moufang and S-WIP?

45. Give an example of S-Moufang neutrosophic 5-groupoid.

46. Define the following:

    i. S-Bol neutrosophic N-groupoids.
    ii. S-P-neutrosophic N-groupoids.
    iii. S-WIP-neutrosophic N-groupoids

    and illustrate them by examples.

47. Define for S-neutrosophic N-quasi loop, S-Moufang neutrosophic N-quasi loop, S-Bol neutrosophic N-quasi loop, S-alternative neutrosophic N-quasi loop and illustrate them with examples.

48. Define

    i. S-Lagrange neutrosophic N-groupoid.
    ii. S-weak Lagrange neutrosophic N-groupoid.
    iii. S-Sylow neutrosophic N-groupoid.
    iv. S-Sylow free neutrosophic N-groupoid.
    v. S-weak Sylow neutrosophic N-groupoid.

49. Illustrate the above 5 definitions with examples and find some interesting results using these definitions.

50. Define S-semi-conjugate sub N-groupoids and illustrate them with examples.

51. Define the 5 definitions given in problem 48 in case of S-neutrosophic bigroupoids and S-neutrosophic groupoids and illustrate these concepts with examples.

52. Define Smarandache seminormal sub N-groupoid and explain it by an example for N = 5.

53. Give an example of a S-neutrosophic bigroupoid of composite order which is simple.

54. Does the S-strong neutrosophic bigroupoid given by S(G) = $\{G_1 \cup G_2, *_1, *_2\}$ where $G_1 = \{\langle Z_{12} \cup I / a *_1 b = 8a + 4b, \text{mod } 12\rangle$ and $G_2 = \{\langle Z_6 \cup I\rangle \times \langle Z_8 \cup I\rangle$ / if (a, b), (c, d) is in $\langle Z_6 \cup I\rangle \times \langle Z_8 \cup I\rangle$; (a, b) $*_2$ (c, d) = (2a + 4c(mod 6), (6b + 2d) (mod 8)} have

    i. S-Normal neutrosophic subbigroupoids.
    ii. S-seminormal neutrosophic subbigroupoids.
    iii. S-Lagrange neutrosophic subbigroupoids.
    iv. p-Sylow S-neutrosophic subbigroupoids.

55. Is the S-neutrosophic bigroupoid given in problem 54, S-Moufang or S-Bol or S-P-groupoid or S-right alternative or S-left alternative?

56. Find all S-neutrosophic biideals of the S-neutrosophic bigroupoid given in problem 54.

57. Is the S-neutrosophic S-bigroupoid given in problem 54, have Cauchy elements and Cauchy neutrosophic elements.

58. Give an example of S-weak neutrosophic N-groupoid.

59. Is S(G) = $\{G_1 \cup G_2 \cup G_3 \cup G_4, *_1, *_2, *_3, *_4\}$ where $G_1 = \{\langle Z_{12} \cup I\rangle / a *_1 b = 8a + 4b \pmod{12}$ for a, b $\in \langle Z_{12} \cup I\rangle\}$, $G_2 = \{\langle Z_{10} \cup I\rangle$ neutrosophic semigroup under multiplication modulo 10}, $G_3$ = S(3) and $G_4 = \{Z_6 \times Z_4$ / if (a, b), (c, d) $\in Z_6 \times Z_4$, (a, b) $*_4$ (c, d) = (a + 5c)(mod 6), (3b + d)(mod 4)} a S-weak neutrosophic 4-groupoid?

60. Is the S-weak neutrosophic 4-groupoid given in problem 59,

    i. S-Sylow.
    ii. S-Lagrange.
    iii. S-weak Lagrange.
    iv. or S-weak Sylow?

    Justify your claim.

61. Define a S-neutrosophic homomorphism between S-weak neutrosophic N-groupoids?

62. Does the S-neutrosophic N-groupoid given in problem 60 have any S-neutrosophic conjugate sub N-groupoids?

63. Given an example of a S-neutrosophic quasi semigroup of order 48, N = 4.

64. What is the main difference between S-neutrosophic N-quasi loop and S-neutrosophic N-quasi semigroup.

65. Give an example of a S-Lagrange neutrosophic N-quasi loop (N = 7).

66. Illustrate by an example the S-Sylow neutrosophic 5-quasi group.

67. Is S-neutrosophic N-quasi loops of finite prime order p, S-Lagrange free?

68. Can a S-neutrosophic N-groupoid of prime order have non-trivial S-neutrosophic sub N-groupoids. Illustrate your claim by an example.

# INDEX

### A

### B

### C

**D**

**F**

**G**

**I**

**K**

**L**

**M**

**N**

**O**

**P**

**Q**

**R**

**S**

**T**

**W**

# ABOUT THE AUTHOR

**Dr.W.B.Vasantha Kandasamy** is an Associate Professor in the Department of Mathematics, Indian Institute of Technology Madras, Chennai, where she lives with her husband Dr.K.Kandasamy and daughters Meena and Kama. Her current interests include Smarandache algebraic structures, fuzzy theory, coding/ communication theory. In the past decade she has guided 11 Ph.D. scholars in the different fields of non-associative algebras, algebraic coding theory, transportation theory, fuzzy groups, and applications of fuzzy theory of the problems faced in chemical industries and cement industries. Currently, four Ph.D. scholars are working under her guidance.

She has to her credit 612 research papers of which 209 are individually authored. Apart from this, she and her students have presented around 329 papers in national and international conferences. She teaches both undergraduate and post-graduate students and has guided over 45 M.Sc. and M.Tech. projects. She has worked in collaboration projects with the Indian Space Research Organization and with the Tamil Nadu State AIDS Control Society. This is her 26th book.

She can be contacted at vasantha@iitm.ac.in
You can visit her work on web at: http://mat.iitm.ac.in/~wbv